\documentclass[10pt,reqno]{amsart}
\usepackage{geometry} 
\usepackage{amsfonts}
\usepackage{amsmath}
\usepackage{amssymb}
\usepackage{color}
\usepackage{tikz}
\usepackage{tikz-cd}
\usepackage{enumitem}
\usepackage{scalerel}
\usepackage{hyperref}
\hypersetup{
  colorlinks   = true, 
  urlcolor     = blue, 
  linkcolor    = blue, 
  citecolor   = red 
}

\let\oldtocsection=\tocsection
\let\oldtocsubsection=\tocsubsection
\let\oldtocsubsubsection=\tocsubsubsection
\renewcommand{\tocsection}[2]{\hspace{0em}\oldtocsection{#1}{#2}}
\renewcommand{\tocsubsection}[2]{\hspace{1em}\oldtocsubsection{#1}{#2}}
\renewcommand{\tocsubsubsection}[2]{\hspace{2em}\oldtocsubsubsection{#1}{#2}}

\usepackage{float}

\numberwithin{equation}{section}
\allowdisplaybreaks

\newtheorem{theorem}{Theorem}[section]
\newtheorem{definition}[theorem]{Definition}
\newtheorem{corollary}[theorem]{Corollary}
\newtheorem{remark}[theorem]{Remark}
\newtheorem{lemma}[theorem]{Lemma}
\newtheorem{proposition}[theorem]{Proposition}

\newtheorem*{method I}{Method I}
\newtheorem*{method II}{Method II}
\newtheorem*{step i}{Step i}
\newtheorem*{step ii}{Step ii}
\newtheorem*{assum}{\it{Assumption}}
\newtheorem*{main result sketch}{\bf{Sketch of main result}}

\newcommand{\scslash}{\scaleto{\slash}{1.6ex}}

\newcommand{\leve}{\big\vert}
\newcommand{\rive}{\big\vert}
\newcommand{\leVe}{\big\Vert}
\newcommand{\riVe}{\big\Vert}

\newcommand{\circg}{\overset{\raisebox{-0.6ex}[0.6ex][0ex]{\scaleto{\circ}{0.6ex}}}{g}\makebox[0ex]{}}    
\newcommand{\circGamma}{\overset{\raisebox{0ex}[0.6ex][0ex]{\scaleto{\circ}{0.6ex}}}{\raisebox{0ex}[1.2ex][0ex]{$\Gamma$}} \makebox[0ex]{}}    
\newcommand{\circnabla}{\overset{\raisebox{0ex}[0.6ex][0ex]{\scaleto{\circ}{0.6ex}}}{\raisebox{0ex}[1.3ex][0ex]{$\nabla$}} \makebox[0ex]{}}    
\newcommand{\circDelta}{\overset{\raisebox{0ex}[0.6ex][0ex]{\scaleto{\circ}{0.6ex}}}{\raisebox{0ex}[1.3ex][0ex]{$\Delta$}} \makebox[0ex]{}}     

\newcommand{\us}{\underline{s}}   
\newcommand{\uC}{\underline{C}}  

\newcommand{\slashg}{g \mkern-8mu \scslash \mkern+1mu \makebox[0ex]{}}  
\newcommand{\dvol}{\mathrm{dvol}}               
\newcommand{\tr}{\mathrm{tr}}              
\newcommand{\slashGamma}{\Gamma \mkern-9mu \raisebox{+0.3ex}[0ex][0ex]{$\scslash$} \mkern+3mu}     

\newcommand{\uL}{\underline{L}}

\newcommand{\ueta}{\underline{\eta}}  
\newcommand{\hatchi}{\hat{\chi}}          
\newcommand{\uchi}{\underline{\chi}}   
\newcommand{\hatuchi}{\hat{\uchi}}      
\newcommand{\uomega}{\underline{\omega}}   

\newcommand{\udelta}{\underline{\delta}}

\renewcommand{\d}{\mathrm{d}}        
\newcommand{\slashd}{\d\mkern-8.5mu \raisebox{+0.3ex}[0ex][0ex]{$\scslash$} \mkern+3mu}      
\newcommand{\slashnabla}{\nabla \mkern-11mu \raisebox{+0.3ex}[0ex][0ex]{$\scslash$}\mkern+5mu \makebox[0ex]{}}
\renewcommand{\div}{\mathrm{div}}
\newcommand{\slashDelta}{\Delta \mkern-10mu \raisebox{+0.3ex}[0ex][0ex]{$\scslash$} \mkern+4mu \makebox[0ex]{}}        
\newcommand{\lie}{\mathcal{L}}           

\newcommand{\sym}{\mathrm{sym}}   

\newcommand{\uf}{\underline{f}}    
\newcommand{\ue}{\underline{e}}       
\newcommand{\uvarepsilon}{\underline{\varepsilon}}   
\newcommand{\Xlt}{\mkern+1mu {}^{t} \mkern-4mu X}  
\newcommand{\flt}{\mkern+1mu{}^{t} \mkern-4mu f}  
\newcommand{\lo}[1]{{l.}\{ #1\}} 
\newcommand{\hi}[1]{{h.}\{ #1\} } 
\newcommand{\dd}[1]{\mathfrak{d} \{ #1 \}} 
\newcommand{\Flt}{{}^{t} \mkern-3mu F} 
\newcommand{\ufrakd}{\underline{\mathfrak{d}}} 
\newcommand{\frakd}{\mathfrak{d}} 
\newcommand{\bdd}[1]{\delta \{ #1 \}} 
\newcommand{\er}[1]{\mathbf{e} \{#1\}} 
\newcommand{\bft}{\mathbf{t}} 

\newcommand{\ufl}[1]{\mkern+1mu{}^{#1}\mkern-4mu \uf}
\newcommand{\uflt}{\ufl{t}}
\newcommand{\relt}{\mkern+1mu {}^t \mkern-1mu \mathrm{re}}
\newcommand{\rel}[1]{\mkern+1mu {}^{#1} \mkern-1mu \mathrm{re}}
\newcommand{\Xl}[1]{\mkern+1mu {}^{#1} \mkern-4mu X}
\newcommand{\Fl}[1]{\mkern+1mu {}^{#1} \mkern-4mu F}

\newcommand{\ufls}{\ufl{s}}
\newcommand{\ucalH}{\underline{\mathcal{H}}}
\newcommand{\dL}{\dot{L}}
\newcommand{\duL}{\dot{\uL}}
\newcommand{\dpartial}{\dot{\partial}}
\newcommand{\dslashg}{\dot{g} \mkern-8mu \scslash \mkern+1mu \makebox[0ex]{}}
\newcommand{\dslashnabla}{\dot{\nabla}  \mkern-11mu \raisebox{+0.3ex}[0ex][0ex]{$\scslash$}\mkern+5mu \makebox[0ex]{}}
\newcommand{\dslashGamma}{\dot{\Gamma}  \mkern-9mu \raisebox{+0.3ex}[0ex][0ex]{$\scslash$} \mkern+3mu \makebox[0ex]{}}
\newcommand{\dtr}{\dot{\tr}}
\newcommand{\dchi}{\dot{\chi}}
\newcommand{\duchi}{\dot{\uchi}}
\newcommand{\deta}{\dot{\eta}}
\newcommand{\duomega}{\dot{\uomega}}
\newcommand{\dslashd}{\dot{\d}\mkern-8.5mu \raisebox{+0.2ex}[0ex][0ex]{$\scslash$} \mkern+3mu}
\newcommand{\db}{\dot{b}}

\newcommand{\ddpartial}{\ddot{\partial}}
\newcommand{\ddL}{\ddot{L}}
\newcommand{\dduL}{\ddot{\uL}}
\newcommand{\dB}{\dot{B}}
\newcommand{\ddslashg}{\ddot{g} \mkern-8mu \scslash \mkern+1mu \makebox[0ex]{} }
\newcommand{\dvarepsilon}{\dot{\varepsilon}}
\newcommand{\ddslashd}{\ddot{\slashd}}
\newcommand{\circdot}{\makebox[2ex]{$\circ\mkern-7mu\cdot$}}
\newcommand{\ddtr}{\ddot{\tr}}
\newcommand{\ddchi}{\ddot{\chi}}
\newcommand{\dduchi}{\ddot{\uchi}}

\newcommand{\dR}{\dot{R}}
\newcommand{\dcirc}{\overset{\raisebox{-0.1ex}[0ex][0ex]{$\cdot$}}{\raisebox{0ex}[0ex][0ex]{\scaleto{\circ}{0.6ex}}}}
\newcommand{\dcircnabla}{\raisebox{-0.1ex}[2.5ex][0ex]{$\overset{\dcirc}{\raisebox{0ex}[1.21ex][0ex]{$\nabla$}}$} \makebox[0ex]{}}
\newcommand{\dcircDelta}{\raisebox{-0.1ex}[2.5ex][0ex]{$\overset{\dcirc}{\raisebox{0ex}[1.215ex][0ex]{$\Delta$}}$} \makebox[0ex]{}}     
\newcommand{\ddcirc}{\overset{\raisebox{-0.1ex}[0ex][0ex]{$\cdot \mkern-1mu \cdot$}}{\raisebox{0ex}[0ex][0ex]{\scaleto{\circ}{0.6ex}}}}
\newcommand{\ddcircnabla}{\raisebox{-0.1ex}[2.5ex][0ex]{$\overset{\ddcirc}{\raisebox{0ex}[1.21ex][0ex]{$\nabla$}}$} \makebox[0ex]{}}
\newcommand{\ddcircDelta}{\raisebox{-0.1ex}[2.5ex][0ex]{$\overset{\ddcirc}{\raisebox{0ex}[1.215ex][0ex]{$\Delta$}}$} \makebox[0ex]{}}

\newcommand{\circtriangle}{\overset{\raisebox{0ex}[0.6ex][0ex]{\scaleto{\circ}{0.6ex}}}{\raisebox{0ex}[1.3ex][0ex]{$\triangle$}} \makebox[0ex]{} \makebox[0ex]{}}

\newcommand{\Sigmal}[1]{\mkern+1mu {}^{#1} \mkern-1mu \Sigma}
\newcommand{\Sl}[1]{\mkern+1mu {}^{#1} \mkern-3mu S}
\newcommand{\fl}[1]{\mkern+1mu {}^{#1} \mkern-4mu f}

\newcommand{\os}{\overline{s}}
\newcommand{\ous}{\overline{\us}}

\newcommand{\ddcircdiv}{\raisebox{-0.1ex}[2.5ex][0ex]{$\overset{\ddcirc}{\raisebox{0ex}[1.21ex][0ex]{$\div$}}$} \makebox[0ex]{}}

\newcommand{\uh}{\underline{h}}

\newcommand{\ud}{\underline{d}}
\newcommand{\ubfd}{\underline{\bf d}}
\newcommand{\uslashd}{\underline{d} \mkern-8mu \raisebox{0.3ex}[0ex][0ex]{$\scslash$} \mkern+2mu}
\newcommand{\uslashbfd}{\underline{\bf d} \mkern-9mu \raisebox{0.3ex}[0ex][0ex]{$\scslash$} \mkern+3mu }

\newcommand{\uhl}[1]{{}^{#1} \mkern-2mu \underline{h}}
\newcommand{\uhlt}{{}^{t} \mkern-2mu \underline{h}}

\newcommand{\hir}[2]{{h_{\boldsymbol{\cdot}}^{#1}}\{ #2\} } 

\newcommand{\hl}[1]{{}^{#1} \mkern-2mu h}
\newcommand{\uls}[2]{{}^{#1} \mkern-2mu {#2}} 
\newcommand{\ubfe}{\underline{\mathbf{e}}}

\newcommand{\upartial}{\underline{\partial}}
\newcommand{\bfs}{\mathbf{s}}
\newcommand{\bfe}{\mathbf{e}}

\title[Marginally trapped surfaces in a perturbed Schwarzschild spacetime]{Marginally trapped surfaces in a perturbed Schwarzschild spacetime}
\author{Pengyu Le}
\date{} 

\begin{document}

\begin{abstract}
The concept of a marginally trapped surface is important in the theory of general relativity. In the Schwarzschild black hole spacetime, its event horizon is foliated by marginally trapped surfaces. In a more general black hole spacetime, the concept of a marginally trapped surface is closely related to various sorts of horizon, for example, the apparent horizon, the trapping boundary, the isolated horizon and the dynamical horizon. In this paper, we study the set of marginally trapped surfaces in a perturbed Schwarzschild spacetime. We show that for every incoming null hypersurface which is nearly spherically symmetric, there exists a unique embedded marginally trapped surface. In order to prove this result, we develop a general method to study the geometry of spacelike surfaces in a double null coordinate system, which can be applied to study other problems for spacelike surfaces in a Lorentzian manifold.
\end{abstract}
\maketitle
\tableofcontents

\section{Introduction}\label{sec 1}
\noindent
The Schwarzschild black hole spacetime found in 1916 \cite{Sc} is the static spherically symmetric vacuum solution of the Einstein equations, soon after Einstein's discovery of his field equations \cite{E1}, \cite{E2}. Its metric reads as follows:
\begin{align}
\nonumber
\d s^2 = - \left( 1-\frac{2m}{r} \right) \d t^2 + \left( 1- \frac{2m}{r} \right)^{-1} \d r^2 + r^2 \left( \d \theta^2 + \sin^2 \theta \d \phi^2 \right).
\end{align}
This metric depends on the parameter $m$ whose physical meaning is the mass of the spacetime. When $m=0$, it becomes the flat Minkowski metric. At first sight, $r=0$ and $r=2m$ look like values for which the metric is singular. Only $r=0$ is a true singularity but $r=2m$ is a coordinate singularity, which can be removed by coordinate transformations. See the classical works \cite{Ed}, \cite{Le}, \cite{Fi}.

Synge (1950 \cite{Sy}), Kruskal (1960 \cite{Kr}) and Szekeres (1960 \cite{Sz}) provided coordinate systems that cover the maximal analytic extension of the Schwarzschild metric. In the Kruskal-Szekeres coordinates $\{u, v ,\theta,\phi \}$, the metric takes the form
\begin{align}
\nonumber
\begin{aligned}
&
\d s^2= -\frac{16m^2}{r} \exp{ \left(\frac{-r}{2m}\right) } \d u \d v + r^2\left(\d \theta^2+\sin^2 \theta \d \phi^2 \right),
\\
&
uv = - (r- 2m)\exp{ \frac{r}{2m} }.
\end{aligned}
\end{align}
The Schwarzschild black hole in the above coordinates can be visualised by figure \ref{fig 1}.
\begin{figure}[H]
\begin{center}
\begin{tikzpicture}
\draw (-3,-3) -- (3,3) node[anchor=north west]{$u=0,r=2m$};
\draw(3,-3) -- (-3,3) node[anchor=north east]{$v=0,r=2m$};
\draw[dotted,thick,domain=-2.5:2.5]  plot(\x, {sqrt(3+\x*\x)}) node[left]{$uv=2m,r=0$};
\draw[dotted,thick,domain=-2.5:2.5]  plot(\x, {-sqrt(3+\x*\x)}) node[left]{$uv=2m,r=0$};
\draw[domain=0:5]  plot(\x, {-0.2*\x}) node[below]{$t=\text{const}$};
\draw[domain=-2:2,variable=\y]  plot({sqrt(5+\y*\y)},\y) node[right]{$r=\text{const}>2m$};
\draw (3,3)--(3.5,3.5) node[right] {event horizon};
\node at (0,1.2) {black hole};
\node at (3.7,0) {exterior region};
\end{tikzpicture}
\end{center}
\caption{The maximal analytic extension of the Schwarzschild spacetime.}
\label{fig 1}
\end{figure}

In 1965 \cite{Pe}, Penrose introduced the concept of a \emph{closed trapped surface}, which is is a closed spacelike surface where the area element decreases pointwise for any infinitesimal displacement along the future null normal direction. Based on this concept, he proved his famous incompleteness theorem.

A concept related to the one of a trapped surface is the concept of a \emph{marginally trapped surface} defined as follows.

\begin{definition}
A spacelike surface $\Sigma$ is called a marginal surface, if one of its future null expansions vanishes identically. 

Furthermore, a spacelike surface $\Sigma$ is called marginally trapped, if one of its future null expansions vanishes identically, and the other future null expansion is non-positive.
\end{definition}

The concept of a marginally trapped surface is closely related to the horizon of a black hole in general relativity. For example, the event horizon of either a Schwarzschild or a Kerr black hole is foliated by marginally trapped surfaces. However it is not the case that the event horizon is foliated by marginally trapped surfaces in a general black hole spacetime. As pointed out in \cite{HE}, the event horizon is a global concept, which depends on the whole future behaviour of the spacetime, thus it is useful to define some different sort of horizon which depends only on the geometry of a spacelike slice of the spacetime. \cite{HE} introduced the concept of a trapped region in a spacelike hypersurface $H$, which is the set of all points, through which a trapped surface in $H$ passes. Moreover, \cite{HE} defined the concept of an apparent horizon as a connected component of the outer boundary of the trapped region, and shows that an apparent horizon shall be marginally trapped.

Besides the concept of an apparent horizon introduced in \cite{HE}, there are various other useful concepts of horizon, for examples, the trapping boundary \cite{Ha}, the isolated horizon \cite{ABF}, the dynamical horizon \cite{AK}. The later two sorts of horizon are examples of the more general concept of a marginally trapped tube introduced in \cite{AG}, which has the topology of $\mathbb{S}^2 \times \mathbb{R}$ and is foliated by marginally trapped spheres. Moreover, \cite{AG} proved the uniqueness of the foliation by marginally trapped surfaces of a dynamical horizon. The work \cite{AMS} proved the existence of a marginally trapped tube under a stability condition of a marginally trapped surface. In \cite{L3}, a more general concept of a marginal tube was introduced, which is foliated by marginal surfaces and where no restriction of the topology is required. It was showed that if every embedded spacelike surface of a marginal tube is a marginal surface, then the marginal tube must be null.

Since the marginally trapped surface is closely related to many useful concepts of horizon in a black hole spacetime, it is natural to ask what is the structure of the set of marginally trapped surfaces in a perturbation of the stationary black hole. The understanding of this set will be useful for the study of the geometry of a perturbed stationary black hole.

The theme of this paper is the study of marginally trapped surfaces in a perturbation of the simplest stationary black hole, the Schwarzschild black hole. We shall consider a perturbation of the Schwarzschild black hole near its event horizon. We employ the double null coordinate system of the Schwarzschild black hole to quantitatively described the perturbation, by comparing the metric components and structure coefficients relative to the double null coordinate system. The precise description of the perturbation is given in definition \ref{def 2.3} in section \ref{sec 2}. Here we just emphasis two points in the perturbation: first that the perturbation is not necessary being vacuum, second that none of the coordinate surface in the double null coordinate system is required to be marginally trapped.

The main result of this paper can be sketchily phrased as follows.
\begin{main result sketch}
Let $(M,g)$ be a Lorentzian manifold which is a perturbation of the Schwarzschild spacetime near the event horizon. In every incoming null hypersurface which is nearly spherically symmetric, there exists a unique marginally trapped surface near the Schwarzschild event horizon.
\end{main result sketch}

The precise version of the main result is given in theorem \ref{thm 10.9} in section \ref{sec 10}. It can be explained geometrically by figure \ref{fig 2}. $\{\us, s\}$ is the double null coordinate system. $C_{s=0}$ is the level set $\{s=0\}$, which is a null hypersurface surface. It is the Schwarzschild event horizon. In every nearly spherically symmetric incoming null hypersurface $\ucalH$ in a perturbed Schwarzschild spacetime, there exists a unique embedded marginally trapped surface $\Sigma$ near $C_{s=0}$.
\begin{figure}[h]
\begin{center}
\begin{tikzpicture}
\draw[dashed] (-1,0) to [out=70,in=-110] (-0.7,0.9);
\draw (-1,0) to [out=-110,in=70] (-1.85+0.5,-2.4+1.4);
\draw[dashed] (1,0) to [out=110,in=-70] (0.7,0.9);
\draw[->] (1,0) to [out=-70,in=110] (1.85-0.5,-2.4+1.4) node[right] {\small $s$}; 
\draw[dashed] (-1,0) to [out=-45,in=135] (-0.5,-0.5);
\draw (-1,0) to [out=135, in= -45] (-2,1);
\draw[dashed] (1,0) to [out=-135,in=45] (0.5,-0.5);
\draw[->] (1,0) to [out=45, in= -135] (1.9,0.9) node[right] {\tiny $C_{s=0}$}
to [out=45,in=-135] (2.3,1.3) node[right] {\small $\us$}; 
\draw[dashed] (-1.5,0.5) to [out=135,in=45] (1.5,0.5);
\draw (1.5,0.5) to [out=-135,in=0] (0,0) to [out=180,in=-45] (-1.5,0.5); 
\draw[dashed] (-1.4,1.3) to [out=-110,in=70] (-1.6,0.7);
\draw (-1.6,0.7) to [out=-110,in=70] (-2.75+0.5,-2.4+1.4);
\draw[dashed] (1.4,1.3) to [out=-70,in=110] (1.6,0.7);
\draw[->] (1.6,0.7) to [out=-70,in=110] (2.75-0.5,-2.4+1.4) node[right]{\small $s$}; 
\node[right] at (2.4-0.4,-1.5+1) {\scriptsize $\ucalH$}; 
\draw[dashed] (-2.3+0.4,-2.2+2.1) to [out=70,in=180] (-0.5,-2.5+2) node[below] {\scriptsize $\Sigma$} to [out=0,in=110] (2.3-0.4,-2.2+2.1);
\draw (2.3-0.4,-2.2+2.1) to [out=-70,in=0] (1,-2+1.9) to [out=180,in=-110] (-2.3+0.4,-2.2+2.1); 
\end{tikzpicture}
\end{center}
\caption{The marginally trapped surface $\Sigma$ in an incoming null hypersurface $\ucalH$.}
\label{fig 2}
\end{figure}

In the following, we give an overview of the strategy to prove the main result and the building blocks of the proof.

The main difficulty to construct a marginally trapped surface is to find a spacelike surface with vanishing outgoing null expansion. Formally, let $\tr \chi(\Sigma)$ denote the outgoing null expansion of a spacelike surface $\Sigma$, we need to solve the equation
\begin{align*}
\tr \chi(\Sigma) =0.
\end{align*}
In order to solve it, we shall translate the above formal equation to a precise equation in analysis. This is done in two steps: first we find the method to parameterise $\Sigma$ by functions, second we calculate the outgoing null expansion in terms of the parameterisation functions of $\Sigma$. Then we translate the formal equation $\tr \chi(\Sigma)=0$ to a precise equation for the parameterisation functions. 

Given the precise form of the equation $\tr \chi(\Sigma)=0$, we use the perturbative method to solve it. We construct an appropriate linearised perturbation of the outgoing null expansion, then apply it to construct approximating solutions of the equation $\tr \chi(\Sigma)=0$, and eventually prove that the approximating solutions converge to the actual solution.

The major building blocks of the proof consist of the followings.
\begin{enumerate}[label=\roman*.]
\item \emph{Parameterisations of spacelike surfaces in section \ref{sec 2}}. We introduce two methods to parameterise a spacelike surface. The first one is simply parameterising a surface by two functions as their graph of $\{\us, s\}$ in the double null coordinate system. 

The second method of parameterisation is less direct than the first one. Suppose that $\Sigma$ is  a spacelike surface embedding in an incoming null hypersurface $\ucalH$. Then the second method also parameterises $\Sigma$ by two functions, where the first function is to parameterise the incoming null hypersurface $\ucalH$, and the second function is to parameterise the position of $\Sigma$ inside $\ucalH$. The transformation from the second method of parameterisation to the first one is studied in details.

Clearly the second method of parameterisation is more natural for studying spacelike surfaces in an incoming null hypersurface, while the first method has the advantage of its directness when evaluating the background geometric quantities on a spacelike surface, like the metric components and structure coefficients.

\item \emph{Formula of the outgoing null expansion in section \ref{sec 3}}. We use a two-step procedure to obtain the formula of the outgoing null expansion, which is naturally coherent with the second method of parameterising spacelike surfaces. Suppose that $\Sigma$ is a spacelike surface embedded in an incoming null hypersurface $\ucalH$. The first step is to obtain the geometric information along $\ucalH$, then the second step is to calculate the outgoing null expansion with this information and the location of $\Sigma$ inside $\ucalH$. The precise formula is formula \eqref{eqn 4.3}. A decomposition of the formula into the first order main part and high order remainder part is introduced in subsection \ref{subsec 4.3}.

The second step has been investigated in prior works \cite{KLR} \cite{L1} \cite{An}, while these works were restricted to the case of $\ucalH$ being the incoming null hypersurface $\uC_{\us}$, the level set of the coordinate function $\us$. Combining with the first step, we extend the formula to spacelike surfaces embedded in a more general class of incoming null hypersurfaces. 

In this more general case, an extra difficulty arises. In the first step obtaining the geometric information along $\ucalH$, we not only need to evaluate the background geometric quantities on $\ucalH$, but also need to know the differential of the parameterisation function $\uh$ of $\ucalH$. This difficulty doesnot appear in the case that $\ucalH =\uC_{\us}$, as the differential of $\uh = \us$ vanishes. This extra difficulty is addressed in subsections \ref{subsec 5.1}, \ref{subsec 5.2}, \ref{subsec 7.1}.

\item \emph{Perturbation and linearised perturbation of the parameterisation of spacelike surfaces in sections \ref{sec 6} and \ref{sec 8}}. More precisely, we study the perturbation and linearised perturbation of the transformation from the second method of parameterising spacelike surfaces to the first method. They are essential for the study of perturbation and linearised perturbation of the outgoing null expansion.

\item \emph{Perturbation and linearised perturbation of the outgoing null expansion in sections \ref{sec 7} and \ref{sec 9}}. They are the most important building blocks for solving the equation $\tr \chi(\Sigma)=0$. 

\end{enumerate}

We emphasis a key structure in the formula of the outgoing null expansion. Suppose that a spacelike surface $\Sigma$ is embedded in an incoming null hypersurface $\ucalH$, and the parameterisation function $f$ gives the location of $\Sigma$ inside $\ucalH$. Then the outgoing null expansion $\tr \chi(\Sigma)$ is given by a quasi-linear elliptic operator on the parameterisation function $f$: let $\vartheta$ be the variable on the sphere and $\circnabla$ be the Levi-Civita connection on the standard round sphere of radius $1$, then
\begin{align*}
\tr \chi(\Sigma) = \tr \chi ( \ucalH; f) = a^{ij} (\ucalH; \vartheta, f, \d f) \circnabla^2_{ij} f + c(\ucalH; \vartheta, f, \d f),
\end{align*}
where $a^{ij} (\ucalH; \vartheta, f, \d f)$, $c(\ucalH; \vartheta, f, \d f)$ are both fully nonlinear. Although the precise formula of $\tr \chi(\Sigma)$ is involved in general, its first order main part reduces to a simple elliptic operator of $f$ in the case of the Schwarzschild spacetime, as showed in formula \eqref{eqn 4.4},
\begin{align*}
- 2r_S^{-2} \circDelta f + \tr \chi'_S,
\end{align*}
where $\circDelta$ is the Laplacian on the standard round sphere of radius $1$. $r_S$ is the radius of the round sphere $\Sigma_{\us,s}$ and $\tr \chi'_S$ is the outgoing null expansion of $\Sigma_{\us,s}$, both in the Schwarzschild spacetime.

We conclude the introduction by pointing out another important feature on the regularity of the solution of $\tr \chi(\ucalH; f)=0$. In the above formula of $\tr \chi$, the terms $a^{ij} (\ucalH; \vartheta, f, \d f)$, $c(\ucalH; \vartheta, f, \d f)$ depend on $\ucalH$ in a complicate way. We introduce a parameterisation function $\ufl{s=0}$ in subsection \ref{subsec 3.1} to determine the location of the incoming null hypersurface $\ucalH$, then we write these terms as $a^{ij}([\ufl{s=0}];  \vartheta, f, \d f)$, $c([\ufl{s=0}]; \vartheta, f, \d f)$. The dependence on $\ufl{s=0}$ of $a^{ij}$ and $c$ is not a pointwise dependence on the values of $\ufl{s=0}$ and its derivatives, but in a functional way. In order to solve the equation $\tr \chi ( \ufl{s=0}, f )= \tr \chi ( \ucalH; f ) = 0$, we need to obtain the regularitites of $a^{ij} ( [ \ufl{s=0}, \vartheta, f, \d f)$ and $c( [ \ufl{s=0}], \vartheta, f, \d f)$. With the estimates obtained in section \ref{sec 5}, we can show that if $\ufl{s=0}$ is a function in the Sobolev space $\mathrm{W}^{n+2,p}$ and $f$ is smooth, then $c([\ufl{s=0}], \vartheta, f, \d f)$ is in the Sobolev space $\mathrm{W}^{n,p}$. Thus the theory of elliptic equations on the sphere tells us that the regularity of the solution $f$ of $\tr \chi ( \ufl{s=0}, f )=0$ cannot be better than the regularity of $\ufl{s=0}$ in general, i.e. roughly speaking, the regularity of the embedding of the marginally trapped surface $\Sigma$ in $\ucalH$ cannot exceed the regularity of the embedding of $\ucalH$ in the spacetime in general. The precise statement on the regularity of the marginally trapped surface is contained in the main theorem \ref{thm 10.9}.

\section{Perturbation of Schwarzschild metric near event horizon}\label{sec 2}
In this section, we introduce the perturbed Schwarzschild metric considered in this paper. Since we are interested in the spacetime near the event horizon, we shall only consider the perturbation of the Schwarzschild metric in a neighbourhood of a piece of the event horizon.

In the coordinate system $\{t,r,\theta,\phi\}$, the Schwarzschild metric $g_{S}$ is
\begin{align*}
g_{S}=-\left( 1-\frac{2m}{r} \right) \d t^2 + \left( 1-\frac{2m}{r} \right)^{-1} \d r ^2 + r^{2} \left( \d \theta^2 + \sin^2 \theta \d \phi^2 \right).
\end{align*}
We denote $2m$ by $r_0$. We consider the following coordinate transformation which is used in \cite{L4}
\begin{align*}
\left\{
\begin{aligned}
&
(r-2m)^{\frac{1}{2}} \exp \frac{t+r}{4m}=\exp\frac{\us}{r_0},
\\
&
(r-2m)^{\frac{1}{2}} \exp \frac{-t+r}{4m}=s \exp \frac{s+r_0}{r_0}.
\end{aligned}
\right.
\end{align*}
Then in the coordinate system $\{\us, s, \theta, \phi\}$, the Schwarzschild metric $g_S$ takes the form
\begin{align*}
g_{S} = 2 \Omega_S^2 \left( \d s\otimes \d \us + \d \us \otimes \d s  \right) + r^2 \left( \d \theta^2 + \sin^2 \theta \d \phi^2 \right),
\end{align*}
where
\begin{align}
\begin{aligned}
&
\Omega_S^2 = \frac{s+r_0}{r}\exp\frac{\us+s+r_0-r}{r_0},
\\
&
(r-r_0)\exp\frac{r}{r_0} = s\exp\frac{\us+s+r_0}{r_0}.
\end{aligned}
\label{eqn 2.1}
\end{align}
This coordinate system $\{\us, s, \theta, \phi\}$ is a double null coordinate system, where the level set $s=0$ is the event horizon $r=r_0$ of the Schwarzschild black hole. 

We denote the level sets of $s$ by $C_s$, the level sets of $\us$ by $\uC_{\us}$ and use $\Sigma_{\us,s}$ to denote the intersection of $\uC_{\us}$ with $C_s$. $\Sigma_{\us,s}$ is a round sphere of radius $r$.

The sphere $\Sigma_{\us, s=0}$ is marginally trapped. We shall consider the perturbation of the Schwarzschild metric in a neighbourhood of the marginally trapped surface $\Sigma_{0,0}$, thus we introduce the so-called $(\kappa,\underline{\kappa})$-neighbourhood of $\Sigma_{0,0}$.

\begin{definition}[$(\kappa,\underline{\kappa})$-neighbourhood $M_{\kappa,\underline{\kappa}}$ of $\Sigma_{0,0}$]\label{def 2.1}
Let $\{\us,s\}$ be the double null coordinates of the Schwarzschild spacetime $\left(\mathcal{S},g_{S}\right)$ introduced above. The $(\kappa,\underline{\kappa})$-neighbourhood $M_{\kappa,\underline{\kappa}}$ of the marginally trapped surface $\Sigma_{0,0}$ is defined by
\begin{align*}
M_{\kappa,\underline{\kappa}} = \left\{ p\in \mathcal{S}: \leve s(p) \rive < \kappa r_0, \leve \us \rive < \underline{\kappa} r_0   \right\}.
\end{align*}
See figure \ref{fig 3}. In this paper, we assume that $\kappa, \underline{\kappa} \leq 0.1$. \footnote{By an elementary estimate, one can derive that $r/r_0 \in (0.8, 1.2)$ in $M_{\kappa, \underline{\kappa}}$. \label{footnote 1}}
\begin{figure}[h]
\begin{center}
\begin{tikzpicture}
\draw[dotted,thick] (1, 0) -- (0,1) -- (-1,0) -- (0,-1) -- cycle;
\draw[->] (0,0) -- (0.5,-0.5) node[right] {$s$};
\draw[->] (0,0) -- (0.5,0.5) node[right] {$\us$};
\draw[fill] (0,0)  circle [radius=0.03];
\draw[->] (-0.5,-0.5) node[left] {\footnotesize $\Sigma_{0,0}$} to [out=30,in=-120] (-0.06,-0.06);
\draw[->] (-0.75,0.75) node[left] {\footnotesize$M_{\kappa,\underline{\kappa}}$: $(\kappa,\underline{\kappa})$-neighbourhood }to [out=-30,in=110] (-0.25,0.25);
\end{tikzpicture}
\end{center}
\caption{$(\kappa,\underline{\kappa})$-neighbourhood $M_{\kappa,\underline{\kappa}}$ of $\Sigma_{0,0}$}
\label{fig 3}
\end{figure}
\end{definition}

Let $\{\us, s, \theta^1, \theta^2\}$ be the coordinate system on $M_{\kappa,\underline{\kappa}}$ inherited from the double null coordinate system of $(M_{\kappa,\underline{\kappa}}, g_{S})$. Consider a class of Lorentzian metrics $g$ on $M_{\kappa,\underline{\kappa}}$ which has the following form
\begin{align*}
g=2\Omega^2  \left( \d s  \otimes \d \us  +\d \us \otimes \d s \right) + \slashg_{ab} \big( \d \theta^a - b^{a} \d s \big) \otimes \big( \d \theta^b - b^{b} \d s \big).
\end{align*}
The pair $(M_{\kappa,\underline{\kappa}}, g)$ is a Lorentzian manifold and $\{\us, s, \vartheta^1, \vartheta^2\}$ is a double null coordinate system of this manifold. We define the associated null frame $\{ L, \uL \}$ where $L, \uL$ are tangential null vector fields of $C_s, \uC_{\us}$ respectively and
\begin{align*}
L\us=1, \quad \uL s=1.
\end{align*}
Hence we have
\begin{align*}
L = \partial_{\us}, \quad \uL= \partial_s + \vec{b} = \partial_s + b^i \partial_i.
\end{align*}
We also introduce the null vectors $L', \uL'$ by
\begin{align*}
L' = \Omega^{-2} L, \quad \uL' = \Omega^{-2} \uL,
\end{align*}
then $\{ \uL, L' \}$ and $\{ \uL', L \}$ are both conjugate null frames on $\Sigma_{\us,s}$ since
\begin{align*}
g(\uL, L') = g(\uL', L) =2.
\end{align*}
With respect to the double null coordinate system $\{ \us, s\}$, we can define the structure coefficients as follows.
\begin{definition}[Structure coefficients]
Let $X,Y$ be tangential vector fields of $\Sigma_{\us, s}$. Define that
\begin{align}
&
\underline{\omega}=\underline{L} \log \Omega, 
&&
\omega= L \log \Omega, 
\nonumber
\\
\nonumber
&
\underline{\chi}(X,Y) = g(\nabla_X \underline{L},Y),  
&& 
\chi(X,Y)= g(\nabla_X L, Y), 
\\
\nonumber
&
\underline{\chi}'(X,Y) =g(\nabla_X \underline{L}',Y),  
&&
\chi'(X,Y)= g(\nabla_X L', Y),  
\\
\nonumber
&
\eta(X) = \frac{1}{2} g ( \nabla_X \underline{L}, L'),  
&&
\underline{\eta}(X) = \frac{1}{2}g(\nabla_X L, \underline{L}').
\end{align}
$\uchi,\chi,\uchi',\chi'$ are the null second fundamental forms in the directions of $\uL,L,\uL',L'$ respectively. $\eta,\ueta$ are the torsions of the null frame $\{\uL,L'\}$ and $\{\uL' ,L\}$ respectively. We can decompose them into trace and trace-free parts with respect to the metric $\slashg$,
\begin{eqnarray}
&\uchi=\hatuchi + \frac{1}{2} \tr \uchi \slashg, & \chi= \hatchi + \frac{1}{2} \tr \chi \slashg,
\nonumber
\\
\nonumber
&\uchi'=\hatuchi' + \frac{1}{2} \tr \uchi' \slashg, & \chi'= \hatchi' + \frac{1}{2} \tr \chi' \slashg.
\end{eqnarray}
The trace-free parts are called shears, and the traces are called null expansions.
\end{definition}

For the Schwarzschild metric, the structure coefficients associated with the double null coordinate system $\{ \us, s\}$ are given by the following formulae,
\begin{align}
\begin{aligned}
&
\partial_s r= \frac{s+r_0}{r} \cdot \frac{r-r_0}{s},
\quad
\partial_{\us} r= \frac{r-r_0}{r}, 
\\
&
\tr \chi_{S}=\frac{2(r-r_0)}{r^2} = \frac{2s}{r^2} \exp\frac{\us+s+r_0-r}{r_0} ,
\\
&
\tr\uchi_{S} = \frac{2(s+r_0)}{r^2} \cdot \frac{r-r_0}{s} = \frac{2(s+r_0)}{r^2} \exp\frac{\us+s+r_0-r}{r_0},
\\
&
\hatchi_{S}=\hatchi'_{S}=0,
\quad 
\hatuchi_{S}=\hatuchi'_{S}=0, 
\quad 
\eta_{S}=\ueta_{S}=0,
\\
&
\omega_{S}=\frac{r_0}{2r^2}, 
\quad 
\uomega_{S}=\frac{1}{2(s+r_0)} +\frac{1}{2r_0} - \Big( \frac{1}{2r} + \frac{1}{2r_0} \Big)  \frac{s+r_0}{r} \exp \frac{\us+s+r_0-r}{r_0}.
\end{aligned}
\label{eqn 2.2}
\end{align}

In the following, we introduce a class of perturbations of the Schwarzschild metric on $M_{\kappa,\underline{\kappa}}$.
\begin{definition}\label{def 2.3}
Let $g_{\epsilon}$ be a Lorentzian metric on $M_{\kappa,\underline{\kappa}}$. In coordinates system $\{ \us, s, \theta^1, \theta^2\}$, it takes the form
\begin{align*}
g_{\epsilon}=2\Omega_{\epsilon}^2  \left( \d s  \otimes \d \us  +\d \us \otimes \d s \right) + \left(\slashg_{\epsilon} \right) _{ab} \big( \d \theta^a - b^{a} \d s\big) \otimes \big( \d \theta^b - b^{b} \d s\big).
\end{align*}
Let $r_{\epsilon}(\us,s)$ be the area radius of the surface $(\Sigma_{\us,s}, \slashg_{\epsilon})$
\begin{align*}
4 \pi r_{\epsilon}^2(\us,s) = \int_{\Sigma_{\us,s}} 1 \cdot \dvol_{\slashg}.
\end{align*}
$g_{\epsilon}$ is a perturbation of the Schwarzschild metric in the following sense: the metric components of $g_{\epsilon}$ are close to the ones of $g_S$, and the structure coefficients for $g_{\epsilon}$ are also close to the ones for $g_S$. Let $n$ be a positive integer. The precise quantitative descriptions are given by the following formulae. For the metric components:\footnote{The following derivatives with respect to $\partial_s, \partial_{\us}$ are all Lie derivatives.}
\begin{align*}
&
\left.
\begin{aligned}
&  1-\epsilon < \leve \frac{r_{\epsilon}}{r_{S}} \rive < 1+ \epsilon, 
\end{aligned}
\right.
\\
&
\left.
\begin{aligned}
&  \leve \circnabla^k \partial_s^l \partial_{\us}^m \left( \log \Omega_{\epsilon} - \log \Omega_{S} \right) \rive < \frac{\epsilon }{r_0^{m+l}},
\end{aligned}
\right.
\\
&
\left.
\begin{aligned}
& \leve \circnabla^k \partial_{s}^l \vec{b}_{\epsilon} \rive_{\circg} \leq \frac{\epsilon \leve\us\rive}{r_0^{l+2}},
&& \leve \circnabla^k \partial_s^l \partial_{\us}^m \vec{b}_{\epsilon} \rive_{\circg} < \frac{\epsilon}{r_0^{m+l+1}},
\end{aligned}
\right.
\\
&
\left.
\begin{aligned}
& \leve \circnabla^k \partial_{s}^l \partial_{\us}^m \left( \slashg_{\epsilon} - \slashg_{S} \right) \rive_{\circg} < \frac{\epsilon}{ r_0^{m+l-2}},
\end{aligned}
\right.
\end{align*}
where $k+m+l \leq n+2$.
For the structure coefficients:
\begin{align*}
&
\left.
\begin{aligned}
& \leve \circnabla^k \partial_{s}^l \partial_{\us}^m   \left(  \tr\chi_{\epsilon} - \tr\chi_{S} \right) \rive_{\circg} < \frac{\epsilon}{r_0^{m+l+1}},
\end{aligned}
\right.
&&
\quad
\left.
\begin{aligned}
& \leve \circnabla^k \partial_{s}^l \partial_{\us}^m  \left( \tr\uchi_{\epsilon} - \tr\uchi_{S} \right) \rive_{\circg} < \frac{\epsilon}{r_0^{m+l+1}},
\end{aligned}
\right.
\\
&
\left.
\begin{aligned}
& \leve \circnabla^k \partial_{s}^l \partial_{\us}^m  \hatchi_{\epsilon} \rive_{\circg} < \frac{\epsilon}{ r_0^{m+l-1}},
\end{aligned}
\right.
&&
\quad
\left.
\begin{aligned}
& \leve \circnabla^k \partial_{s}^l \partial_{\us}^m  \hatuchi_{\epsilon} \rive_{\circg} < \frac{\epsilon}{ r_0^{m+l-1} }.
\end{aligned}
\right.
\\
&
\left.
\begin{aligned}
& \leve \circnabla^k \partial_{s}^l \partial_{\us}^m  \eta_{\epsilon} \rive_{\circg} < \frac{\epsilon }{r_0^{m+l}},
\end{aligned}
\right.
&&
\quad
\left.
\begin{aligned}
& \leve \circnabla^k \partial_{s}^l \partial_{\us}^m  \ueta_{\epsilon} \rive_{\circg} < \frac{\epsilon}{r_0^{m+l}},
\end{aligned}
\right.
\\
&
\left.
\begin{aligned}
& \leve \circnabla^k \partial_{s}^l \partial_{\us}^m (\omega_{\epsilon} -\omega_{S}) \rive_{\circg} < \frac{\epsilon}{r_0^{m+l+1}},
\end{aligned}
\right.
&&
\quad
\left.
\begin{aligned}
& \leve \circnabla^k \partial_{s}^l \partial_{\us}^m  (\uomega_{\epsilon} -  \uomega_{S}) \rive_{\circg} < \frac{\epsilon }{r_0^{m+l+1}},
\end{aligned}
\right.
\end{align*}
where $k+m+l \leq n+1$.
In the above formulae, $\circg$ is the metric of the standard round sphere with radius one on the surface $\Sigma_{\us,s}$ and $\circnabla$ is the corresponding covariant derivatives of $\circg$.
\end{definition}

Note that in the above definition, $(M_{\kappa, \underline{\kappa}}, g_{\epsilon})$ is not necessary vacuum. Moreover, there is no guarantee that the Schwarzschild event horizon $C_{s=0}$ contains a marginally trapped surface in $(M_{\kappa, \underline{\kappa}}, g_{\epsilon})$ anymore, and none of the coordinate surface $\Sigma_{\us,s}$ is assumed to be marginally trapped.

For the rest of this paper, we will simply use $M, g$ to denote $M_{\kappa,\underline{\kappa}}, g_{\epsilon}$, and omit $\epsilon$ in the lower indices of the metric components and structure coefficients to simplify the notations.

\section{Parameterisation of spacelike surface}\label{sec 3}
In this section, we introduce two methods to parametrise a spacelike surface $\Sigma$ in $(M,g)$. Moreover, we will describe the transformation between these two parameterisations and give an estimate for the parameterisation transformation, which is proposition \ref{pro 3.3}.

\subsection{Two methods to parametrise spacelike surface}\label{subsec 3.1}
The first kind of parameterisation is to simply consider the surface $\Sigma$ as a graph over $\mathbb{S}^2$ in the double null coordinate system.

\begin{figure}[h]
\begin{center}
\begin{tikzpicture}
\draw[dashed] (-1,0)
to [out=70, in=180] (0,0.5)
to [out=0,in=110] (1,0);
\draw (1,0)
to [out=-70,in=0] (0,-0.8) node[below]{\tiny $\Sigma_{0,0}$}
to [out=180,in=-110] (-1,0); 
\draw[dashed] (-1,0) to [out=70,in=-110] (-0.7,0.9);
\draw (-1,0) to [out=-110,in=70] (-1.9,-2.5);
\draw[dashed] (1,0) to [out=110,in=-70] (0.7,0.9);
\draw[->] (1,0) to [out=-70,in=110] (1.9,-2.5) node[right] {\small $s$}; 
\draw[dashed] (-1,0) to [out=-45,in=135] (-0.5,-0.5);
\draw (-1,0) to [out=135, in= -45] (-2,1);
\draw[dashed] (1,0) to [out=-135,in=45] (0.5,-0.5);
\draw[->] (1,0) to [out=45, in= -135] (2,1) node[right] {\small $\us$}; 
\node[below right] at (1.4,0.55) {\tiny $C_{s=0}$};
\node[above right] at (1.7,-2) {\tiny $\uC_{\us=0}$};
\draw[dashed] (-2.7,-2.2+1) to [out=70,in=180] (-0.5,-2.5+1) to [out=0,in=110] (2.7,-2.2+1);
\draw (2.7,-2.2+1) node[right] {\scriptsize $\Sigma=\{(\us,s,\vartheta): \us=\uf(\vartheta), s= f(\vartheta)\}$} to [out=-70,in=0] (1,-2+0.8) to [out=180,in=-110] (-2.7,-2.2+1); 
\end{tikzpicture}
\end{center}
\caption{The first kind of parameterisation of $\Sigma$.}
\label{fig 4}
\end{figure}
As demonstrated in figure \ref{fig 4}, the surface $\Sigma$ is the graph of the pair of functions $(\uf, f)$ over the domain of variable $\vartheta$ in the double coordinate system $\{ \us, s, \vartheta \}$,
\begin{align*}
\Sigma=\{(\us,s,\vartheta): \us=\uf(\vartheta), s= f(\vartheta)\}.
\end{align*}
We call $(\uf, f)$ the first kind of parameterisation of $\Sigma$.

The second kind of parameterisation is to consider the incoming null hypersurface $\ucalH$ containing $\Sigma$ and the embedding of $\Sigma$ in the incoming null hypersurface. This parameterisation can be visualised by figure \ref{fig 5}.
\begin{figure}[h]
\begin{center}
\begin{tikzpicture}
\draw[dashed] (-1,0)
to [out=70, in=180] (0,0.5)
to [out=0,in=110] (1,0);
\draw (1,0)
to [out=-70,in=0] (0,-0.8) node[below]{\tiny $\Sigma_{0,0}$}
to [out=180,in=-110] (-1,0); 
\draw[dashed] (-1,0) to [out=70,in=-110] (-0.7,0.9);
\draw (-1,0) to [out=-110,in=70] (-1.85,-2.4);
\draw[dashed] (1,0) to [out=110,in=-70] (0.7,0.9);
\draw[->] (1,0) to [out=-70,in=110] (1.85,-2.4) node[right] {\small $s$}; 
\node[above right] at (1.6,-1.9) {\tiny $\uC_{\us=0}$};
\draw[dashed] (-1,0) to [out=-45,in=135] (-0.5,-0.5);
\draw (-1,0) to [out=135, in= -45] (-2,1);
\draw[dashed] (1,0) to [out=-135,in=45] (0.5,-0.5);
\draw[->] (1,0) to [out=45, in= -135] (1.8,0.8) node[right] {\tiny $C_{s=0}$}to [out=45,in=-135] (2.3,1.3) node[right] {\small $\us$}; 
\draw[dashed] (-1.5,0.5) to [out=135,in=45] (1.5,0.5);
\draw (1.5,0.5) to [out=-135,in=0] (0,0) node[below]{\scriptsize $\Sigma_0$} to [out=180,in=-45] (-1.5,0.5); 
\draw[dashed] (-1.4,1.3) to [out=-110,in=70] (-1.6,0.7);
\draw (-1.6,0.7) to [out=-110,in=70] (-2.75,-2.4);
\draw[dashed] (1.4,1.3) to [out=-70,in=110] (1.6,0.7);
\draw[->] (1.6,0.7) to [out=-70,in=110] (2.75,-2.4) node[right]{\small $s$}; 
\node[right] at (2.4,-1.5) {\scriptsize $\ucalH$}; 
\draw[dashed] (-2.3,-2.2+1) to [out=70,in=180] (-0.5,-2.5+1) node[below] {\scriptsize $\Sigma$} to [out=0,in=110] (2.3,-2.2+1);
\draw (2.3,-2.2+1) to [out=-70,in=0] (1,-2+0.9) to [out=180,in=-110] (-2.3,-2.2+1); 
\draw[->] (-1.73,0.73) to [out=-110,in=70] (-2.45,-1.2);
\node[left] at (-2.3,-0.8) {$f$};
\draw[->] (-1.05,-0.1) to [out=135,in=-45] (-1.65,0.5);
\node[below] at (-1.55,0.1) {$\ufl{s=0}$};
\end{tikzpicture}
\end{center}
\caption{The second kind of parameterisation of $\Sigma$.}
\label{fig 5}
\end{figure}

Let $\Sigma_0$ be the intersection of $\ucalH$ with $C_{s=0}$. Assume that $\Sigma_0$ is parametrised by a function $\ufl{s=0}$ as its graph of $\us$ over the domain of $\vartheta$ in $\{\us, \vartheta\}$ coordinate system in $C_{s=0}$,
\begin{align*}
\Sigma_0 = \{ (\us,s=0,\vartheta): \us= \ufl{s=0}(\vartheta) \}.
\end{align*}
The restriction of the double null coordinates $\{ s, \vartheta \}$ on the incoming null hypersurface $\ucalH$ is a coordinate system. Then we assume that $\Sigma$ is parametrised by a function $f$ as its graph of $s$ over the domain of $\vartheta$ in $\{ s, \vartheta \}$ coordinate system on $\ucalH$
\begin{align*}
\Sigma = \{ (s, \vartheta) \in \ucalH: s= f(\vartheta) \}.
\end{align*}
We define $(\ufl{s=0}, f)$ to be the second kind parameterisation of $\Sigma$.

For the rest of the paper, we will simply call the first or second kind of parameterisation as the first or second parameterisation for the sake of brevity.

\subsection{Transformation from the second parameterisation to the first}\label{subsec 3.2}
We want to know the transformation from the second parameterisation to the first one. Assume that $\Sigma$ has the second parameterisation $(\ufl{s=0}, f)$. Since the first parameterisation shares the same parameterisation function $f$ for the $s$ coordinate, we need only to determine the parametrisation function $\uf$ for the $\us$ coordinate. In the following, we introduce two methods to obtain $\uf$.

\begin{method I}
$\Sigma$ is embedded in the incoming null hypersurface $\ucalH$. Let $\Sigma_{s}$ be the intersection of $\ucalH$ with $C_s$. Suppose that $\Sigma_s$ has the first parameterisation $(\ufls, s)$
\begin{align*}
\Sigma_s = \{ (\us, s, \vartheta): \us = \ufls(\vartheta) \}.
\end{align*}
We introduce the parameterisation function $\uh$ for $\ucalH$
\begin{align*}
\uh(s,\vartheta) = \ufls(\vartheta).
\end{align*}
Then $\ucalH$ is parametrised as
\begin{align*}
\ucalH= \{ (\us, s, \vartheta): \us = \uh(s,\vartheta), s\in (-\kappa r_0, \kappa r_0) \}.
\end{align*}
The above parameterisation of the incoming null hypersurface $\ucalH$ is also used in the work \cite{L4}. It can be illustrated by figure \ref{fig 6}
\begin{figure}[h]
\begin{center}
\begin{tikzpicture}
\draw[dashed] (-1,0)
to [out=70, in=180] (0,0.5)
to [out=0,in=110] (1,0);
\draw (1,0)
to [out=-70,in=0] (0,-0.8) node[below]{\tiny $\Sigma_{0,0}$}
to [out=180,in=-110] (-1,0); 
\draw[dashed] (-1,0) to [out=70,in=-110] (-0.7,0.9);
\draw (-1,0) to [out=-110,in=70] (-2.4,-4);
\draw[dashed] (1,0) to [out=110,in=-70] (0.7,0.9);
\draw[->] (1,0) to [out=-70,in=110] (2.4,-4) node[right] {\small $s$}; 
\node[above right] at (1.3,-1.3) {\tiny $\uC_{\us=0}$};
\draw[dashed] (-1,0) to [out=-45,in=135] (-0.5,-0.5);
\draw (-1,0) to [out=135, in= -45] (-2,1);
\draw[dashed] (1,0) to [out=-135,in=45] (0.5,-0.5);
\draw[->] (1,0) to [out=45, in= -135] (1.8,0.8) node[right] {\tiny $C_{s=0}$}to [out=45,in=-135] (2.3,1.3) node[right] {\small $\us$}; 
\draw[dashed] (-1.5,0.5) to [out=135,in=45] (1.5,0.5);
\draw (1.5,0.5) to [out=-135,in=0] (0,0) node[below]{\scriptsize $\Sigma_0$} to [out=180,in=-45] (-1.5,0.5); 
\draw[dashed] (-1.4,1.3) to [out=-110,in=70] (-1.6,0.7);
\draw (-1.6,0.7) to [out=-110,in=70] (-3,-3);
\draw[dashed] (1.4,1.3) to [out=-70,in=110] (1.6,0.7);
\draw[->] (1.6,0.7) to [out=-70,in=110] (3,-3) node[right]{\small $s$}; 
\node[right] at (2.1,-0.7) {\scriptsize $\ucalH, \ \us = \uh(s,\vartheta) = \ufls(\vartheta) $};
\draw[dashed] (-1.5,-3.5) -- (-2,-3);
\draw (-2,-3) -- (-3.5,-1.5);
\draw[dashed] (1.5,-3.5) -- (2,-3);
\draw[->] (2,-3) -- (3.5,-1.5); 
\node[right] at (3,-2) {\tiny $C_s$};
\draw[dashed] (-2.5,-2.5) to [out=135,in=45] (2.5,-2.5);
\draw (2.6,-2.4) to [out=-135,in=0] (0,-3.4) node[below]{\scriptsize $\Sigma_s$} to [out=180,in=-45] (-2.6,-2.4); 
\end{tikzpicture}
\end{center}
\caption{Parametrisation of $\ucalH$.}
\label{fig 6}
\end{figure}

\noindent
Since $\Sigma$ is the graph $s=f(\vartheta)$ in $\ucalH$, then in the double coordinate system $\Sigma$ is given by
\begin{align*}
\{ (\us, s, \vartheta): \ s= f(\vartheta),  \us=\uh(f(\vartheta), \vartheta) \}.
\end{align*}
Therefore the parameterisation function $\uf$ is given by
\begin{align*}
\uf(\vartheta) = \uh(f(\vartheta), \vartheta) = \left( \ufl{s= f(\vartheta)} \right) (\vartheta).
\end{align*}
See figure \ref{fig 7}.
\begin{figure}[h]
\begin{center}
\begin{tikzpicture}
\draw (-1,0) to [out=135, in= -45] (-2,1);
\draw[->] (1,0) to [out=45, in= -135] (1.8,0.8) node[right] {\tiny $C_{s=0}$}to [out=45,in=-135] (2.3,1.3) node[right] {\small $\us$};
\draw[dashed] (-1.5,0.5) to [out=135,in=45] (1.5,0.5);
\draw (1.5,0.5) to [out=-135,in=0] (0,0) node[below]{\scriptsize $\Sigma_0$} to [out=180,in=-45] (-1.5,0.5); 
\draw[dashed] (-1.4,1.3) to [out=-110,in=70] (-1.6,0.7);
\draw (-1.6,0.7) to [out=-110,in=70] (-3,-3);
\draw[dashed] (1.4,1.3) to [out=-70,in=110] (1.6,0.7);
\draw[->] (1.6,0.7) to [out=-70,in=110] (3,-3) node[right]{\small $s$}; 
\node[right] at (2.4,-1.5) {\scriptsize $\ucalH$};
\draw (-2,-3) -- (-3.5,-1.5);
\draw[->] (2,-3) -- (3.5,-1.5); 
\node[right] at (3,-2) {\tiny $C_s$};
\draw[dashed] (-2.5,-2.5) to [out=135,in=45] (2.5,-2.5);
\draw (2.6,-2.4) to [out=-135,in=0] (0,-3.4) node[below]{\scriptsize $\Sigma_s$} to [out=180,in=-45] (-2.6,-2.4); 
\draw[dashed] (-2.3+0.23,-2.2+1.6) to [out=70,in=180] (-0.5,-2.5+1.6) to [out=0,in=110] (2.3-0.23,-2.2+1.6) node[right] { $\Sigma, \ \us=\uh(f(\vartheta),\vartheta)$};
\draw (2.3-0.23,-2.2+1.6) to [out=-70,in=0] (1,-2+1.5) to [out=180,in=-110] (-2.3+0.23,-2.2+1.6); 
\draw[->] (-1.7,0.7) to [out=-110,in=70] (-2.45+0.28,-1.2+0.63);
\node[left] at (-2,-0.2) {$s=f(\vartheta)$}; 
\end{tikzpicture}
\end{center}
\caption{Method I}
\label{fig 7}
\end{figure}
Thus in order to obtain $\uf$, it is sufficient to know $\uh$ or $\ufls$. We apply the equations derived in \cite{L4}, that $\ufls$ satisfies
\begin{align}
\label{eqn 3.1}
\dpartial_{s} \ufls=  -b^i  \dpartial_i \ufls + \Omega^2 \left(\slashg^{-1}\right)^{ij} \dpartial_i \ufls \dpartial_j \ufls,
\end{align}
and equivalently $\uh$ satisfies
\begin{align}
\label{eqn 3.2}
\dpartial_{s} \uh=  -b^i  \dpartial_i \uh + \Omega^2 \left(\slashg^{-1}\right)^{ij} \dpartial_i \uh \dpartial_j \uh.
\end{align}
Here the notations $\dpartial_s, \dpartial_i$ denote the coordinate derivatives of the coordinate system $\{ s, \vartheta \}$ on $\ucalH$.
\end{method I}

\begin{method II}
Introduce a family of surfaces $\{S_t\}$ in $\ucalH$, where $S_t$ has the second parameterisation $(\ufl{s=0}, \flt = t f)$. $S_{t=0}$ is simply $\Sigma_0$ and $S_{t=1}$ is $\Sigma$. Thus $\{S_t\}$ is a family of surfaces deforming from $\Sigma_0$ to $\Sigma$. See figure \ref{fig 8}.
\begin{figure}[h]
\begin{center}
\begin{tikzpicture}
\draw (-1,0) to [out=135, in= -45] (-2.3,1.3);
\draw[->] (1,0) to [out=45, in= -135] (1.8,0.8) node[right] {\tiny $C_{s=0}$}to [out=45,in=-135] (2.3,1.3) node[right] {\small $\us$};
\draw[dashed] (-1.5,0.5) to [out=135,in=45] (1.5,0.5);
\draw (1.5,0.5) to [out=-135,in=0] (0,0) node[above]{\scriptsize $S_{t=0}:\Sigma_0$} to [out=180,in=-45] (-1.5,0.5); 
\draw[dashed] (-1.4,1.3) to [out=-110,in=70] (-1.6,0.7);
\draw (-1.6,0.7) to [out=-110,in=70] (-3,-3);
\draw[dashed] (1.4,1.3) to [out=-70,in=110] (1.6,0.7);
\draw[->] (1.6,0.7) to [out=-70,in=110] (3,-3) node[right]{\small $s$}; 
\node[right] at (2.4,-1.5) {\scriptsize $\ucalH$};
\draw[dashed] (-2.3+0.23,-2.2+1.6) to [out=70,in=180] (-0.4,-0.7) node[below] {\scriptsize $S_t$} to [out=0,in=110] (2.3-0.23,-2.2+1.6);
\draw (2.3-0.23,-2.2+1.6) to [out=-70,in=0] (0.7,-0.5) to [out=180,in=-110] (-2.3+0.23,-2.2+1.6); 
\draw[->] (-1.71,0.71) to [out=-110,in=70] (-2.45+0.25,-1.2+0.6);
\node[left] at (-2,-0.3) {$\flt$}; 
\draw[dashed] (-2.3-0.4,-2.2) to [out=70,in=180] (-0.5,-2.5) node[below] {\scriptsize $S_{t=1}: \Sigma$} to [out=0,in=110] (2.3+0.4,-2.2);
\draw (2.3+0.4,-2.2) to [out=-70,in=0] (1,-2) to [out=180,in=-110] (-2.3-0.4,-2.2); 
\draw[->] (-1.7-0.3,0.7+0.3) to [out=-110,in=70] (-2.3-0.8,-2.2+0.4);
\draw (-2.3-0.4,-2.2) to [out=135, in= -45] (-2.3-0.4-0.7,-2.2+0.7);
\node[left] at (-2.3-0.55,-2.2+1-0.2) {$f$}; 
\end{tikzpicture}
\end{center}
\caption{The family of surfaces $\{ S_t \}$.}
\label{fig 8}
\end{figure}

\noindent
Assume that the first parameterisation of $S_t$ is $(\uflt, \flt)$, then $\uflt$ satisfies the following first order nonlinear equation
\begin{align}
\label{eqn 3.3}
\partial_t \uflt = f \cdot \left[ 1- \left( b^i + \uvarepsilon^i \right) t f_i \right]^{-1} \cdot \left[ \uvarepsilon - \left( b^i + \uvarepsilon^i \right) \partial_i \uflt \right]
\end{align}
where $f_i$ is the partial derivative of $f$, and $\uvarepsilon, \uvarepsilon^i$ are given by the following formulae
\begin{align*}
\begin{aligned}
&
\underline{\varepsilon}^k = \underline{e}^k + \underline{\varepsilon} e^k,
\quad
\underline{\varepsilon} = \frac{ -|\underline{e}|^2}{(2\Omega^2 + e\cdot \underline{e}) + \sqrt{(2\Omega^2 + e\cdot \underline{e})^2 -|e|^2 |\underline{e}|^2}},
\\
&
|e|^2 = \slashg_{ij}e^ie^j,
\quad
|\underline{e}|^2 = \slashg_{ij} \underline{e}^i \underline{e}^j, 
\quad
e\cdot \underline{e} =\slashg_{ij} e^i \underline{e}^j,
\\
&
e^k =-2\Omega^2 \cdot tf_i \left(B^{-1}\right)_j^i \left(\slashg^{-1}\right)^{jk}, 
\quad
\underline{e}^k = -2\Omega^2 \cdot \uflt_i \left(B^{-1}\right)_j^i \left( \slashg^{-1} \right)^{jk},
\\
&
B_i^j= \delta_i^j - \flt_i b^j.
\end{aligned}
\end{align*}
We give the derivation of equation \eqref{eqn 3.3} in appendix \ref{appen eqn 3.3}. Note that $\uvarepsilon, \uvarepsilon^k, e^k, \ue^k, B_i^j$ above depend on $t$, but for the sake of brevity, we donot indicate the $t$-dependence in the symbols denoting these notations. 

In order to obtain the parameterisation function $\uf$ of $\Sigma$, we solve equation \eqref{eqn 3.3} for $\uflt$ with the initial value
\begin{align*}
\ufl{t=0} = \ufl{s=0},
\end{align*}
then $\uf$ is
\begin{align*}
\uf= \ufl{t=1}.
\end{align*}
Note that in the special case that $f\equiv s$ a constant function, equation \eqref{eqn 3.3} leads exactly to equation \eqref{eqn 3.1}.
\end{method II}

Comparing above two methods, method II is more direct. In method I, we solve equation \eqref{eqn 3.2} for the parameterisation $\uh$ of $\ucalH$ first, then restrict it to $\Sigma$ to obtain $\uflt$. In the following, we introduce a variant of method I and II, which derives an equation for $\uflt$ from equation \eqref{eqn 3.2}. We introduce the following lemma first.
\begin{lemma}\label{lem 3.1}
Let $\upsilon$ be a function on $\ucalH$ and satisfies the equation
\begin{align*}
\dpartial_s \upsilon = F(s, \upsilon, \dpartial_i \upsilon).
\end{align*}
Suppose that $\{S_t\}$ is a family of surfaces embedded in $\ucalH$, where $S_t$ is parameterised by $s= \flt(\vartheta)$ in $(s,\vartheta)$ coordinates on $\ucalH$. Define the function $u(t,\vartheta)$ by
\begin{align*}
u:(t, \vartheta) \mapsto (s,\vartheta) \mapsto \upsilon(s, \vartheta),
\quad
(s,\vartheta) = (\flt(\vartheta), \vartheta),
\quad
u(t,\vartheta) = \upsilon(\flt(\vartheta), \vartheta).
\end{align*}
Then
\begin{align*}
\partial_t u = \partial_t \flt \cdot \dpartial_s \upsilon,
\quad
\partial_i u = \partial_i \flt \cdot \dpartial_s \upsilon + \dpartial_i \upsilon,
\end{align*}
and $u$ satisfies the following equation
\begin{align*}
\partial_t u = \partial_t \flt \cdot F\big(\flt, u,\, \partial_i u - \partial_i \flt \cdot \left( \partial_t \flt \right)^{-1} \partial_t u \big).
\end{align*}
\end{lemma}
The proof of the above lemma is straightforward. Applying it to $\uh$ on $\ucalH$ and $\{S_t\}$ in method II, and noting that $\flt= t f$, we obtain an equation for the parameterisation function $\uflt(\vartheta) = \uh( \flt(\vartheta), \vartheta)$ of $S_t$
\begin{align*}
\partial_t \uflt 
&= 
f \cdot F\big( tf,\, \uflt,\, \uflt_i - t f^{-1} f_i \cdot \partial_t \uflt  \big)
\\
&=
f \cdot \left[
b^i  \big( \uflt_i - t f^{-1} f_i \cdot \partial_t \uflt \big)
\right.
\\
&\phantom{= f \cdot \Big[ \ }
\left.
+ \Omega^2 ( \slashg^{-1} )^{ij} \cdot \big( \uflt_i - t f^{-1} f_i \cdot \partial_t \uflt \big) \cdot \big( \uflt_j - t f^{-1} f_j \cdot \partial_t \uflt \big)
\right].
\end{align*}
The above equation is a first order nonlinear equation of $\uflt$. In this paper, we will not solve this equation to obtain $\uflt$. However, we will use the above idea to obtain the restriction of the differential of $\uh$, i.e. $\dpartial_i \uh$, on $\Sigma$ in subsection \ref{subsec 5.1}.

\subsection{Estimate of the parameterisation transformation}\label{subsec 3.3}
In previous subsection, we introduce two methods to obtain the first parameterisation $(\uf, f)$ from the second one $(\ufl{s=0} ,f)$. In the following, we will estimate this parameterisation transformation through method II.

Before considering the general case, we shall mention that the estimate in the special case $f\equiv s$ being a constant function is already obtained in \cite{L4}. We state the result here, which is essentially theorem 3.3 in \cite{L4}.
\begin{proposition}\label{pro 3.2}
Let $\ucalH$ be an incoming null hypersurface in $(M,g)$ which is parameterised by $\us = \uh(s, \vartheta)$ as its graph in the $\{\us, s, \vartheta \}$ coordinate system. Let $\ufl{s=0}$ be the parameterisation function of $\Sigma_0$, which the intersection of $\ucalH$ with $C_{s=0}$,
\begin{align*}
\ufl{s=0}(\vartheta) = \uh(0,\vartheta).
\end{align*}
Suppose that $\ufl{s=0}$ satisfies the following estimates\footnote{Subscript $o$ is short for oscillation, and $m$ is short for mean value, since $\udelta_o$ bounds the norm of the differential and $\udelta_m$ bounds the mean value.}
\begin{align*}
\leVe \dslashd \ufl{s=0} \riVe^{n+1,p}_{\circg} \leq \udelta_o r_0,
\quad
\leve \overline{\ufl{s=0}}^{\circg} \rive \leq \udelta_m r_0,
\end{align*}
where $n\geq 1,p>2$ or $n\geq 2, p>1$.

There exist a small positive constant $\delta$ depending on $n,p$, and constants $c_o, c_{m,m}, c_{m,o}$ also depending on $n,p$, such that if $\epsilon, \udelta_o, \udelta_m$ are suitably bounded that $\epsilon, \udelta_o, \epsilon \udelta_m \leq \delta$, then $\uh$ satisfies the following estimates
\begin{align*}
&
\leVe \dslashd \uh \riVe_{\circg}^{n+1,p} 
\leq
c_o \leVe \dslashd \ufl{s=0} \riVe_{\circg}^{n+1,p} 
\leq
c_o \udelta_o r_0,
\\
&
\leve \overline{\uh(s,\cdot)}^{\circg} - \overline{\ufl{s=0}}^{\circg} \rive 
\leq
 \frac{c_{m,m}}{r_0} \epsilon \big\vert \overline{\ufl{s=0}}^{\circg} \big\vert \leVe \dslashd \ufl{s=0} \riVe^{n+1,p}  +  \frac{c_{m,o}}{r_0} \left( \leVe \dslashd \ufl{s=0} \riVe^{n+1,p} \right)^2
 \\
 &\hspace{80pt}
\leq 
( c_{m,m} \epsilon \udelta_m \udelta_o + c_{m,o} \udelta_o^2 ) r_0,
\end{align*}
for all $s\in (-\kappa r_0, \kappa r_0)$.
\end{proposition}
This proposition is proved by estimating solutions of equations \eqref{eqn 3.1} \eqref{eqn 3.2}. We apply the Laplacian $\dcircDelta$ to equation \eqref{eqn 3.2} and derive the following equation\footnote{We always use the dot on the top to mean that the corresponding object is associated with the $\{ s,\vartheta \}$ coordinate system on $\ucalH$. For example, $\dpartial_s$ is the coordinate derivative in $\{ s,\vartheta \}$ coordinate system on $\ucalH$. Let $\Sigma_s$ be the coordinate surface of constant $s$ in $\{ s,\vartheta \}$ coordinate system on $\ucalH$, then
 $\dcircnabla$ is the covariant derivative of $(\Sigma_s, \circg)$, and $\dcircDelta$ is the Laplacian of $(\Sigma_s, \circg)$.}
\begin{align*}
\dpartial_s \dcircDelta \uh
=&
-b^i \dpartial_i (\dcircDelta \uh)
+ 2\Omega^2 {\slashg^{-1}}^{ij} \uh_j \dpartial_i (\dcircDelta\uh)
+ \text{lower order terms involving $\dslashd \uh, \dcircnabla^2 \uh$}.
\end{align*}
Then integrate the above equation and use Gronwall's inequality to obtain the estimates of $\uh$.

The general case is treated by a similar manner. Applying the Laplacian $\ddcircDelta$ of $\left( S_t, \circg \right)$ to equation \eqref{eqn 3.3} of the parameterisation function $\uflt$,\footnote{We shall use two dots on the top to mean that the corresponding object is defined on $S_t$. For example $\ddcircnabla$ denotes the covariant derivative on $(S_t, \circg)$, $\ddpartial_i$ denotes the coordinate derivative of $\vartheta$ on $S_t$.} we can derive an equation for $\ddcircDelta \uflt$. Then integrate this equation and use Gronwall's inequality to obtain the estimates of $\uflt$.

We rewrite equation \eqref{eqn 3.3} as
\begin{align}
\label{eqn 3.4}
\begin{aligned}
&
\partial_t \uflt 
= 
F\left( f,\ t b^i f_i,\ t e^i f_i,\ t \ue^i f_i,\ \uvarepsilon,\  b^i \, \uflt_i,\  e^i \, \uflt_i,\  \ue^i \, \uflt_i  \right),
\\
&
F= f \cdot [ 1-tb^i f_i - t\ue^i f_i - t e^i f_i \cdot \uvarepsilon ]^{-1} \cdot \left[ \uvarepsilon - b^i \, \uflt_i - \ue^i \, \uflt_i - e^i \, \uflt_i \cdot \uvarepsilon \right].
\end{aligned}
\end{align}
Then applying the Laplacian $\ddcircDelta$ to equation \eqref{eqn 3.4}, we obtain the following equation of $\ddcircDelta \uflt$
\begin{align}
\label{eqn 3.5}
\partial_t \big( \ddcircDelta \uflt \big)
=
\Xlt^i \ddpartial_i \big( \ddcircDelta \uflt \big) + \relt,
\end{align}
where
\begin{align*}
\Xlt^i
=
&
\partial_{t \ue^i f_i}F \cdot \left[-2\Omega^2 \left(B^{-1}\right)_j^i \left( \slashg^{-1} \right)^{jk} \cdot t f_k  \right]
+
\partial_{b^i \, \uflt_i} F \cdot b^i
\\
&
+
\partial_{\ue^i \, \uflt_i} F \cdot \left[ \ue^i -2\Omega^2 \left(B^{-1}\right)_j^i \left( \slashg^{-1} \right)^{jk} \, \uflt_k \right]
\\
&
+
\partial_{\uvarepsilon} F \cdot 
\partial_{\vert \ue \vert^2} \uvarepsilon \cdot \left[ 8\Omega^4 \left( B^{-1} \right)_j^i  \left( \slashg^{-1} \right)^{jk} \left(B^{-1}\right)_k^l \, \uflt_l \right]
\\
&
+
\partial_{\uvarepsilon} F \cdot \partial_{e\cdot \ue} \uvarepsilon \cdot \left[  4\Omega^4 \left(B^{-1}\right)_k^l \left( \slashg^{-1} \right)^{jk} \left( B^{-1} \right)_j^i \cdot t f_l   \right],
\end{align*}
and
\begin{align*}
\relt
&=
\partial_{f} F \cdot \ddcircDelta f 
+
\partial_{t f_i \ue^i} F \cdot \left\{ \ddcircDelta \left( t f_i \ue^i \right) - \left[-2\Omega^2 \left(B^{-1}\right)_j^i \left( \slashg^{-1} \right)^{jk} \cdot t f_k \right] \big(\ddcircDelta \uflt \big)_i  \right\}
\\
&\phantom{=}
+
\partial_{t e^i f_i} F \cdot \ddcircDelta \left( te^i f_i \right)
+
\partial_{t b^i f_i} F \cdot \ddcircDelta\left( t b^i f_i \right)
\\
&\phantom{=}
+ 
\partial_{\uvarepsilon} F \cdot \partial_{\vert \ue \vert^2} \uvarepsilon \cdot \left\{ \ddcircDelta \vert \ue \vert^2 - 8\Omega^4\ \uflt_l \left(B^{-1}\right)_k^l \left( \slashg^{-1} \right)^{jk} \left( B^{-1} \right)_j^i \big( \ddcircDelta \uflt \big)_i \right\}
\\
&\phantom{=}
+
\partial_{\uvarepsilon} F \cdot \partial_{e \cdot \ue} \uvarepsilon \cdot \left\{ \ddcircDelta \left( e \cdot \ue \right) -4\Omega^4\ t f_l \left(B^{-1}\right)_k^l \left( \slashg^{-1} \right)^{jk} \left( B^{-1} \right)_j^i \big( \ddcircDelta \uflt \big)_i  \right\}
\\
&\phantom{=}
+
\partial_{e^i \uflt_i} F \cdot \left\{ \ddcircDelta \big( e^i \, \uflt_i \big) - e^i  \big( \ddcircDelta\uflt \big)_i \right\}
+
\partial_{b^i \uflt_i} F \cdot \left\{ \ddcircDelta \big( b^i \, \uflt_i \big) - b^i  \big( \ddcircDelta\uflt \big)_i \right\}
\\
&\phantom{=}
+
\partial_{\ue^i \uflt_i }F \cdot \left\{ \ddcircDelta\big( \ue^i \, \uflt_i \big) - \left[ \ue^i -2\Omega^2 \left(B^{-1}\right)_j^i \left( \slashg^{-1} \right)^{jk} \uflt_k \right] \big( \ddcircDelta \uflt \big)_i \right\}
\\
&\phantom{=}
+
\partial^2_{\mathbf{a} \mathbf{b}} F \cdot \ddcircnabla_i \mathbf{a} \ddcircnabla^i \mathbf{b},
\\
\mathbf{a}, \mathbf{b}:
&\phantom{=}
f,\ t b^i f_i,\ t e^i f_i,\ t \ue^i f_i,\ \uvarepsilon,\  b^i \, \uflt_i,\  e^i \, \uflt_i,\  \ue^i \, \uflt_i.
\end{align*}
We shall elaborate on $\relt$. 
\begin{enumerate}[label=\alph*.]
\item Firstly in the formula of $\relt$, terms like $\ddcircnabla \vec{b}, \ddcircnabla \Omega, \ddcircnabla \slashg$ are given by
\begin{align*}
\ddcircnabla_i a = \circnabla_i a + \uflt_i \cdot \partial_{\us} a + t f_i \cdot \partial_s a,
\end{align*}
where $a$ is some background quantity like the metric components or their derivatives with respect to $\circnabla, \partial_s, \partial_{\us}$ as in definition \ref{def 2.3}. The higher order covariant derivatives by $\ddcircnabla$ are calculated successively by the above rule, for example the second order covariant derivative $\ddcircnabla^2_{ij} a$ is
\begin{align*}
\ddcircnabla_i \ddcircnabla_j a
=&
\circnabla^2_{ij} a + (tf)_i \cdot \partial_s \circnabla_j  a + \uflt_i \cdot \partial_{\us} \circnabla_j a 
\\
&
+ \ddcircnabla^2_{ij} f \cdot \partial_s a 
+ (tf)_j \cdot \circnabla_i \partial_s a 
+ (tf)_j \uflt_i   \cdot \partial_s \partial_{\us} a 
+ (tf)_i (tf)_j \cdot \partial_s^2 a
\\
&
+ \ddcircnabla^2_{ij} \uflt \cdot \partial_{\us} a
+ \uflt_j \cdot \circnabla_i \partial_{\us} a
+ (tf)_i \uflt_j \cdot \partial_s \partial_{\us} a
+ \uflt_i \uflt_j \cdot \partial_{\us}^2 a.
\end{align*}

\item Secondly, $\relt$ doesnot contain any top order derivatives of $\uflt$. The highest order derivative of $\uflt$ in $\relt$ is of second order.

\item Thirdly, the top order derivative of $f$ in $\relt$ is of third order, which comes from $\ddcircDelta$ applying to $f_i$ in $t b^i f_i$, $t e^i f_i$, $t \ue^i f_i$, $e^i \left( \uflt \right)^i$, $\vert e \vert^2$, $e \cdot \ue$ and $B_i^j$, $\left( B^{-1} \right)_i^j$.

\item Lastly, $\relt$ is a quadratic nonlinear term. Formally, if we write $\uflt$, $f_i$, $b^i$ as $\delta$, then the lowest degree terms of $\delta$ in $\relt$ are quadratic terms.

\end{enumerate}

Heuristically, we can make the following analogies for $F, \Xlt, \relt$
\begin{align*}
F: 
&\ 
r_0^{-2} f \left( \big\vert \ddslashd \uflt \big\vert_{\circg}^2 + \epsilon \vert \uflt \vert \cdot \big\vert \ddslashd \uflt \big\vert_{\circg} \right),
\\
\Xlt: 
&\
f \Big( b^i + r_0^{-2} \big( \circg^{-1} \big)^{ij} \uflt_i \Big) \ddpartial_i,
\\
\relt: 
&\
r_0^{-2}
\Big(  \vert f \vert + \big\vert \ddslashd f \big\vert_{\circg} + \big\vert \ddcircnabla^2 f \big\vert_{\circg}  \Big)
\cdot
\Big( \epsilon \vert \uflt \vert + \big\vert \ddslashd \uflt \big\vert_{\circg} + \big\vert \ddcircnabla^2 \uflt \big\vert_{\circg} \Big)^2
\\
&\
+ r_0^{-3}
\vert f \vert 
\cdot
\Big( \epsilon \vert \uflt \vert \cdot \big\vert \ddslashd \uflt \big\vert + \big\vert \ddslashd \uflt \big\vert^2_{\circg} \Big) 
\cdot 
\Big\vert t \ddcircnabla^3 f \Big\vert_{\circg}.
\end{align*}
Note in $\relt$, the top order derivative of $f$ is one order higher than $\uflt$, thus it results in the consequence that the regularity of $\uflt$ will be at least one order less than $f$ when integrating equations \eqref{eqn 3.4} \eqref{eqn 3.5}.

Now we state the estimate of the parameterisation transformation in the general case obtained from equations \eqref{eqn 3.4} \eqref{eqn 3.5}.
\begin{proposition}\label{pro 3.3}
Let $\Sigma$ be a spacelike surface embedded in $(M,g)$ with the second parameterisation $\left(\ufl{s=0}, f \right)$. Assume that the parameterisation functions $\ufl{s=0}, f$ satisfy the estimates
\begin{align*}
\leVe \dslashd \ufl{s=0} \riVe^{n,p} \leq \udelta_{o} r_0, 
\quad 
\big\vert \overline{\ufl{s=0}}^{\circg} \big\vert \leq \udelta_m r_0,
\quad
\Vert \ddslashd f \Vert^{n+1,p} \leq \delta_{o} r_0,
\quad
\overline{f}^{\circg} = \os,
\end{align*}
where $n\geq 2, p>2$ or $n\geq 3, p>1$.

Assume that the first parameterisation of $\Sigma$ is $(\uf,f)$. There exist a small positive constant $\delta$ depending on $n,p$, and constants $c_o, c_{m,m}, c_{m,o}$ also depending on $n,p$, such that if $\epsilon, \udelta_o, \udelta_m, \delta_o$ are suitably bounded that $\epsilon, \udelta_o, \epsilon \udelta_m, \delta_o \leq \delta$, then the parameterisation function $\uf$ satisfies the following estimates
\begin{align}
\label{eqn 3.6}
\begin{aligned}
&
\leVe \ddslashd \uf \riVe^{n,p} 
\leq
c_o \udelta_o r_0, 
\\
&
\big\vert \overline{\uf}^{\circg} - \overline{\ufl{s=0}}^{\circg} \big\vert 
\leq
\left( c_{m,m} \epsilon \udelta_m \udelta_o + c_{m,o} \udelta_o^2 \right) (|\os|/r_0+ \delta_o)  r_0.
\end{aligned}
\end{align}
\end{proposition}
We emphasis again that from equations \eqref{eqn 3.4} \eqref{eqn 3.5}, the regularity of $\uflt$ is one order less than $f$, even we increase the regularity of $\ufl{s=0}$ the same as $f$. We shall sketch the proof here. The rest of the proof details will be presented in appendix \ref{appen pro 3.3}.
\begin{proof}[Proof sketch]
The proposition is proved using bootstrap argument. We construct the family of surfaces $\{S_t\}$ as in method II. It will be shown that estimates \eqref{eqn 3.6} hold for every parameterisation function $\uflt, t\in [0,1]$.

We shall first assume the following inequalities for $\delta$ and $c_o, c_{m,m}, c_{m,o}$,
\begin{align*}
\delta \leq \frac{1}{2}, \quad \left( c_o + c_{m,m} + c_{m,o} \right) \delta \leq 1.
\end{align*}
Clearly there exist some choice of $\delta, c_o>1, c_{m,m}, c_{m,o}$ and a small neighbourhood interval of $t=0$ in which estimates \eqref{eqn 3.6} hold. Therefore we introduce the following bootstrap assumption. 
\begin{assum}
Estimates \eqref{eqn 3.6} hold for $\uflt$ in the closed interval $[0, t_a]$.
\end{assum}
The goal is to prove that the assumption is true on the interval $ [0,1]$. We shall show that for carefully chosen $\delta,c_o, c_{m,m}, c_{m,o}$, the inequalities in estimates \eqref{eqn 3.6} can be improved to strict inequalities at the end point $t=t_a$. Then this implies that the maximal interval where the assumption is valid is $[0,1]$.

In order to integrate equations \eqref{eqn 3.4} \eqref{eqn 3.5}, we need to estimate $F,\Xlt,\relt$. We introduce the notations $\ud_o, \ud_m$ to simplify some formulae in the proof,
\begin{align}
\ud_o = c_o \udelta_o, 
\quad 
\ud_m = \big[ 1+ c_{m,m} (|\os|/r_0+ \delta_o) \epsilon \udelta_o \big] \udelta_m + c_{m,o} (|\os|/r_0+ \delta_o) \udelta_o^2.
\label{eqn 3.7}
\end{align}

By the bootstrap assumption and bounds of metric components in definition \ref{def 2.3}, we can show that these terms satisfy the following estimates\footnote{We abuse the notation $c(n,p)$ to denote any constant depending on $n,p$. It is not necessary that different $c(n,p)$ denote the same constant in the proof, not even in the same formula. \label{footnote 6}}
\begin{align}
\begin{aligned}
&
\vert F \vert
\leq
c(n,p) (|\os|/r_0+ \delta_o) ( \ud_o +  \epsilon \ud_m ) \ud_o r_0,
\\
&
\Vert \Xlt \Vert^{n,p}
\leq
c(n,p) (|\os|/r_0+ \delta_o) ( \ud_o +  \epsilon \ud_m )
\leq
c(n,p),
\\
&
\Vert \relt \Vert^{n-1,p}
\leq
c(n,p) (|\os|/r_0+ \delta_o) ( \ud_o +  \epsilon \ud_m ) \ud_o r_0.
\end{aligned}
\label{eqn 3.8}
\end{align}
Therefore by Gronwall's inequality and the theory of Laplace equation on the sphere, we obtain that
\begin{align*}
&
\begin{aligned}
\leve \overline{\uflt}^{\circg} - \overline{\ufl{s=0}}^{\circg} \rive
\leq
&
\int_0^{t} \big| \overline{F}^{\circg} \big| \d t
\\
\leq 
&
c(n,p) (|\os|/r_0+ \delta_o) \left[ c_o \udelta_o +  \epsilon (1+ c_{m,m} \epsilon \udelta_o) \udelta_m + c_{m,o} \epsilon \udelta_o^2 \right] c_o \udelta_o r_0
\\
\leq 
&
c(n,p) (|\os|/r_0+ \delta_o) c_o (1 + c_{m,m} \epsilon \udelta_o) \cdot \epsilon \udelta_m \udelta_o r_0
\\
&
+
c(n,p) (|\os|/r_0+ \delta_o) c_o (c_o + c_{m,o} \epsilon \udelta_o ) \cdot \udelta_o^2 r_0 ,
\end{aligned}
\\
&
\begin{aligned}
\leVe \ddslashd \uflt \riVe^{n,p}
\leq
&
c(n,p) \leVe \ddcircDelta \uflt \riVe^{n-1,p}
\leq
c(n,p)\left\{ \leVe \dcircDelta \ufl{s=0} \riVe^{n-1,p} + \int_0^t \Vert \relt \Vert^{n-1,p} \d t \right\}
\\
\leq 
&
c(n,p) \udelta_0 r_0 
\\
&
+ c(n,p) (|\os|/r_0 + \delta_o) \left[ c_o \udelta_o +  \epsilon (1+ c_{m,m} \epsilon \udelta_o) \udelta_m + c_{m,o} \epsilon \udelta_o^2 \right] c_o \udelta_o r_0.
\end{aligned}
\end{align*}
Recall that in definition \ref{def 2.3}, we assume that $|\os|/r_0 \leq \kappa \leq 0.1$.
Then in order to close the bootstrap argument, we require $\delta, c_o, c_{m,m}, c_{m,o}$ satisfy the following inequalities
\begin{align*}
c(n,p) c_o (1 + c_{m,m} \delta^2) < c_{m,m},
\\
c(n,p) c_o (c_o + c_{m,o} \delta^2 ) < c_{m,o},
\\
c(n,p) + c(n,p) \left[ c_o \delta +  (1+ c_{m,m} \delta^2) \delta + c_{m,o} \delta^3 \right] c_o < c_o.
\end{align*}
In the above inequalities, it can be assumed that all constants $c(n,p)$ are the same. We choose $( c_o, c_{m,m}, c_{m,o} ) = (2 c(n,p), 3 c(n,p)^2, 5c(n,p)^3)$, then choose $\delta$ sufficiently small such that the above inequalities hold. Then the bootstrap argument is closed.
\end{proof}
\begin{remark}
We shall see in appendix \ref{appen pro 3.3} that the proof of proprosition \ref{pro 3.3} requires the $L^{\infty}$ bounds of the metric components $b, \Omega, \slashg$ and their derivatives up to the $(n+1)$-th order.
\end{remark}

We already explained why there is one order less in the regularity of $\uf$ than $f$. It comes from the $3$rd order derivatives of $f$ in $\relt$ of equation \eqref{eqn 3.5}. However we can actually improve the regularity of $\uf$ to the same order as $f$, by taking equation \eqref{eqn 3.2} of the parameterisation function $\uh$ of $\ucalH$ into account additionally. In fact, $\uf$ satisfies the following improved estimate.
\begin{proposition}\label{pro 3.5}
Under the same setting of proposition \ref{pro 3.3}, assume additionally that 
\begin{align*}
\leVe \dslashd \ufl{s=0} \riVe^{n+1,p} \leq \udelta_o r_o.
\end{align*}
Then there exist suitable $\delta$ and $c_{o'}$ depending on $n,p$ such that
\begin{align*}
\leVe \ddslashd \uf \riVe^{n+1,p} 
\leq 
c_{o'} \leVe \dslashd \ufl{s=0} \riVe^{n+1,p} 
\leq 
c_{o'} \udelta_o r_0.
\end{align*}
\end{proposition}
The improved estimate of $\uf$ is a corollary of later proposition \ref{pro 5.1} on the restriction of the differential of $\uh$ on $\Sigma$, thus we leave the proof after proving proposition \ref{pro 5.1}. Here we just remark that the above improved estimate of $\uf$ requires the $L^{\infty}$ bounds of the metric components up to the $(n+2)$-th order derivatives.

\section{Formula of the outgoing null expansion}\label{sec 4}
In this section, we derive the general formula for the outgoing null expansion of a spacelike surface $\Sigma$ in $(M,g)$. Then we treat the case in the Schwarzschild spacetime as an example. In order to clarify the structure of the outgoing null expansion, we shall present a decomposition of the formula into the first order main part and high order remainder part in the last subsection.

\subsection{A two-step procedure to calculate the outgoing null expansion}
Suppose that $\Sigma$ has the second parameterisation $(\ufl{s=0}, f)$. The outgoing null expansion of $\Sigma$ is obtained through a two-step procedure. Suppose that $\Sigma$ is embedded in an incoming null hypersurface $\ucalH$.
\begin{step i}
Let $\{\Sigma_s\}$ be the intersection of $\Sigma$ with $C_s$. Then $\{\Sigma_s\}$ foliates $\ucalH$. We calculate the metric components and structure coefficients associated with this foliation $\{\Sigma_s\}$ of $\ucalH$.
\end{step i}
\begin{step ii}
$\Sigma$ is embedded in $\ucalH$ as the graph $s = f(\vartheta)$ in $\{s, \vartheta\}$ coordinate system. Then we calculate its outgoing null expansion in terms of the parametrisation function $f$ and geometric quantities associated with $\{\Sigma_s\}$ in step i.
\end{step ii}

We carry out these two steps in the following. The essential formulae for this procedure are given in \cite{L4} and will be applied here without derivations.

\begin{step i}
$\Sigma_s$ has the second parameterisation $(\ufl{s=0}, s)$. Its first parameterisation is $(\ufls, s)$ where $\ufls(\vartheta) = \uh(s, \vartheta)$ as stated in subsection \ref{subsec 3.1} Method I. The coordinate vectors of $\{s, \vartheta\}$ coordinate system on $\ucalH$ are
\begin{align*}
\dpartial_s = \partial_s + \dpartial_s \uh\cdot \partial_{\us},
\quad
\dpartial_i = \partial_i + \dpartial_i \uh\cdot \partial_{\us}.
\end{align*}
We use $\cdot$ on the top to indicate the corresponding notation being associated with $\Sigma_s$ or the $\{s, \vartheta\}$ coordinate system on $\ucalH$.

Introduce the conjugate null frame $\{\duL, \dL'\}$ on $\Sigma_s$. Here for the sake of brevity, we donot emphasis the dependence of $s$ in the symbols denoting the frame vectors. 
\begin{align*}
\left\{
\begin{aligned}
&
\dL'=L',
\\
&
\duL=\uL+\uvarepsilon L + \uvarepsilon^i \partial_i,
\end{aligned}
\right.
\end{align*}
where
\begin{align*}
\uvarepsilon=-\Omega^2 (\slashg^{-1} )^{ij} \uh_i \uh_j = -\Omega^2 \vert \dslashd \uh \vert_{\slashg}^2,
\quad
\uvarepsilon^i = -2\Omega^2 ( \slashg^{-1} )^{ij} \uh_j.
\end{align*}
The shifting vector $\db$ between $\dpartial_s$ and $\duL$ is given by
\begin{align*}
\duL = \dpartial_s + \db^i \dpartial_i
\quad
\Rightarrow
\quad
\db^i = b^i - 2 \Omega^2 ( \slashg^{-1} )^{ij} \uh_j.
\end{align*}

Let $\dslashg$ be the intrinsic metric on $\Sigma_s$, then 
\begin{align*}
\dslashg_{ij} = \slashg_{ij},
\quad
( \dslashg^{-1} )^{ij} = ( \slashg^{-1} )^{ij}.
\end{align*}
The degenerated metric on $\ucalH$ in $\{s, \vartheta\}$ coordinate system is
\begin{align*}
g\vert_{\ucalH} =  \dslashg_{ij} \big( \d \vartheta^i - \db^i \d s  \big) \otimes \big( \d \vartheta^j - \db^j \d s \big).
\end{align*}

The structure coefficients on $\Sigma_s$ associated with $\{\duL, \dL' \}$ are given by the following formulae:
\begin{align}
\nonumber
\dchi'_{ij} =& \chi'_{ij},  \quad \dtr \dchi' = \tr \chi',
\\
\nonumber
\duchi_{ij}
=&
\uchi_{ij} -\Omega^2 \vert \dslashd \uh \vert^2_{\slashg} \chi_{ij}  -2\Omega^2 \slashnabla_{ij}^2\uh -4\Omega^2 \sym \left\{\ueta \otimes \dslashd \uh \right\}_{ij} 
\\
\nonumber
&
- 4\omega\Omega^2 (\dslashd \uh \otimes \dslashd \uh )_{ij} 
+ 4\Omega^2 \sym \left\{\dslashd\uh \otimes ( \chi \cdot \dslashd \uh ) \right\}_{ij},
\\
\nonumber
\dtr \duchi =& \left( \dslashg^{-1} \right)^{ij} \duchi_{ij}
=  \tr \uchi - 2\Omega^2 \slashDelta \uh - \Omega^2 \vert \dslashd \uh \vert_{\slashg}^2 \tr \chi -4\Omega^2 \ueta \cdot \dslashd \uh 
\\
&\hspace{14ex}
- 4 \Omega^2 \omega \vert \dslashd \uh \vert^2_{\slashg} + 4 \Omega^2 \chi ( \slashnabla \uh, \slashnabla \uh),
\label{eqn 4.1}
\\
\nonumber
\deta_i =& \eta_i + ( \chi \cdot \dslashd \uh )_i,
\\
\nonumber
\duomega 
=& 
\uomega - 2\Omega^2 \eta \cdot \dslashd \uh - \Omega^2 \chi ( \slashnabla \uh, \slashnabla \uh).
\end{align}
In the above formulae, we use $\cdot$ to denote the inner product with respect to $\slashg$, and use $\dtr$ to denote the trace with respect to $\dslashg= \slashg$. $\dslashd$ is differential operator on $\Sigma_s$. $\slashnabla$ in $\slashnabla \uh, \slashnabla^2_{ij} \uh$ is the pull back of the covariant derivative $\slashnabla$ of $(\Sigma_{\us,s}, \slashg)$ to $\Sigma_s$.\footnote{Here we abuse the notation $\slashnabla$ to denote both the covariant derivative of $(\Sigma_{\us,s}, \slashg)$ and its pull back to $\Sigma_s$. Which meaning $\slashnabla$ represents in a concrete formula depends on where the differentiated function, vector field or tensor field is defined. For example, if a vector field $V$ is defined on $\Sigma_s$, then $\slashnabla$ in $\slashnabla V$ is the pull back of the covariant derivative of $(\Sigma_{\us,s}, \slashg)$ on $\Sigma_s$. If $\slashnabla$ can be interpreted in both way in a formula, we will state the precise meaning of $\slashnabla$ in that formula to avoid ambiguity.  \label{footnote 7}}
$\slashDelta$ in $\slashDelta \uh$ is the operator $\left( \slashg^{-1} \right)^{ij} \slashnabla^2_{ij}$.

The precise meaning of the pull back $\slashnabla$ on $\Sigma_s$ is as follows: let $\slashGamma_{ij}^k$ be the Christoffel symbol of the covariant connection $\slashnabla$ of $(\Sigma_{\us,s}, \slashg)$, then
\begin{align*}
&
\text{$\phi$: a function on $\Sigma_s$}
&&
\slashnabla_i \phi = \dpartial_i \phi, \quad \slashnabla^i \phi= ( \slashg^{-1} )^{ij} \dpartial_j \phi,
\\
&
\text{$V$: a vector field on $\Sigma_s$}
&&
\slashnabla_i V^k =  \dpartial_i V^k + \slashGamma_{ij}^k V^j,
\\
&
\text{$T$: a tensor field on $\Sigma_s$}
&&
\slashnabla_{i} T_{i_1 \cdots i_k}^{j_1 \cdots j_l}
=
\dpartial_i T_{i_1 \cdots i_k}^{j_1 \cdots j_l} 
-  \slashGamma_{i i_m}^{r}  T_{i_1 \cdots \underset{\hat{i_m}}{r}\cdots i_k}^{j_1 \cdots j_l}
+ \slashGamma_{i s}^{j_n} T_{i_1 \cdots i_k}^{j_1 \cdots \overset{\hat{j_n}}{s}\cdots  j_l}.
\end{align*}
\end{step i}

Before proceeding with step ii, we introduce the covariant derivative of $(\Sigma_s, \dslashg)$ denoted by $\dslashnabla$. Let $\dslashGamma_{ij}^k$ be the Christoffel symbol of $\dslashnabla$. It is given by the following formula,
\begin{align*}
\dslashGamma_{ij}^k = \slashGamma_{ij}^k + ( \slashg^{-1} )^{kl} \big( \dpartial_i \uh \cdot \chi_{jl} + \dpartial_j \uh \cdot \chi_{il} - \dpartial_l \uh \cdot \chi_{ij}  \big).
\end{align*}
We introduce the tensor $\triangle_{ij}^k$ to denote the difference of $\dslashGamma_{ij}^k$ with $\slashGamma_{ij}^k$,
\begin{align}\label{eqn 4.2}
\triangle_{ij}^k = ( \slashg^{-1} )^{kl} \big( \dpartial_i \uh \cdot \chi_{jl} + \dpartial_j \uh \cdot \chi_{il} - \dpartial_l \uh \cdot \chi_{ij}  \big).
\end{align}
$\dslashnabla$ is the covariant derivative of $(\Sigma_s, \dslashg)$, and we also use it to denote the pull back of $\dslashnabla$ to $\Sigma$. The precise interpretation of $\dslashnabla$ shall be understood in the context.\footnote{This is similar to the pull back of $\slashnabla$. See footnote \ref{footnote 7}. Which meaning $\dslashnabla$ is interpreted as depends on where the differentiated function or field is defined. For example, if $V$ is a vector field defined on $\Sigma$, then $\dslashnabla$ in $\dslashnabla V$ should be interpreted as the pull back of the covariant derivative of $(\Sigma_s, \dslashg)$ on $\Sigma$. If it can be interpreted in both way, we will point out the precise meaning of $\slashnabla$ to avoid ambiguity.} 
Let $\ddpartial_i$ denote the partial derivative in the coordinate system $\{\vartheta\}$ on $\Sigma$, then we have the following formulae for the pull back of $\dslashnabla$ on $\Sigma$
\begin{align*}
&
\text{$\phi$: a function on $\Sigma$}
&&
\dslashnabla_i \phi = \ddpartial_i \phi, \quad \dslashnabla^i \phi= ( \dslashg^{-1} )^{ij} \ddpartial_j \phi = ( \slashg^{-1} )^{ij} \ddpartial_j \phi ,
\\
&
\text{$V$: a vector field on $\Sigma$}
&&
\dslashnabla_i V^k =  \ddpartial_i V^k + \dslashGamma_{ij}^k V^j,
\\
&
\text{$T$: a tensor field on $\Sigma$}
&&
\dslashnabla_{i} T_{i_1 \cdots i_k}^{j_1 \cdots j_l}
=
\ddpartial_i T_{i_1 \cdots i_k}^{j_1 \cdots j_l} 
-  \dslashGamma_{i i_m}^{r}  T_{i_1 \cdots \underset{\hat{i_m}}{r}\cdots i_k}^{j_1 \cdots j_l}
+ \dslashGamma_{i s}^{j_n} T_{i_1 \cdots i_k}^{j_1 \cdots \overset{\hat{j_n}}{s}\cdots  j_l}.
\end{align*}
Now we return to proceed on step ii.

\begin{step ii}
$\Sigma$ is the graph of $s= f(\vartheta)$ in the $\{s, \vartheta\}$ coordinate system on $\ucalH$. The tangential frame vector of $\Sigma$ is given by
\begin{align*}
\ddpartial_i = \dpartial_i + f_i \cdot \dpartial_s = \dB_i^j \dpartial_i + f_i \duL,
\quad
\dB_i^j = \delta_i^j - f_i \cdot \db^j.
\end{align*}
We use $\ddot{a}$ to indicate the corresponding notation being associated with $\Sigma$.

Let $\ddslashg$ be the intrinsic metric on $\Sigma$,
\begin{align*}
\ddslashg_{ij} 
=
\dB_i^k \dB_i^l \dslashg_{kl}
=
\slashg_{ij} - \big( \slashg_{ik} \db^i f_j + \slashg_{jl} \db^l f_i \big) + f_i f_j \slashg_{kl} \db^k \db^l.
\end{align*}

Introduce the conjugate null frame $\{ \dduL, \ddL' \}$ on $\Sigma$
\begin{align*}
\left\{
\begin{aligned}
&
\ddL' = \dL' + \dvarepsilon' \duL + \dvarepsilon'^i \dpartial_i,
\\
&
\dduL = \duL,
\end{aligned}
\right.
\end{align*}
where
\begin{align*}
\dvarepsilon'^i
=&
-2 ( \dslashg^{-1} )^{ik} \big( \dB^{-1} \big)_k^j f_j,
\\
\dvarepsilon'
=&
-\vert \ddslashd f \vert_{\ddslashg}^2
=
- ( \ddslashg^{-1} )^{ij} f_i f_j.
\end{align*}
$\ddslashd$ is the differential operator on $\Sigma$.

We concern about the outgoing null expansion in the direction of $\ddL'$. We have the following formulae of it
\begin{align}
\begin{aligned}
\ddchi'_{ij} 
=&
\dchi'_{ij} + \dvarepsilon' \duchi_{ij} + ( \db \circdot \vec{\dvarepsilon}' -2 ) \dslashnabla^2_{ij} f
\\
&
+2\sym \left\{  
\ddslashd f \otimes \left[ \dslashnabla \db \circdot \vec{\dvarepsilon}' - \duchi (\vec{\dvarepsilon}' ) -\dvarepsilon' \duchi (\db ) - \dchi' (\db ) -2 \deta \right]
\right\}_{ij}
\\
&
+ \left[
2\duchi (\db,\vec{\dvarepsilon}' ) +\dvarepsilon' \duchi (\db,\db ) +\dchi' (\db,\db ) +4\deta (\db ) 
\right.
\\
&
\hskip35ex
\left.
-\dslashnabla_{\db} \db \circdot \vec{\dvarepsilon}' -\dpartial_s \db \circdot \vec{\dvarepsilon}' -4\duomega
\right] f_i f_j,
\\
\ddtr \ddchi' =& ( \ddslashg^{-1} )^{ij} \ddchi'_{ij}.
\end{aligned}
\label{eqn 4.3}
\end{align}
In the above formulae, we use $\circdot$ to denote the inner product with respect to $\dslashg$, for example
\begin{align*}
\big( \dslashnabla \db \circdot \vec{\dvarepsilon}' \big)_i = \dslashg_{kl} \cdot \dslashnabla_i \db^k \cdot \dvarepsilon'^l.
\end{align*}
Numerically, $\circdot$ is the same as the inner product with respect to $\slashg$, thus we will simply write $\cdot$ for $\circdot$ later. We use $\ddtr$ to denote the trace with respect to $\ddslashg$.
\end{step ii}
We finish the description of the two-step procedure to obtain the outgoing null expansion of $\Sigma$.

\subsection{Formula in Schwarzschild spacetime}\label{subsec 4.2}
In the following, we demonstrate the procedure in the Schwarzschild spacetime. Suppose that $\Sigma$ has the second parameterisation $( \ufl{s=0}, f )$ in the Schwarzschild spacetime. We calculate its outgoing null expansion of $\Sigma$ in the following.
\begin{step i}
$\ucalH$ is parametrised by $\uh$ as its graph of $\us = \uh(s,\vartheta)$. The coordinate frame vectors of $\{s, \vartheta\}$ coordinate system on $\ucalH$ are
\begin{align*}
\dpartial_s = \partial_s + \uh_s \partial_{\us},
\quad
\dpartial_i = \partial_i + \uh_i \partial_{\us}.
\end{align*}
The conjugate null frame $\{ \duL, \dL' \}$ of $\Sigma_s$ is
\begin{align*}
\left\{
\begin{aligned}
&
\dL'=L',
\\
&
\duL=\uL+\uvarepsilon L + \uvarepsilon^i \partial_i,
\end{aligned}
\right.
\end{align*}
where
\begin{align*}
\uvarepsilon=-r^{-2} \Omega_S^2  (\circg^{-1} )^{ij} \uh_i \uh_j = -r^{-2} \Omega_S^2 \vert \dslashd \uh \vert_{\circg}^2,
\quad
\uvarepsilon^i = -2 r^{-2} \Omega_S^2 ( \circg^{-1} )^{ij} \uh_j.
\end{align*}
The shifting vector $\db$ between $\dpartial_s$ and $\duL$ is
\begin{align*}
\duL = \dpartial_s + \db^i \dpartial_i
\quad
\Rightarrow
\quad
\db^i = - 2 r^{-2} \Omega_S^2 ( \circg^{-1} )^{ij} \uh_j.
\end{align*}
The intrinsic metric $\dslashg$ on $\Sigma_s$ is
\begin{align*}
\dslashg_{ij} = r^2 \circg_{ij},
\quad
( \dslashg^{-1} )^{ij} = r^{-2} ( \circg^{-1} )^{ij}.
\end{align*}
The degenerated metric on $\ucalH$ in $\{s, \vartheta\}$ coordinate system is
\begin{align*}
g_S\vert_{\ucalH} =  r^2 \circg_{ij} \big( \d \vartheta^i - \db^i \d s  \big) \otimes \big( \d \vartheta^j - \db^j \d s \big).
\end{align*}

The structure coefficients on $\Sigma_s$ with respect to $\{\duL, \dL' \}$ are given by the following formulae:
\begin{align*}
\dtr \dchi' =& \tr \chi_S', \quad \widehat{\dchi'} =0,
\\
\duchi_{ij}
=&
\uchi_{ij} 
-2\Omega_S^2 \circnabla_{ij}^2\uh
+ \Omega_S^2 \, \tr \chi_S \vert \dslashd \uh \vert^2_{\circg} \cdot \circg_{ij}  
- 4\omega_S\, \Omega_S^2 (\dslashd \uh \otimes \dslashd \uh )_{ij},
\\
\dtr \duchi 
=&
\tr \uchi_S 
- 2r^{-2} \Omega_S^2 \circDelta \uh 
+ r^{-2} \Omega_S^2  \tr \chi_S  \vert \dslashd \uh \vert_{\circg}^2 
- 4 r^{-2} \Omega_S^2 \omega_S \vert \dslashd \uh \vert^2_{\circg},
\\
\deta_i =& \frac{1}{2} \tr \chi_S \uh_i,
\\
\duomega 
=& 
\uomega_S - \frac{1}{2} r^{-2} \Omega_S^2 \tr \chi_S \vert \dslashd \uh \vert^2_{\circg}.
\end{align*}
\end{step i}

\begin{step ii}
$\Sigma$ is the graph $s= f(\vartheta)$ in the $\{s, \vartheta\}$ coordinate system on $\ucalH$. The tangential frame vector of $\Sigma$ is
\begin{align*}
\ddpartial_i = \dpartial_i + f_i \cdot \dpartial_s = \dB_i^j \dpartial_i + f_i \duL,
\quad
\dB_i^j = \delta_i^j - f_i \cdot \db^j.
\end{align*}
The intrinsic metric $\ddslashg$ on $\Sigma$ is
\begin{align*}
\ddslashg_{ij} 
=
\dB_i^k \dB_i^l \dslashg_{kl}
=
r^2 \circg_{ij} + 2 \Omega_S^2 (   \uh_i f_j + \uh_j f_i ) + 4 r^{-2} \Omega_S^4 \vert \dslashd \uh \vert^2_{\circg} \cdot f_i f_j .
\end{align*}
The conjugate null frame $\{ \dduL, \ddL' \}$ on $\Sigma$ is
\begin{align*}
\left\{
\begin{aligned}
&
\ddL' = \dL' + \dvarepsilon' \duL + \dvarepsilon'^i \dpartial_i,
\\
&
\dduL = \duL,
\end{aligned}
\right.
\end{align*}
where
\begin{align*}
\dvarepsilon'^i
=&
-2 r^{-2} ( \circg^{-1} )^{ik} \big( \dB^{-1} \big)_k^j f_j,
\\
\dvarepsilon'
=&
-\vert \ddslashd f \vert_{\ddslashg}^2
=
-( \ddslashg^{-1} )^{ij} f_i f_j.
\end{align*}
The outgoing null expansion in the direction of $\ddL'$ is
\begin{align}
\begin{aligned}
\ddchi'_{ij} 
=&
\frac{1}{2} r^2 \tr \chi'_S \cdot \circg_{ij}  -2 \circnabla^2_{ij} f + \cdots
\\
\ddtr \ddchi' 
=&
\tr \chi'_S  -2 r^{-2} \circDelta f + \cdots
\end{aligned}
\label{eqn 4.4}
\end{align}
Here we donot write all terms in the formulae of $\ddchi'_{ij}$ and $\ddtr \ddchi'$. We will see in the next subsection that the terms written explicitly here are the main parts of $\ddchi'_{ij}$ and $\ddtr \ddchi'$.
\end{step ii}

\subsection{Decomposition of the outgoing null expansion}\label{subsec 4.3}
We introduce the following decomposition of the outgoing null expansion of $\Sigma$ into first order main part and high order remainder part. 

Denote the first order main part of $\ddtr \ddchi'$ by $\lo{\ddtr \ddchi'}$ and the high order remainder part by $\hi{\ddtr \ddchi'}$,
\begin{align*}
&
\lo{\ddtr \ddchi'}
=
\tr \chi'_S|_{\Sigma} - 2 \left( r_S|_{\Sigma}\right)^{-2} \circDelta f,
\footnotemark
\\
&
\hi{\ddtr \ddchi'} = \ddtr \ddchi' - \lo{\ddtr \ddchi'}.
\end{align*}
\footnotetext{We view $\tr \chi'_S, r_S$ in the Schwarzschild spacetime as functions on $M$, then their restrictions to $\Sigma$ are simply their values on $\Sigma$.}

We give the detailed formula of $\hi{\ddtr \ddchi'}$ in the following. First introduce the following decompositions
\begin{align*}
\ddslashg:
&
\left\{
\begin{aligned}
\lo{\ddslashg} 
&=
\left( r_S|_{\Sigma} \right)^2 \circg, 
\\
\hi{\ddslashg_{ij}}
&=
\hir{1}{\ddslashg_{ij}}
+
\hir{2}{\ddslashg_{ij}},
\
\left\{
\begin{aligned}
\hir{1}{\ddslashg_{ij}} 
&=
\slashg_{ij} - ( r_S|_{\Sigma} )^2 \circg_{ij},
\\
\hir{2}{\ddslashg_{ij}} 
&=
- \big( \slashg_{ik} \db^i f_j + \slashg_{jl} \db^l f_i \big) + f_i f_j \slashg_{kl} \db^k \db^l,
\end{aligned}
\right.
\end{aligned}
\right.
\\
\ddslashg^{-1}:
&
\left\{
\begin{aligned}
\lo{\ddslashg^{-1}} 
&=
( r_S|_{\Sigma} )^{-2}\cdot  \circg^{-1} , 
\\
\hi{(\ddslashg^{-1})^{ij}}
&=
\hir{1}{(\ddslashg^{-1})^{ij}} 
+
\hir{2}{(\ddslashg^{-1})^{ij}} ,
\\
\hir{a}{(\ddslashg^{-1})^{ij}} 
&
=
\left\{
\begin{aligned}
&
- \left( r_S|_{\Sigma} \right)^{-2} \big( \circg^{-1} \big)^{ik} \cdot \hir{a}{ \ddslashg_{kl} }  \cdot \left(\ddslashg^{-1}\right)^{lj},
\\
&
- \left(\ddslashg^{-1}\right)^{ik}\cdot \hir{a}{ \ddslashg_{kl} }  \cdot \left( r_S|_{\Sigma} \right)^{-2} \big( \circg^{-1} \big)^{lj},
\end{aligned}
\right.
\quad
a=1,2.
\end{aligned}
\right.
\end{align*}
and
\begin{align*}
\\
\ddchi':
&
\left\{
\begin{aligned}
\lo{\ddchi'}
&=
\chi_S'|_{\Sigma}- 2 \circnabla^2 f,
\\
\hi{\ddchi'_{ij}}
&=
\hir{1}{\ddchi'_{ij}}
+
\hir{2}{\ddchi'_{ij}},
\\
\hir{1}{\ddchi'_{ij}}
&=
\chi'_{ij} - (\chi'_S|_{\Sigma})_{ij},
\\
\hir{2}{\ddchi'_{ij}}
&=
- 2 \big( \dslashnabla^2_{ij} f - \circnabla^2_{ij} f \big)
+ \db \circdot \vec{\dvarepsilon}' \dslashnabla^2_{ij} f
\\
&\phantom{=}
+2\sym \left\{  
\ddslashd f \otimes \left[ \dslashnabla \db \circdot \vec{\dvarepsilon}' - \duchi\big(\vec{\dvarepsilon}'\big) -\dvarepsilon' \duchi\big(\db\big) - \dchi'\big(\db\big) -2 \deta \right]
\right\}_{ij}
\nonumber
\\
&\phantom{=}
+ \left[
2\duchi (\db,\vec{\dvarepsilon}' ) +\dvarepsilon' \duchi (\db,\db ) +\dchi' (\db,\db ) +4\deta (\db ) 
\right.
\nonumber
\\
&
\phantom{= + \left[ 2\duchi(\db,\vec{\dvarepsilon}' ) +\dvarepsilon' \duchi (\db,\db ) \right.}
\left.
-\dslashnabla_{\db} \db \circdot \vec{\dvarepsilon}' -\dpartial_s \db \circdot \vec{\dvarepsilon}' -4\duomega
\right] f_i f_j,
\end{aligned}
\right.
\end{align*}
then the high order remainder part $\hi{\ddtr \ddchi'}$ is given by
\begin{align*}
&
\hi{\ddtr \ddchi'}
= 
( \ddslashg^{-1} )^{ij} \cdot \hi{\ddchi'_{ij}} 
+
\hi{ ( \ddslashg^{-1} )^{ij}} \cdot \ddchi'_{ij},
\\
&
\hir{a}{\ddtr \ddchi'}
= 
( \ddslashg^{-1} )^{ij} \cdot \hir{a}{\ddchi'_{ij}} 
+
\hir{a}{ ( \ddslashg^{-1} )^{ij}} \cdot \ddchi'_{ij},
\quad
a=1,2.
\end{align*}
Note in $\hi{\ddchi'}$, the formula contains the term $\dslashnabla^2_{ij} f - \circnabla^2_{ij} f$, which involves the difference between the Christoffel symbols $\dslashGamma$ and $\circGamma$. We give the formula of $\dslashGamma-\circGamma$ here. Introduce $\circtriangle = \slashGamma - \circGamma$
\begin{align*}
\circtriangle_{ij}^k
=
\frac{1}{2} ( \slashg^{-1} )^{kl} \big( \circnabla_i \slashg_{jl} + \circnabla_j \slashg_{il} - \circnabla_l \slashg_{ij} \big),
\end{align*}
thus $\dslashGamma - \circGamma$ is the sum of $\triangle$ in formula \eqref{eqn 4.2} and $\circtriangle$
\begin{align*}
\dslashGamma_{ij}^k - \circGamma_{ij}^k
=
\triangle_{ij}^k + \circtriangle_{ij}^k
=&
( \slashg^{-1} )^{kl} \big( \dpartial_i \uh \cdot \chi_{jl} + \dpartial_j \uh \cdot \chi_{il} - \dpartial_l \uh \cdot \chi_{ij}  \big)
\\
&
+
\frac{1}{2} ( \slashg^{-1} )^{kl} \big( \circnabla_i \slashg_{jl} + \circnabla_j \slashg_{il} - \circnabla_l \slashg_{ij} \big).
\end{align*}
Therefore the formula of $\dslashnabla^2 f - \circnabla^2 f$ is
\begin{align*}
\dslashnabla^2_{ij} f - \circnabla^2_{ij} f
= 
- \triangle_{ij}^k \cdot f_k - \circtriangle_{ij}^k \cdot f_k.
\end{align*}

We explain the motivation behind the above decompositions. Heuristically, assume that the differentials of the parametrisation functions $\ddslashd f$, $\dslashd \ufl{s=0} $, $\ddslashd \uf$ are of the size $\delta$. Furthermore assume that $\epsilon\cdot \ufl{s=0}, \epsilon \uf$ are also of size $\delta$. Then we sort the terms of sizes $\epsilon^2, \epsilon\delta, \delta^2$ and higher orders into the high order remainder part, and make the first order main part as simple as possible. Most terms in the high order remainder parts fit the above scheme, except terms
\begin{align*}
\hir{1}{\ddslashg}: \slashg_{ij} - \left( r_S|_{\Sigma} \right)^2 \circg_{ij} ,
\quad
\hir{1}{\ddchi'}: \chi'_{ij} - ( \chi'_S|_{\Sigma} )_{ij},
\end{align*}
and the terms $\hir{1}{\ddslashg^{-1}},\,\hir{1}{\ddtr \ddchi'}$ inherited from the above two terms. However comparing them with the first order main part, these two terms are still higher order small terms in the following sense:
\begin{align*}
\slashg_{ij} - \left( r_S|_{\Sigma} \right)^2 \circg_{ij} \lesssim \underbrace{ \epsilon r_0^2\ \lesssim\  r_0^2 }_{\text{smaller by the order of $\epsilon$}} \approx \left( r_S|_{\Sigma} \right)^2 \circg_{ij} = \lo{\ddslashg_{ij}},
\end{align*}
and for $f\sim \kappa r_0$
\begin{align*}
\chi'_{ij} - ( \chi'_S|_{\Sigma} )_{ij} \lesssim \underbrace{ \epsilon r_0 \lesssim \kappa r_0}_{\text{smaller by the order of $\epsilon$}} \approx ( \chi'_S|_{\Sigma} )_{ij}.
\end{align*}
There is no a priori reason that we must use this kind of decomposition instead of using another one such like
\begin{align*}
\lo{\ddslashg} = \slashg_{ij}, \quad \lo{\ddchi'} = \chi' - 2 \circnabla^2 f.
\end{align*}
However we shall see later that when we consider the perturbation and the linearised perturbations of $\ddtr \ddchi'$, the decomposition choosed here is convenient to work with.

\section{Estimate of the outgoing null expansion}\label{sec 5}
In this section, we shall estimate the outgoing null expansion of $\Sigma$ in terms of the bounds of its second parameterisation $(\ufl{s=0}, f)$. Given these bounds, the estimate of the first parameterisation function $\uf$ is obtained by proposition \ref{pro 3.3}. However we will see that these estimates of $\ufl{s=0}, f, \uf$ arenot sufficient to estimate the outgoing null expansion.

Recall the parameterisation function $\uh(s,\vartheta)$ of $\ucalH$ introduced in subsection \ref{subsec 3.2} method I. The estimates of $(\dslashd \uh)|_{\Sigma}$ and $(\circnabla^2 \uh)|_{\Sigma}$ are essential for estimating the outgoing null expansion. These estimates are obtained in propositions \ref{pro 5.1}, \ref{pro 5.3}.

Then with the estimates of $f, \uf$, $(\dslashd \uh)|_{\Sigma}$, $(\circnabla^2 \uh)|_{\Sigma}$, we can estimate the outgoing null expansion. The result is given in proposition \ref{pro 5.4}.

\subsection{Estimate of differential of parameterisation function $\uh$}\label{subsec 5.1}
The formula of the high order remainder term $\hi{\ddtr \ddchi'}$ involves $\db$, $\dvarepsilon'$, $\vec{\dvarepsilon}'$, $\duchi$, $\deta$, $\duomega$, $\dslashnabla \db$, $\dslashnabla^2 f$, $\dpartial_s \db$. These terms involve $\dslashd \uh$ and $\circnabla^2 \uh$. Therefore in order to obtain the estimate of $\ddtr \ddchi'$, it is necessary to obtain estimates of $\dslashd \uh$ and $\circnabla^2 \uh$ on $\Sigma$ first. We discuss the estimate of the differential $\dslashd \uh$ on $\Sigma$ in this subsection, and leave the estimate of the Hessian $\circnabla^2 \uh$ for the next subsection.

Before proceeding with the estimates, it is necessary to clarify two notations: $(\dslashd \uh)|_{\Sigma}$ and $\ddslashd \big(\uh|_{\Sigma}\big)$ or simply $\ddslashd \uh$. Recall that the partial derivative $\dpartial_i$ and the differential $\dslashd$ are associated with $\{s,\vartheta\}$ coordinate system on $\ucalH$, while the partial derivative $\ddpartial_i$ and the differential $\ddslashd$ are associated with $\vartheta$ coordinate system on $\Sigma$. In fact, we have
\begin{align*}
\ddpartial_i = \dpartial_i + f_i  \cdot \dpartial_s,
\end{align*}
as introduced in subsection \ref{subsec 4.2} step ii. Comparing the components of $(\dslashd \uh)|_{\Sigma}$ and $\ddslashd \big(\uh|_{\Sigma}\big)$,
\begin{align*}
(\dslashd \uh)_i|_{\Sigma} = (\dpartial_i \uh)|_{\Sigma},
\quad
(\ddslashd (\uh|_{\Sigma} ))_i 
=
\ddpartial_i (\uh|_{\Sigma})
=
(\dpartial_i \uh)|_{\Sigma} + f_i \cdot (\dpartial_s \uh)|_{\Sigma}.\footnotemark
\end{align*}
\footnotetext{$\uh|_{\Sigma}$ is actually $\uf$, hence $\ddpartial_i \big(\uh|_{\Sigma}\big)= \ddpartial_i \uf$. \label{footnote 10}}

Return to the estimate of $(\dslashd \uh)|_{\Sigma}$. Recall that $\uh(s,\vartheta)$ satisfies equation \eqref{eqn 3.2}
\begin{align}
\dpartial_{s} \uh=  -b^i  \dpartial_i \uh + \Omega^2 \left(\slashg^{-1}\right)^{ij} \dpartial_i \uh \dpartial_j \uh.
\tag{\ref{eqn 3.2}}
\end{align}
In the coordinate system $\{s,\vartheta\}$ on $\ucalH$, $\Sigma$ is parameterised as $\Sigma= \{ (s,\vartheta) = (f(\vartheta),\vartheta)\}$, hence $(\dslashd \uh)_i|_{\Sigma}$ is simply the restriction of $\dpartial_i \uh$ on $\Sigma$, i.e.
\begin{align*}
(\dslashd \uh)_i|_{\Sigma}(\vartheta) = \dpartial_i \uh ( f(\vartheta), \vartheta).
\end{align*}

We briefly explain the method to obtain the estimate of $\dpartial_i \uh(f(\vartheta),\vartheta)$. The main tool is lemma \ref{lem 3.1} in subsection \ref{subsec 3.2}:
\begin{enumerate}
\item[a.] differentiate equation \eqref{eqn 3.2} to obtain equations for $\dpartial_i \uh(s,\vartheta),i=1,2$ on $\ucalH$;
\item[b.] introduce the family of surfaces $\{S_t\}$ as in subsection \ref{subsec 3.2} method II, that $S_t=\{(s,\vartheta) = (t f(\vartheta),\vartheta)\}$ in $\ucalH$;
\item[c.] apply lemma \ref{lem 3.1} to obtain the equation for $\dpartial_i \uh(tf(\vartheta),\vartheta)$ along $\{S_t\}$;
\item[d.] integrate the equation to obtain the estimate of $\dpartial_i \uh ( f(\vartheta), \vartheta)$.
\end{enumerate}
When carrying out the above procedures, it encounters the issue that $\dpartial_i \uh(s,\vartheta)$ isnot globally well-defined functions on $\ucalH$ since there exists no coordinate system $\vartheta$ covering the whole sphere. Therefore we employ the rotational vector field derivatives instead of the coordinate derivatives. Referring to appendix \ref{appen R} for the basics of the rotational vector field derivatives.

Let $\dR_i$ be the rotational vector field on $\Sigma_s = \ucalH \cap C_s$ induced by an isometric embedding of $(\Sigma_{0,0}, \circg)$ as in appendix \ref{appen R}, then $[\dpartial_s, \dR_i]=0$. We have the following formula of $\dR_i$
\begin{align*}
\dR_i = R_i + (\dR_i \uh) \cdot \partial_{\us}.
\end{align*}
Differentiate equation \eqref{eqn 3.2} in the direction of the rotational vector field $\dR_i$,
\begin{align}
\dpartial_s \dR_k \uh
=
&
- b^i \dpartial_i ( \dR_k \uh ) 
+ 2 \Omega^2 ( \slashg^{-1} )^{ij} \dpartial_j \uh\ \dpartial_i ( \dR_k \uh )
- [ R_k, b]^i \dpartial_i \uh
- ( \dR_k \uh ) \partial_{\us} b^i \dpartial_i \uh
\nonumber
\\
&
+ \big[ \lie_{R_k} (\Omega^2 \slashg^{-1} ) \big]^{ij} \dpartial_i \uh \dpartial_j \uh
+ (\dR_k \uh ) \big[ \partial_{\us} (\Omega^2 \slashg^{-1} ) \big]^{ij} \dpartial_i \uh \dpartial_j \uh.
\label{eqn 5.1}
\end{align}
We rewrite the above equation in terms of the rotational vector field components. Denote the rotational derivative $\dR_k \uh$ by $\uh_{\dR,k}$, then the above equation is equivalent to
\begin{align}
\dpartial_s \uh_{\dR,k}
=
&
- b^i \dpartial_i ( \uh_{\dR,k} ) 
+ 2 \Omega^2 ( \slashg^{-1} )^{R,\overline{i}j} \uh_{\dR,j} \dpartial_i ( \uh_{\dR,k} )
- [ R_k, b ]^{R,i} \uh_{\dR,i}
- ( \partial_{\us} b )^{R,i} \uh_{\dR,i} \uh_{\dR,k}
\nonumber
\\
&
+ \big[ \lie_{R_k} (\Omega^2 \slashg^{-1} ) \big]^{R, ij} \uh_{\dR,i} \uh_{\dR,j}
+ \big[ \partial_{\us} (\Omega^2 \slashg^{-1} ) \big]^{R,ij} \uh_{\dR,i} \uh_{\dR,j} \uh_{\dR,k}.
\footnotemark
\label{eqn 5.2}
\end{align}
\footnotetext{See the meaning of mixed components $\left( \slashg^{-1} \right)^{R,\overline{i}j}$ in appendix \ref{appen R}.}

Equations \eqref{eqn 5.2} with $k=1,2,3$ form a system of equations for $\uh_{\dR,k}$. Then applying lemma \ref{lem 3.1} to this system and the family of surface $\{ S_t \}$, we obtain the equations for $\uh_{\dR,k}(t f(\vartheta), \vartheta)$. Denote $\uh_{\dR,k}(t f(\vartheta), \vartheta)$ by $\uhl{t}_{\dR,k}(\vartheta)$,
\begin{align*}
\partial_t \uhlt_{\dR,k}
=
&
- f b^i \big[ \ddpartial_i \uhlt_{\dR,k} - t f_i f^{-1} \partial_t \uhlt_{\dR,k}  \big]
\\
&
+ 2 f \Omega^2 ( \slashg^{-1} )^{R,\overline{i}j} \, \uhlt_{\dR,j} \big[ \ddpartial_i \uhlt_{\dR,k} - t f_i f^{-1} \partial_t \uhlt_{\dR,k}  \big]
\\
&
- f [ R_k, b ]^{R,i}  \, \uhl{t}_{\dR,i} 
- f ( \partial_{\us} b )^{R,i} \, \uhlt_{\dR,i} \uhlt_{\dR,k}
\nonumber
\\
&
+ f \big[ \lie_{R_k} (\Omega^2 \slashg^{-1} ) \big]^{R, ij}  \, \uhlt_{\dR,i} \uhlt_{\dR,j}
+ f \big[ \partial_{\us} (\Omega^2 \slashg^{-1} ) \big]^{R,ij} \, \uhlt_{\dR,i} \uhlt_{\dR,j} \uhlt_{\dR,k},
\end{align*}
We rewrite the above equation as follows
\begin{align}
\begin{aligned}
&
\partial_t\, \uhlt_{R,k}
=
\Xl{t}_{\uh}^i \ddpartial_i\, \uhlt_{\dR,k} + \rel{t}_{\uh,k},
\\
&
\uhl{t=0}_{\dR,k} = \ufl{s=0}_{\dR,k}= \dR_k \left( \ufl{s=0} \right),
\end{aligned}
\label{eqn 5.3}
\end{align}
where
\begin{align*}
\Xlt_{\uh} 
=
f \big[ 1- t b^m f_m - 2 t \Omega^2 ( \slashg^{-1} )^{R,\overline{m}j} f_m \,\uhlt_{\dR,j}  \big]^{-1}
\cdot
\big[ - b^i + 2 \Omega^2 ( \slashg^{-1} )^{R,\overline{i}j}\,\uhlt_{\dR,j}  \big]
\ddpartial_i,
\end{align*}
and
\begin{align*}
\relt_{\uh,k}
=&
f \big[ 1- t b^m f_m - 2 t \Omega^2 \left( \slashg^{-1} \right)^{R,\overline{m}j} f_m \,\uhlt_{\dR,j}  \big]^{-1}
\cdot
\\
&
\left\{
- [ R_k, b ]^{R,i}  \, \uhl{t}_{\dR,i} 
- ( \partial_{\us} b )^{R,i} \, \uhlt_{\dR,i} \uhlt_{\dR,k}
\right.
\\
&
\phantom{\{}
\left.
+ \big[ \lie_{R_k} (\Omega^2 \slashg^{-1} ) \big]^{R, ij}  \, \uhlt_{\dR,i} \uhlt_{\dR,j}
+ \big[ \partial_{\us} (\Omega^2 \slashg^{-1} ) \big]^{R,ij} \, \uhlt_{\dR,i} \uhlt_{\dR,j} \uhlt_{\dR,k} 
\right\}.
\end{align*}
We shall integrate equation \eqref{eqn 5.3} to obtain the estimate for the differential $\dslashd \uh$ on $\Sigma$.
\begin{proposition}\label{pro 5.1}
Let $\Sigma$ be a spacelike surface embedding in an incoming null hypersurface $\ucalH$ in $(M,g)$. Assume that $\Sigma$ has the second parameterisation $\big( \ufl{s=0}, f \big)$ and $\ucalH$ is parameterised by $\uh$ as in subsection \ref{subsec 3.2} method I. Suppose that the parameterisation functions $\ufl{s=0}, f$ satisfy the estimates,
\begin{align*}
\leVe \dslashd \ufl{s=0} \riVe^{n+1,p} \leq \udelta_{o} r_0, 
\quad 
\big\vert \overline{\ufl{s=0}}^{\circg} \big\vert \leq \udelta_m r_0,
\quad
\leVe \ddslashd f \riVe^{n+1,p} \leq \delta_{o} r_0,
\quad
\overline{f}^{\circg} = \os,
\end{align*}
where $n\geq 2, p>2$ or $n\geq 3, p>1$.

There exist a small positive constant $\delta$ and a constant $c_{\uh}$ both depending on $n,p$, such that if $\epsilon, \udelta_o, \udelta_m, \delta_o$ are suitably bounded that $\epsilon, \udelta_o, \epsilon\udelta_m, \delta_o \leq \delta$, then the differential $(\dslashd \uh )|_{\Sigma}$ satisfies the following estimate
\begin{align}
\leVe ( \dslashd \uh )|_{\Sigma} \riVe^{n+1,p} \leq c_{\uh} \leVe \dslashd \ufl{s=0} \riVe^{n+1,p} \leq c_{\uh} \udelta_o r_0.
\label{eqn 5.4}
\end{align}
\end{proposition}
\begin{proof}
The proof is similar to the one of proposition \ref{pro 3.3}. We use the bootstrap argument. Introduce the following bootstrap assumption.
\begin{assum}
Estimate \eqref{eqn 5.4} holds for $( \dslashd \uh )|_{S_t}$ in the closed interval $[0,t_a]$.
\end{assum}
By continuity, if $c_{\uh}>1$, then since $( \dslashd \uh )|_{S_0} = \dslashd \ufl{s=0}$, there exists some small interval $[0,t_a]$ such that the assumption holds.

In the following, we shall show that if the bootstrap assumption is true, then for suitably chosen $\delta$ and $c_{\uh}$ independent of $t_a$, estimate \eqref{eqn 5.4} can be improved to the strict inequality at $t=t_a$.

Assume that $\delta\leq 1/2$ and $c_{\uh} \delta \leq 1$. Furthermore, assume that $\delta$ is suitably small such that proposition \ref{pro 3.3} holds. Let $\ud_o = c_o \udelta_o, \ud_m = \big[ 1+c_{m,m} (|\os|/r_0+ \delta_o)  \epsilon \udelta_o \big] \udelta_m + c_{m,o} (|\os|/r_0+ \delta_o)  \udelta_o^2$ as in equations \eqref{eqn 3.7}. Then by proposition \ref{pro 3.3},
\begin{align*}
\leVe \ddslashd \uflt \riVe^{n,p} \leq \ud_o r_0,
\quad
\leve \overline{\uflt}^{\circg} \rive \leq \ud_m r_0.
\end{align*}
By the bootstrap assumption and the estimates of $\uflt, f$, we have\footnote{We abuse the notation $c(n,p)$ as in footnote \ref{footnote 6}.}
\begin{align*}
&
\leVe \Xlt_{\uh} \riVe^{n+1,p}
\leq
c(n,p)(|\os|/r_0 + \delta_o) \big( \epsilon ( \ud_m+ \ud_o)  + c_{\uh} \udelta_o \big)
\leq
c(n,p),
\\
&
\leVe \relt_{\uh,k} \riVe^{n+1,p}
\leq
c(n,p) (|\os|/r_0 + \delta_o)  \big( \epsilon (\ud_m+\ud_o) \cdot c_{\uh} \udelta_o + \epsilon c_{\uh}^2 \udelta_o^2 + c_{\uh}^3 \udelta_o^3 \big) r_0,
\end{align*}
which follow from
\begin{align}
\begin{aligned}
\Vert \vec{b} \Vert_{S_t}^{n+1,p}, \Vert [R_k, \vec{b}] \Vert_{S_t}^{n+1,p} \leq c(n,p) \epsilon ( \ud_m + \ud_o ) r_0^{-1},
\\
\leVe \partial_{\us} b \riVe_{S_t}^{n+1,p} \leq c(n,p) \epsilon r_0^{-2},
\\
\leVe \lie_{R_k} (\Omega^2 \slashg^{-1}) \riVe_{S_t}^{n+1,p} \leq c(n,p) \epsilon r_0^{-2},
\\
\leVe \partial_{\us} (\Omega^2 \slashg^{-1}) \riVe_{S_t}^{n+1,p} \leq c(n,p) r_0^{-3}.
\end{aligned}
\label{eqn 5.5}
\end{align}
Then by Gronwall's inequality, integrate equation \eqref{eqn 5.3} to obtain
\begin{align*}
\leVe \big(\dslashd \uh\big)|_{S_t} \riVe^{n+1,p}
\leq
& 
c(n,p) \leVe \uhlt_{\dR,k} \riVe^{n+1,p}
\\
\leq
&
c(n,p) \left\{
\leVe \ufl{s=0}_{\dR,k} \riVe^{n+1,p}
+ 
\int_0^t \leVe \relt_{\uh,k} \riVe^{n+1,p} \d t
\right\}
\\
\leq
&
c(n,p) \left[ \udelta_o + \epsilon (\ud_m+ \ud_o) \cdot c_{\uh} \udelta_o + \epsilon c_{\uh}^2 \udelta_o^2 + c_{\uh}^3 \udelta_o^3  \right] r_0.
\end{align*}
Therefore in order to strengthen estimate \eqref{eqn 5.4} to strict inequality at the end point $t_a$, it is sufficient to require that
\begin{align*}
c(n,p) \left( 1 +  \epsilon (\ud_m + \ud_o) c_{\uh}  + \epsilon \udelta_o c_{\uh}^2 + \udelta_o^2 c_{\uh}^3 \right)
\leq
c(n,p)  \left( 1 +  \delta c_{\uh}  + \delta^2 c_{\uh}^2 + \delta^2 c_{\uh}^3 \right)
< c_{\uh}.
\end{align*}
Thus we choose $c_{\uh}=2c(n,p)$, and $\delta$ suitably small such that the above inequality holds. For such $\delta$ and $c_{\uh}$, the bootstrap argument is closed and the proposition is proved.
\end{proof}
\begin{remark}
The above proof requires the $L^{\infty}$ bounds of the metric components $b, \Omega, \slashg$ and their derivatives up to $(n+2)$-th order in estimates \eqref{eqn 5.5}. The proof also makes use of the estimates of $\uf$ obtained in proposition \ref{pro 3.3}. However, in return, proposition \ref{pro 5.1} gives a better estimate for $\uf$ with a higher regularity, which is proposition \ref{pro 3.5}. We present the proof of proposition \ref{pro 3.5} in the following.
\end{remark}
\begin{proof}[Proof of proposition \ref{pro 3.5}]
As we mentioned in footnote \ref{footnote 10},
\begin{align*}
\ddslashd \uf = \ddslashd (\uh|_{\Sigma}) = (\dslashd \uh )|_{\Sigma} + \ddslashd f \cdot (\dpartial_s \uh )|_{\Sigma}.
\end{align*}
Then substituting $\dpartial_s \uh$ from equation \eqref{eqn 3.2},
\begin{align*}
\ddslashd \uf 
= 
(\dslashd \uh )|_{\Sigma} 
+ 
\ddslashd f \cdot \big( -b^i  \dpartial_i \uh + \Omega^2 (\slashg^{-1} )^{ij} \dpartial_i \uh \dpartial_j \uh \big)\Big|_{\Sigma}.
\end{align*}
Therefore by estimate \eqref{eqn 5.4} of $( \dslashd \uh )|_{\Sigma}$ and
\begin{align*}
\Vert \ddslashd f \Vert^{n+1,p} \leq \delta_o r_0,
\quad
\Vert \vec{b} \Vert_{\Sigma}^{n+1,p} \leq c(n,p) \epsilon ( \ud_m + \ud_o) r_0^{-1},
\quad
\Vert \Omega^2 \slashg^{-1} \Vert_{\Sigma}^{n+1,p} \leq c(n,p) r_0^{-2},
\end{align*}
we obtain the estimate for $\ddslashd \uf$, that for $\delta$ suitably small
\begin{align*}
\leVe \ddslashd \uf \leVe^{n+1,p}
\leq
c_{\uh} \udelta_o r_0 + c(n,p) \delta_o \big(\epsilon (\ud_m+\ud_o) c_{\uh} \udelta_o  + c_{\uh}^2 \udelta_o^2 \big) r_0
\leq
(c_{\uh} + c(n,p)) \udelta_o r_0.
\end{align*}
Then proposition \ref{pro 3.5} is proved.
\end{proof}

\subsection{Estimate of Hessian of parameterisation function $\uh$}\label{subsec 5.2}
We follow the similar route to estimate the Hessian $\circnabla^2 \uh$ on $\Sigma$ as in the previous subsection. First derive the equation for $( \circnabla^2 \uh )|_{\Sigma}$ from equation \eqref{eqn 3.2},  then integrate the equation to obtain the estimate. Again we employ the rotational vector field derivatives. The Hessian $\circnabla^2 \uh$ satisfies the following formula
\begin{align*}
\circnabla^2 \uh(\dR_l , \dR_k)= \dR_l \dR_k \uh + \epsilon_{lij} x_k x_i \big(\dR_j \uh \big)
\end{align*}
Denote $ \dR_l \dR_k \uh$ by $\uh_{\dR,lk}$, then the above formula can be written as
\begin{align*}
( \circnabla^2 \uh )_{\dR, lk} = \uh_{\dR,lk} +  \epsilon_{lij} x_k x_i \uh_{\dR,j}.
\end{align*}
Note that $\uh_{\dR,lk}$ is not symmetric in the two indices $k,l$. By this formula, it is sufficient to estimate $\uh_{\dR,lk}|_{\Sigma}$ in order to estimate $( \circnabla^2 \uh )|_{\Sigma}$.

Differentiating equation \eqref{eqn 5.1} in the direction of the rotational vector fields, then by the Leibniz rule of the Lie derivatives and the formula $\lie_{X} \d = \d \lie_{X}$, we obtain\footnote{Use the square bracket $[X, \xi]$ to denote the Lie derivatives in the direction of $X$, i.e. $[X, \xi]= \lie_X \xi$}
\begin{align}
\dpartial_s \uh_{\dR,lk} 
=
&
- b^i \dpartial_i ( \uh_{\dR,lk} ) 
+ 2 \Omega^2 ( \slashg^{-1} )^{R,\overline{i}j} \uh_{\dR,j} \dpartial_i ( \uh_{\dR,lk} )
\nonumber
\\
Q_{lk}
&
\left\{
\begin{aligned}
&
- \big[\dR_l, b\big]^{\dR,i} \uh_{\dR,ik}  
\\
&
+ 
\left\{ 
2 \big[ \dR_l, \Omega^2 ( \slashg^{-1} ) \big]^{\dR,ij}   \uh_{\dR,j}  \uh_{\dR,ik} 
+ 
2 \Omega^2 ( \slashg^{-1} )^{R,ij} \uh_{\dR,ik} \uh_{\dR,jl}
\right\}
\\
&
-
\left\{
\big[ \dR_l, [ R_k, b ]\big]^{\dR,i} \uh_{\dR,i} 
+
[ R_k, b ]^{R,i} \uh_{\dR,il}
\right\}
\\
&
- 
\left\{
\big[ \dR_l,  \partial_{\us} b \big]^{\dR,i} \uh_{\dR,i} \uh_{\dR,k}
+ ( \partial_{\us} b )^{R,i} \uh_{\dR,il} \uh_{\dR,k}
+ ( \partial_{\us} b )^{R,i} \uh_{\dR,i} \uh_{\dR,lk}
\right\}
\\
&
+
\left\{
\big[\dR_l, [ R_k, \Omega^2 \slashg^{-1} ] \big]^{\dR, ij} \uh_{\dR,i} \uh_{\dR,j}
+ [ R_k, \Omega^2 \slashg^{-1} ]^{R, ij} \uh_{\dR,il} \uh_{\dR,j}
\right.
\\
&\phantom{+ \Big\{[}
\left.
+ [ R_k, \Omega^2 \slashg^{-1} ]^{R, ij} \uh_{\dR,i} \uh_{\dR,jl}
\right\}
\\
&
+ 
\left\{
\big[\dR_l, \partial_{\us} (\Omega^2 \slashg^{-1} ) \big]^{\dR,ij} \uh_{\dR,i} \uh_{\dR,j} \uh_{\dR,k}
+ [ \partial_{\us} (\Omega^2 \slashg^{-1} ) ]^{R,ij} \uh_{\dR,il} \uh_{\dR,j} \uh_{\dR,k}
\right.
\\
&\phantom{+ \Big\{[}
\left.
+ [ \partial_{\us} (\Omega^2 \slashg^{-1} ) ]^{R,ij} \uh_{\dR,i} \uh_{\dR,jl} \uh_{\dR,k}
+ [ \partial_{\us} (\Omega^2 \slashg^{-1} ) ]^{R,ij} \uh_{\dR,i} \uh_{\dR,j} \uh_{\dR,lk}
\right\}.
\end{aligned}
\right.
\label{eqn 5.6}
\end{align}
In $Q_{lk}$, $\dR_l = R_l + \uh_{\dR,l} \partial_{\us}$, and numerically
\begin{align*}
[\dR_l, \xi]^{\dR,i_1\cdots i_k} = [ R_l, \xi]^{R, i_1\cdots i_k} + \uh_{\dR,l} \big( \partial_{\us} \xi \big)^{R,i_1\cdots i_k},
\end{align*}
where $\xi$ could be $b$, $\Omega^2 \slashg^{-1}$, $[R_k, b]$, $\partial_{\us} b$, $[ R_k, \Omega^2 \slashg^{-1} ]$, $\partial_{\us} (\Omega^2 \slashg^{-1} )$.

Applying lemma \ref{lem 3.1} to equation \eqref{eqn 5.6} and the family of surfaces $\{S_t\}$, and denoting $\uh_{\dR,lk}|_{S_t}$ by $\uhlt_{\dR,lk}$, we obtain the equation for $\uhlt_{\dR,lk}$
\begin{align}
\begin{aligned}
&
\partial_t \uhlt_{\dR,lk}
=
\Xlt_{\uh}^i \ddpartial_i \uhlt_{\dR,lk}
+
\relt_{\uh,lk},
\\
&
\uhl{t=0}_{\dR,lk} = \ufl{s=0}_{\dR,lk},
\end{aligned}
\label{eqn 5.7}
\end{align}
where
\begin{align*}
\Xlt_{\uh} 
=
f \big[ 1- t b^m f_m - 2 t \Omega^2 \left( \slashg^{-1} \right)^{R,\overline{m}j} f_m \,\uhlt_{\dR,j}  \big]^{-1}
\cdot
\big[ - b^i + 2 \Omega^2 \left( \slashg^{-1} \right)^{R,\overline{i}j}\,\uhlt_{\dR,j}  \big]
\ddpartial_i,
\end{align*}
and
\begin{align*}
\relt_{\uh,lk}
=&
f \big[ 1- t b^m f_m - 2 t \Omega^2 \left( \slashg^{-1} \right)^{R,\overline{m}j} f_m \,\uhlt_{\dR,j} \big]^{-1}
\cdot
Q_{lk}|_{S_t}
\end{align*}
Then we can integrate equation \eqref{eqn 5.7} to get the estimate for $( \circnabla^2 \uh )|_{\Sigma}$.
\begin{proposition}\label{pro 5.3}
Under the same setting of proposition \ref{pro 5.1}, there exist a small positive constant $\delta$ and a constant $c_{\uh,2}$ both depending on $n,p$, such that if $\epsilon, \udelta_o, \epsilon \udelta_m, \delta_o \leq \delta$, then the Hessian $( \circnabla^2 \uh )|_{\Sigma}$ satisfies the following estimate
\begin{align*}
\leVe ( \circnabla^2 \uh )|_{\Sigma} \riVe^{n,p} \leq c_{\uh,2} \leVe \dslashd \ufl{s=0} \riVe^{n+1,p} \leq c_{\uh,2} \udelta_o r_0.
\end{align*}
\end{proposition}
\begin{proof}
The proof is essentially the same as the proof proposition \ref{pro 5.1}, thus we just list some important points here. Choose $c_{\uh,2} \geq c_{\uh}$, then $\leVe (\dslashd \uh )|_{\Sigma} \riVe^{n+1,p}$ and $\leVe (\circnabla^2 \uh )|_{\Sigma} \riVe^{n,p}$ can be bounded by the same quantity $c_{\uh,2} \udelta_o r_0$. This is just a technical assumption which simplifies the expressions of some estimates in the proof.

The key is to show that $\relt_{\uh,lk}$ satisfies the estimate
\begin{align*}
\leVe \relt_{\uh,lk} \riVe^{n,p}
\leq
c(n,p)(|\os|/r_0 + \delta_o) 
\big( 
\epsilon (\ud_m+\ud_o) \cdot c_{\uh,2} \udelta_o + c_{\uh,2}^2 \udelta_o^2 
\big) r_0.
\end{align*}
This is proved in the same manner as in the proof of proposition \ref{pro 5.1}, thus the detail is omitted. We remark that this proposition also requires the $L^{\infty}$ bounds of the metric components up to their $(n+2)$-th order derivatives, the same as proposition \ref{pro 5.1}.
\end{proof}

\subsection{Estimate of the outgoing null expansion}
We already obtained the estimates for $\uf$, $\big(\dslashd \uh\big)|_{\Sigma}$ and $\big( \circnabla^2 \uh \big)|_{\Sigma}$. We are ready to estimate the outgoing null expansion now.

Estimate the first order main part $\lo{\ddtr \ddchi'}$ and the high order remainder part $\hi{\ddtr \ddchi'}$ separately. Following the notations in the subsection \ref{subsec 4.3},
\begin{align*}
&
\lo{\ddtr \ddchi'}
=
\tr \chi'_S|_{\Sigma} - 2 \left( r_S|_{\Sigma}\right)^{-2} \circDelta f,
\\
&
\hi{\ddtr \ddchi'} = \ddtr \ddchi' - \lo{\ddtr \ddchi'}
\end{align*}

From the values of $\tr \chi_S|_{\Sigma}$ in formulae \eqref{eqn 2.2}, and $\Omega_S^2, r_S$ in formulae \eqref{eqn 2.1}, we have
\begin{align*}
\tr \chi'_S|_{\Sigma} = \frac{s}{s + r_0} \cdot \frac{2}{r_S},
\quad
r_S|_{\Sigma} = r_0 + s + \frac{s \, \us}{r_0+s} + O\left(\frac{s \, \us^2}{r_0^2}\right),
\end{align*}
Then we get the value of $\lo{\ddtr \ddchi'}$ from the above. 

In the following, we shall estimate $\hi{\ddtr \ddchi'}$. Given the bounds of parametrisation functions $(\ufl{s=0}, f)$ in proposition \ref{pro 3.3}, i.e.
\begin{align*}
\leVe \dslashd \ufl{s=0} \riVe^{n+1,p} \leq \udelta_{o} r_0, 
\quad
\overline{\ufl{s=0}}^{\circg} = \ous,
\quad 
\Big\vert \overline{\ufl{s=0}}^{\circg} \Big\vert \leq \udelta_m r_0,
\quad
\leVe \ddslashd f \riVe^{n+1,p} \leq \delta_{o} r_0,
\quad
\overline{f}^{\circg} = \os,
\end{align*}
where $n\geq 2, p>2$ or $n\geq 3, p>1$, let
\begin{align}
\ud_o = c_o \udelta_o, 
\quad
\uslashd_m = (|\os|/r_0+ \delta_o) (c_{m,m} \epsilon \udelta_m \udelta_o + c_{m,o} \udelta_o^2),
\label{eqn 5.8}
\footnotemark
\end{align}
\footnotetext{The simpler notation $\ud_m$ is already used to denote $\big[ 1+c_{m,m} (|\os|/r_0+ \delta_o)  \epsilon \udelta_o \big] \udelta_m + c_{m,o} (|\os|/r_0+ \delta_o)  \udelta_o^2$ in formulae \eqref{eqn 3.7} in the proof of proposition \ref{pro 3.3}.}
then the parameterisation function $\uf$ satisfies the estimates
\begin{align*}
\leVe \ddslashd \uf \riVe^{n,p} \leq \ud_o r_0,
\quad
\leve \overline{\uf}^{\circg} - \ous \rive \leq \uslashd_m r_0.
\end{align*}
Introduce the notation $\ud_{\uh}= \max\{ c_{\uh}, c_{\uh,2}  \} \udelta_o$, then y propositions \ref{pro 5.1}, \ref{pro 5.3},
\begin{align*}
\leVe ( \dslashd \uh )|_{\Sigma} \riVe^{n+1,p}
\leq
\ud_{\uh} r_0,
\quad
\leVe ( \circnabla^2 \uh )|_{\Sigma} \riVe^{n,p}
\leq
\ud_{\uh}  r_0,
\end{align*}

With the above estimates of $f, \uf$, $\big( \dslashd \uh \big)|_{\Sigma}$, $\big( \circnabla^2 \uh \big)|_{\Sigma}$, we shall prove the following proposition on the estimate of the high order remainder part $\hi{\ddtr \ddchi'}$.
\begin{proposition}\label{pro 5.4}
Let $\Sigma$ be a spacelike surface in $(M,g)$. Assume that it has the second parameterisation $(\ufl{s=0}, f)$, and the parameterisation functions satisfy the estimates
\begin{align*}
\leVe \dslashd \ufl{s=0} \riVe^{n+1,p} \leq \udelta_{o} r_0, 
\quad
\overline{\ufl{s=0}}^{\circg} = \ous,
\quad 
\leve \overline{\ufl{s=0}}^{\circg} \rive \leq \udelta_m r_0,
\quad
\Vert \ddslashd f \Vert^{n+1,p} \leq \delta_{o} r_0,
\quad
\overline{f}^{\circg} = \os,
\end{align*}
where $n\geq 2, p>2$ or $n\geq 3, p>1$.

There exists a small positive constant $\delta$ depending on $n,p$, such that if $\epsilon, \udelta_o, \epsilon \udelta_m, \delta_o \leq \delta$, then $\hi{\ddtr \ddchi'}$ satisfies the estimate
\begin{align*}
&
\leVe \hir{1}{\ddtr \ddchi'} \riVe^{n,p}
\leq
c(n,p) \epsilon r_0^{-1},
\\
&
\leVe \hir{2}{\ddtr \ddchi'} \riVe^{n,p}
\leq
c(n,p) ( \epsilon + \delta_o +  \ud_{\uh} |\os|/r_0) \delta_o r_0^{-1}.
\end{align*}
\end{proposition}
\begin{proof}
By the $L^{\infty}$ bounds in definition \ref{def 2.3}, the estimates of $f$ and $\uf$ in proposition \ref{pro 3.3}, the estimate of $( \dslashd \uh )|_{\Sigma}$ in proposition \ref{pro 5.1}, we have
\begin{align*}
&
\Vert \lo{\ddslashg} \Vert^{n,p}
\leq
c(n,p) r_0^2,
\quad 
\Vert \hir{1}{\ddslashg} \Vert^{n,p} 
\leq
c(n,p) \epsilon r_0^2,
\\
&
\Vert \db^i  \Vert_{\Sigma}^{n,p} 
\leq
c(n,p) [\epsilon (\ud_m + \ud_o) + \ud_{\uh} ] r_0^{-1},
\\
&
\Vert \hir{2}{\ddslashg} \Vert^{n,p} 
\leq
c(n,p) [ \epsilon (\ud_m + \ud_o) + \ud_{\uh} ] \delta_o r_0^2,
\end{align*}
and
\begin{align*}
&
\Vert \lo{\ddchi'} \Vert^{n,p}
\leq
c(n,p)(|\os|+\delta_o r_0),
\quad
\Vert \hir{1}{\ddchi'} \Vert^{n,p} 
\leq
c(n,p) \epsilon r_0.
\end{align*}
Note $\big( \dslashd \uh \big)|_{\Sigma}$ appears in $\db^i= b^i - 2 \Omega^2 \left( \slashg^{-1} \right)^{ij} \uh_j$ and $\hir{2}{\ddslashg}$.

There remains the estimate of $\hir{2}{\ddchi'}$. We list the following estimates for the terms in the formula of $\hir{2}{\ddchi'}$.
\begin{align*}
&
\begin{aligned}
&
\Vert \circtriangle \Vert_{\Sigma}^{n,p} \leq c(n,p) \epsilon,
&&
\Vert \triangle \Vert_{\Sigma}^{n,p} \leq c(n,p) (\epsilon + |\os|/r_0 + \delta_o )\ud_{\uh},
\\
&
\Vert \vec{\dvarepsilon}' \Vert^{n,p} \leq c(n,p) \delta_o r_0^{-1},
&&
\Vert \dvarepsilon' \Vert^{n,p} \leq c(n,p) \delta_o^2,
\\
&
\Vert \duchi \Vert_{\Sigma}^{n,p} \leq c(n,p) r_0,
&&
\Vert \dchi' \Vert_{\Sigma}^{n,p} \leq c(n,p) (\epsilon+ |\os|/r_0 + \delta_o) r_0,
\\
&
\Vert \deta \Vert_{\Sigma}^{n,p} \leq c(n,p) [\epsilon + (|\os|/r_0 + \delta_o) \ud_h],
\end{aligned}
\\
&
\Vert \duomega \Vert_{\Sigma}^{n,p} \leq c(n,p) [ \epsilon + \ud_m + \ud_o + \epsilon \ud_{\uh} + (|\os|/r_0 + \delta_o) \ud_h^2 ] r_0^{-1},
\\
&
\begin{aligned}
&
\Vert \dslashnabla \db \Vert_{\Sigma}^{n,p} \leq c(n,p) [ \epsilon (\ud_m + \ud_o ) +  \ud_{\uh} ] r_0^{-1},
&&
\Vert \dpartial_s \db \Vert_{\Sigma}^{n,p} \leq c(n,p) [ \epsilon (\ud_m + \ud_o ) +  \ud_{\uh} ]  r_0^{-2}.
\end{aligned}
\end{align*}
Note that $\big( \dslashd \uh \big)_{\Sigma}$ appears in $\triangle$, $\vec{\dvarepsilon}'$, $\dvarepsilon'$, $\duchi$, $\deta$, $\duomega$, and $\big( \circnabla^2  \uh \big)_{\Sigma}$ appears in $\duchi$, $\dslashnabla \db$, $\dpartial_{\us} \db$. We shall explain more on the estimates of $\dslashnabla \db$, $\dpartial_{\us} \db$. Their estimates follow from
\begin{align*}
\dslashnabla_k \db^i
=&
\circnabla_k \big[ (b|_{\Sigma})^i - 2(\Omega^2 \slashg^{-1})|_{\Sigma}^{ij} \, \uh_j \big]
+
\big( \triangle_{kl}^i + \circtriangle_{kl}^i \big) \big( b^l - 2\Omega^2 ( \slashg^{-1} )^{lj} \uh_j \big),
\\
\dpartial_s \db^i
=&
( \partial_s b^i + \dpartial_s \uh \cdot \partial_{\us} b^i  ) 
- 2 [ \partial_s ( \Omega^2 \slashg^{-1} ) ]^{R,\overline{i}j} \uh_{\dR,j} 
\\
&
-2 \dpartial_s \uh \cdot [ \partial_{\us} ( \Omega^2 \slashg^{-1} ) ]^{R,\overline{i}j} \uh_{\dR,j} 
- 2 \Omega^2 ( \slashg^{-1} )^{R,\overline{i}j} \dpartial_s \uh_{\dR,j},
\end{align*}
and equations \eqref{eqn 3.2}, \eqref{eqn 5.2} of $\dpartial_s \uh$, $\dpartial_s \uh_{\dR,k}$.

Assembling the above listed estimates into $\hir{2}{\ddchi'}$, we obtain that
\begin{align*}
\leVe \hir{2}{\ddchi'} \riVe^{n,p}
\leq
c(n,p) (\epsilon + \delta_o + \ud_{\uh} |\os|/r_0 ) \delta_o r_0.
\end{align*}
Then substituting the estimates of $\hir{a}{\ddslashg}$, $\hir{a}{\ddchi'}$ into the formulae of $\hir{1}{\ddchi'}$ and $\hir{2}{\ddchi'}$, the proposition is proved.
\end{proof}
\begin{remark}
Note that proposition \ref{pro 5.4} requires the $L^{\infty}$ bounds of the structure coefficients up to $n$th order derivatives. It also requires the $L^{\infty}$ bounds of the metric components up to $(n+2)$-th order derivatives, which follows from the estimate of $( \circnabla^2 \uh )|_{\Sigma}$.
\end{remark}

\section{Perturbation of parameterisation of spacelike surface}\label{sec 6}
In this section, we study the following problem. Let $\Sigmal{a}, a=1,2$ be two spacelike surfaces in $(M,g)$. Suppose that  the second parameterisation of $\Sigmal{a}$ is $\left( \ufl{a,s=0}, \fl{a} \right)$ and the first parametrisation is $\left( \ufl{a}, \fl{a} \right)$. Define the perturbation functions
\begin{align*}
\dd{\ufl{s=0}} = \ufl{2,s=0} - \ufl{1,s=0},
\quad
\dd{\uf} = \ufl{2} - \ufl{1},
\quad
\dd{f} = \fl{2} - \fl{1}.
\end{align*}
We will show how to obtain $\dd{\uf}$ from $\dd{\ufl{s=0}}$ and $\dd{f}$. We shall also give an estimate for $\dd{\uf}$ in terms of the bounds of $\dd{\ufl{s=0}}$ and $\dd{f}$.

The results in this section will be applied to the perturbation of the outgoing null expansion on $\Sigmal{a}$ in section \ref{sec 7}.

\subsection{Equation of perturbation function}
We apply method I in section \ref{sec 3} to each $\Sigmal{a}$ to obtain $\ufl{a}$. Let $\{\Sl{a}_t\}$ be the family of surfaces with second parametrisation $(\ufl{s=0} , t \fl{a})$. Suppose that the first parameterisation of $\Sl{a}_t$ is $(\ufl{a,t}, t \fl{a} )$. See figure \ref{fig 9}.
\begin{figure}[h]
\begin{center}
\begin{tikzpicture}
\draw[dashed] (-0.9,-0.1) to [out=135,in=180] (0,0.4) node[above]{\scriptsize $\Sl{1}_{t=0}$}  to[out=0,in=45] (0.9,-0.1);
\draw (1,0) to [out=-135,in=-45] (-1,0);
\draw[dashed] (-1,0) to [out=70,in=-110] (-0.7,0.9);
\draw (-1,0) to [out=-110,in=70] (-2.4,-4);
\draw[dashed] (1,0) to [out=110,in=-70] (0.7,0.9);
\draw[->] (1,0) to [out=-70,in=110] (2.4,-4) node[right] {\small $s$}; 
\draw[dashed] (-1.5,-1.5) to [out=70,in=180] (-0.4,-1.7) node[below] {\scriptsize $\Sl{1}_t$} to [out=0,in=110] (1.5,-1.5);
\draw (1.5,-1.5) to [out=-70,in=0] (0.7,-1.5) to [out=180,in=-110] (-1.5,-1.5); 
\draw[dashed] (-2.2,-3.5) to [out=70,in=180] (-0.5,-3.8) node[below] {\scriptsize $\Sl{1}_{t=1}: \Sigma$} to [out=0,in=110] (2.2,-3.5);
\draw (2.2,-3.5) to [out=-70,in=0] (1,-3.3) to [out=180,in=-110] (-2.2,-3.5); 
\draw[dashed] (-1,0) to [out=-45,in=135] (-0.5,-0.5);
\draw (-1,0) to [out=135, in= -45] (-2.3,1.3);
\draw[dashed] (1,0) to [out=-135,in=45] (0.5,-0.5);
\draw[->] (1,0) to [out=45, in= -135] (1.8,0.8) 
to [out=45,in=-135] (2.3,1.3) node[right] {\small $\us$};
\draw[dashed] (-1.5,0.5) to [out=135,in=180] (0,1.1) node[above]{\scriptsize $\Sl{2}_{t=0}$}  to[out=0,in=45] (1.5,0.5);
\draw (1.5,0.5) to [out=-135,in=0] (0,0) to [out=180,in=-45] (-1.5,0.5); 
\draw[dashed] (-1.4,1.3) to [out=-110,in=70] (-1.6,0.7);
\draw (-1.6,0.7) to [out=-110,in=70] (-3,-3);
\draw[dashed] (1.4,1.3) to [out=-70,in=110] (1.6,0.7);
\draw[->] (1.6,0.7) to [out=-70,in=110] (3,-3) node[right]{\small $s$}; 
\draw[dashed] (-2.3+0.23,-2.2+1.6) to [out=70,in=180] (-0.4,-0.7) node[below] {\scriptsize $\Sl{2}_t$} to [out=0,in=110] (2.3-0.23,-2.2+1.6);
\draw (2.3-0.23,-2.2+1.6) to [out=-70,in=0] (0.7,-0.5) to [out=180,in=-110] (-2.3+0.23,-2.2+1.6); 
\draw[dashed] (-2.3-0.4,-2.2) to [out=70,in=180] (-0.5,-2.5) node[below] {\scriptsize $\Sl{2}_{t=1}: \Sigma$} to [out=0,in=110] (2.3+0.4,-2.2);
\draw (2.3+0.4,-2.2) to [out=-70,in=0] (1,-2) to [out=180,in=-110] (-2.3-0.4,-2.2); 
\draw[->,thick] (-1,0) -- (-1.65,0.65);
\draw[->,thick] (-1.5,-1.5) -- (-2.3+0.23,-2.2+1.6);
\draw[->,thick] (-2.2,-3.5) -- (-2.3-0.4,-2.2);
\end{tikzpicture}
\end{center}
\caption{Perturbation of surfaces $\Sigmal{a}$ and the family of surfaces $\{ \Sl{a}_t \}$.}
\label{fig 9}
\end{figure}

We have equations \eqref{eqn 3.4} \eqref{eqn 3.5} for $\ufl{a,t}$
\begin{align*}
&
\partial_t \ufl{a,t} 
= 
F\left(  \fl{a},\ t b^i \cdot \fl{a}_i,\ t e^i \cdot \fl{a}_i,\ t \ue^i \cdot \fl{a}_i,\ \uvarepsilon,\  b^i \, \ufl{a,t}_i,\  e^i \, \ufl{a,t}_i,\  \ue^i \, \ufl{a,t}_i  \right),
\tag{\ref{eqn 3.4}}
\\
&
\partial_t \big( \ddcircDelta \ufl{a,t} \big)
=
\Xl{a,t}^i \ddpartial_i \big( \ddcircDelta \ufl{a,t} \big) + \rel{a,t},
\tag{\ref{eqn 3.5}}
\end{align*}
where $\Xl{a,t}, \rel{a,t}$ are the corresponding vector and function on $\Sigmal{a}_t$. We can derive the equations for $\dd{\uflt}$. Introduce that
\begin{align*}
&
\Fl{a,t}= F\left( \fl{a},\ t b^i \cdot \fl{a}_i,\ t e^i \cdot \fl{a}_i,\ t \ue^i \cdot \fl{a}_i,\ \uvarepsilon,\  b^i \, \ufl{a,t}_i,\  e^i \, \ufl{a,t}_i,\  \ue^i \, \ufl{a,t}_i  \right),
\\
&
\dd{\Flt} = \Fl{2,t} - \Fl{1,t},
\quad
\dd{\Xlt} = \Xl{2,t} - \Xl{1,t},
\quad
\dd{\relt} = \rel{2,t} - \rel{1,t}.
\end{align*}
Then $\dd{\uflt}$ satisfies equations
\begin{align}
&
\partial_t \dd{\uflt} 
= 
\dd{\Flt},
\label{eqn 6.1}
\\
&
\partial_t \big( \ddcircDelta \dd{\uflt} \big)
=
\Xl{1,t}^i \ddpartial_i \big( \ddcircDelta \dd{\uflt} \big)+\dd{\Xlt}^i \ddpartial_i \big( \ddcircDelta \ufl{2,t} \big) + \dd{\relt}.
\label{eqn 6.2}
\end{align}
The initial condition of above equations is
\begin{align*}
\dd{\ufl{t=0}} = \dd{\ufl{s=0}}.
\end{align*}
Note that in equation \eqref{eqn 6.2}, the term $\dd{\Xlt}^i \ddpartial_i \Big( \ddcircDelta \ufl{2,t} \Big)$ on right hand side involves 3rd order derivative of $\ufl{2,t}$. Equation \eqref{eqn 6.2} is a propagation equation for 2nd derivative of $\dd{\uflt}$, therefore when integrating equation \eqref{eqn 6.2} to obtain estimates for $\dd{\uflt}$, the regularity of $\dd{\uflt}$ will be one order less than $\ufl{2,t}$. See proposition \ref{pro 6.1} next subsection.

\subsection{Estimate of perturbation function}\label{subsec 6.2}
We shall use equations \eqref{eqn 6.1} \eqref{eqn 6.2} to obtain estimates for $\dd{\uflt}$. A special case where $\fl{1}= \fl{2} \equiv s$ being a constant function is already treated in theorem 4.2 \cite{L4}. Using the similar method as in \cite{L4}, the following proposition in the general case is obtained.
\begin{proposition}\label{pro 6.1}
Let $\Sigmal{a}, a=1,2$ be two spacelike surfaces in $(M,g)$. Suppose that $\Sigmal{a}$ has the second parameterisation $(\ufl{a,s=0}, \fl{a})$. Assume that the parameterisation functions satisfy the estimates
\begin{align*}
&
\leVe \dslashd \ufl{a,s=0} \riVe^{n,p} \leq \udelta_{o} r_0, 
\quad 
\Big\vert \overline{\ufl{a,s=0}}^{\circg} \Big\vert \leq \udelta_m r_0,
\\
&
\leVe \ddslashd \fl{a} \riVe^{n+1,p} \leq \delta_{o} r_0,
\quad
\overline{\fl{a}}^{\circg} = \os_a, 
\quad
\os= \max\{ |\os_1|, |\os_2| \},
\end{align*}
and the perturbation functions satisfy
\begin{align*}
&
\leVe \dslashd \dd{\ufl{s=0}} \riVe^{n-1,p} \leq \ufrakd_o r_0,
\quad
\leve \overline{\dd{\ufl{s=0}}}^{\circg} \rive \leq \ufrakd_m r_0,
\\
&
\leVe \ddslashd \dd{f} \riVe^{n,p} \leq \frakd_o r_0,
\quad
\leve \overline{\dd{f}}^{\circg} \rive \leq \frakd_m r_0,
\end{align*}
where $n\geq 2, p>2$ or $n\geq 3, p>1$.

Assume that the first parameterisation of $\Sigmal{a}$ is $(\ufl{a}, \fl{a})$. There exist a small positive constant $\delta$ depending on $n,p$, and constants $c_o^{\ufrakd_o}, c_o^{\ufrakd_m}, c_m^{\ufrakd_o}, c_m^{\ufrakd_m}, c^{\frakd}$ depending on $n,p$, such that if $\epsilon, \udelta_o, \udelta_m, \delta_o, \ufrakd_o, \ufrakd_m, \frakd_o, \frakd_m$ are suitably bounded that
\begin{align*}
\epsilon,\ \udelta_o,\ \epsilon \udelta_m,\ \delta_o,\ \ufrakd_o,\ \epsilon \ufrakd_m,\ \frakd_o,\ \frakd_m \leq \delta,
\end{align*}
then the perturbation function $\dd{\uf}$ satisfies the following estimates
\begin{align}
\label{eqn 6.3}
\begin{aligned}
&
\begin{aligned}
\leVe \ddslashd \dd{\uf} \riVe^{n-1,p}
\leq&
c^{\ufrakd_o}_o \ufrakd_o r_0 
+ c^{\ufrakd_m}_{o}  \big(|\os|/r_0+\delta_o\big)  (\udelta_o^2 + \epsilon \udelta_o) \ufrakd_m r_0 
\\
&
+ c^{\frakd}( \udelta_o^2 + \epsilon \udelta_m \udelta_o) (\frakd_m + \frakd_o) r_0 ,
\end{aligned}
\\
&
\begin{aligned}
\leve \overline{\dd{\uf}}^{\circg} - \overline{\dd{\ufl{s=0}}}^{\circg} \rive
\leq
&
c^{\ufrakd_m}_m  \big(|\os|/r_0+\delta_o\big)  (\udelta_o^2 + \epsilon \udelta_o) \ufrakd_m r_0 
\\
&
+ c^{\ufrakd_o}_{m}  \big(|\os|/r_0+\delta_o\big) ( \udelta_o + \epsilon \udelta_m) \ufrakd_{o} r_0 
\\
&
+c^{\frakd}( \udelta_o^2 + \epsilon \udelta_m \udelta_o) (\frakd_m + \frakd_o) r_0.
\end{aligned}
\end{aligned}
\end{align}
\end{proposition}
We shall use bootstrap arguments and Gronwall's inequality to integrate equations \eqref{eqn 6.1} \eqref{eqn 6.2}. The keys in the proof are the estimates of $\dd{\Flt}$, $\dd{\Xlt}$ and $\dd{\relt}$. We give the proof sketch here, and the rest of details will be presented in appendix \ref{appen pro 6.1}.
\begin{proof}[Proof sketch]
We prove that estimates \eqref{eqn 6.3} hold for all $\dd{\uflt}, t\in[0,1]$. First assume that $\delta$ and $c_o^{\ufrakd_o}, c_o^{\ufrakd_m}, c_m^{\ufrakd_o}, c_m^{\ufrakd_m}, c^{\frakd}$ satisfy
\begin{align*}
\delta \leq \frac{1}{2},
\quad
(c_o^{\ufrakd_o} + c_o^{\ufrakd_m} + c_m^{\ufrakd_o} + c_m^{\ufrakd_m} + c^{\frakd}) \delta \leq 1,
\end{align*}
and $\delta$ is sufficiently small such that proposition \ref{pro 3.3} hold.

By continuity, there exist some choices of $\delta$, $c_o^{\ufrakd_o}>1, c_o^{\ufrakd_m}, c_m^{\ufrakd_o}, c_m^{\ufrakd_m}, c^{\frakd}$ and some small neighbourhood interval of $t=0$ where estimates \eqref{eqn 6.3} hold. We introduce the following bootstrap assumption.
\begin{assum}
Estimates \eqref{eqn 6.3} hold for $\dd{\uflt}$ in the closed interval $[0, t_a]$.
\end{assum}
We prove that there exist constants $\delta$ and $c_o^{\ufrakd_o}, c_o^{\ufrakd_m}, c_m^{\ufrakd_o}, c_m^{\ufrakd_m}, c^{\frakd}$ independent of $t_a$, such that estimates \ref{eqn 6.3} can be improved to strict inequalities at the end point $t=t_a$ by the bootstrap assumption. We use the notations $\ud_o, \ud_m$ introduced in formula \eqref{eqn 3.7} in the proof of proposition \ref{pro 3.3}, and introduce the notations $\ubfd_o, \uslashbfd_m, \ubfd_m$ to simplify formulae and estimates in the proof, 
\begin{align}
\begin{aligned}
&
\ubfd_o 
=
c^{\ufrakd_o}_o \ufrakd_o  
+ c^{\ufrakd_m}_{o} \big(|\os|/r_0+\delta_o\big) (\udelta_o^2 + \epsilon \udelta_o) \ufrakd_m 
+ c^{\frakd}( \udelta_o^2 + \epsilon \udelta_m \udelta_o) (\frakd_m + \frakd_o),
\\
&
\begin{aligned}
\uslashbfd_m 
= &
c^{\ufrakd_m}_m  \big(|\os|/r_0+\delta_o\big) (\udelta_o^2 + \epsilon \udelta_o) \ufrakd_m  
+ c^{\ufrakd_o}_{m} \big(|\os|/r_0+\delta_o\big) ( \udelta_o + \epsilon \udelta_m) \ufrakd_{o} 
\\
&
+c^{\frakd}( \udelta_o^2 + \epsilon \udelta_m \udelta_o) (\frakd_m + \frakd_o),
\end{aligned}
\\
&
\ubfd_m = \ufrakd_m + \uslashbfd_m.
\end{aligned}
\label{eqn 6.4}
\end{align}

By the bootstrap assumption on the interval $[0,t_a]$, we can show that $\dd{\Flt}, \dd{\Xlt}, \dd{\relt}$ satisfy the following estimates
\begin{align}
\begin{aligned}
&
\begin{aligned}
\leve \dd{\Flt} \rive
\leq
&
c(n,p) \big(|\os|/r_0+\delta_o\big) 
\left[(\ud_o^2 + \epsilon \ud_o) \ubfd_m
+
(\ud_o + \epsilon \ud_m) \ubfd_o \right]
\\
&
+
c(n,p) (\ud_o^2 + \epsilon \ud_o \ud_m) (\frakd_m + \frakd_o) r_0,
\end{aligned}
\\
&
\begin{aligned}
\leVe \dd{\Xlt} \riVe^{n-1,p}
\leq
&
c(n,p)  \big(|\os|/r_0+\delta_o\big) 
\left[ (\ud_o + \epsilon) \ubfd_m
+
\ubfd_o \right]
\\
&
+
c(n,p) (\ud_o + \epsilon \ud_m) (\frakd_m + \frakd_o),
\end{aligned}
\\
&
\begin{aligned}
\leVe \dd{\relt} \riVe^{n-2,p}
\leq
&
c(n,p) \big(|\os|/r_0+\delta_o\big) 
\left[ (\ud_o^2 + \epsilon \ud_o) \ubfd_m
+
(\ud_o + \epsilon \ud_m) \ubfd_o \right]
\\
&
+
c(n,p) (\ud_o^2 + \epsilon \ud_o \ud_m) (\frakd_m + \frakd_o) r_0.
\end{aligned}
\end{aligned}
\label{eqn 6.5}
\end{align}
Therefore integrating equation \eqref{eqn 6.1}, we obtain that
\begin{align*}
\leve \overline{\dd{\uflt}}^{\circg} - \overline{\dd{\ufl{s=0}}}^{\circg} \rive
\leq
&
\int_0^t \leve \overline{\dd{\Flt}}^{\circg} \rive \d t
\\
\leq
&
c(n,p) \big(|\os|/r_0+\delta_o\big) 
\left[
(\ud_o^2 + \epsilon \ud_o) \ubfd_m
+
(\ud_o + \epsilon \ud_m) \ubfd_o 
\right] r_0
\\
&
+
c(n,p) (\ud_o^2 + \epsilon \ud_o \ud_m) (\frakd_m + \frakd_o) r_0.
\end{align*}
Applying Gronwall's inequality to equation \eqref{eqn 6.2}, we obtain that
\begin{align*}
&
\leVe \ddcircDelta \dd{\uflt} \riVe^{n-2,p}
\\
&
\leq
c(n,p) \left\{ \leVe \dcircDelta \dd{\ufl{s=0}} \riVe^{n-2,p} +  \int_0^t \leVe \dd{\Xlt}^i \ddpartial_i \big( \ddcircDelta \ufl{2,t} \big) + \dd{\relt} \riVe^{n-2,p} \d t \right\}.
\end{align*}
The above integrated term has the upper bound
\begin{align}
\begin{aligned}
&
c(n,p)
 \big(|\os|/r_0+\delta_o\big)  
 \left[ (\ud_o^2 + \epsilon \ud_o) \ubfd_m
+
(\ud_o + \epsilon \ud_m) \ubfd_o \right] r_0
\\
&
+
c(n,p) (\ud_o^2 + \epsilon \ud_o \ud_m) (\frakd_m + \frakd_o) r_0,
\end{aligned}
\label{eqn 6.6}
\end{align}
thus
\begin{align*}
\leVe \ddcircDelta \dd{\uflt} \riVe^{n-2,p}
\leq
&
c(n,p) \ufrakd_o r_0 
\\
&
+
c(n,p)  \big(|\os|/r_0+\delta_o\big) 
\left[ (\ud_o^2 + \epsilon \ud_o) \ubfd_m
+
c(n,p) (\ud_o + \epsilon \ud_m) \ubfd_o \right]
r_0
\\
&
+
c(n,p)
(\ud_o^2 + \epsilon \ud_o \ud_m) (\frakd_m + \frakd_o) r_0.
\end{align*}
Collecting estimates for $\dd{\uflt}$, substituting $\ud_o, \ud_m, \ubfd_o, \ubfd_m$ and using the assumptions $\delta \leq \frac{1}{2}$, $(c_o^{\ufrakd_o} + c_o^{\ufrakd_m} + c_m^{\ufrakd_o} + c_m^{\ufrakd_m} + c^{\frakd}) \delta \leq 1$ to simplify formulae, we obtain that
\begin{align*}
\begin{aligned}
\leve \overline{\dd{\uflt}}^{\circg} - \overline{\dd{\ufl{s=0}}}^{\circg} \rive
\leq
&
c(n,p)  \big(|\os|/r_0+\delta_o\big) 
\big[
(\udelta_o^2 + \epsilon \udelta_o) \ufrakd_m
+
(\udelta_o + \epsilon \udelta_m) c_o^{\ufrakd_o} \cdot \ufrakd_o 
\big]
r_0
\\
&
+
c(n,p) (\udelta_o^2 + \epsilon \udelta_o \udelta_m) (\frakd_m + \frakd_o) r_0,
\end{aligned}
\end{align*}
and
\begin{align*}
\begin{aligned}
\leVe \ddslashd \dd{\uflt} \riVe^{n-1,p}
\leq
&
c(n,p) \leVe \ddcircDelta \dd{\uflt} \riVe^{n-2,p}
\\
\leq
&
c(n,p) \ufrakd_o r_0 
+
c(n,p) \big(|\os|/r_0+\delta_o\big) (\udelta_o^2 + \epsilon \udelta_o) \ufrakd_m r_0
\\
&
+
c(n,p)
(\udelta_o^2 + \epsilon \udelta_o \udelta_m) (\frakd_m + \frakd_o) r_0.
\end{aligned}
\end{align*}
Therefore it is sufficient to choose $\delta$ and $c_o^{\ufrakd_o}$, $c_o^{\ufrakd_m}$, $c_m^{\ufrakd_o}$, $c_m^{\ufrakd_m}$, $c^{\frakd}$ satisfying the following inequalities
\begin{align*}
&
c(n,p) < c_o^{\ufrakd_o},
\quad
c(n,p) < c_o^{\ufrakd_m},
\\
&
c(n,p) < c_m^{\ufrakd_m},
\quad
c(n,p) c_o^{\ufrakd_o}< c_m^{\ufrakd_o},
\quad
c(n,p) < c^{\frakd}.
\end{align*}
We can assume that the above constants $c(n,p)$ are the same, thus $c_o^{\ufrakd_o}=c_o^{\ufrakd_m}=c_m^{\ufrakd_m}=c^{\frakd} =2c(n,p), c_m^{\ufrakd_o} = 3c(n,p)^2$ solve the above inequalities. Then choose $\delta$ sufficiently small such that the assumptions on $\delta$ hold. Therefore the bootstrap argument is closed.
\end{proof}

Reviewing proposition \ref{pro 6.1} and its proof, we see that the regularity of $\ddcircDelta \dd{\uflt}$ doesnot surpass the term $\dd{\Xlt}^i \ddpartial_i \big( \ddcircDelta \ufl{2,t} \big)$. This is similar to the issue of the regularity of $\ufl{a,t}$ in proposition \ref{pro 3.3}. 

However if we assume that $\ufl{2,t}$ being constant, i.e. $\Sigmal{2}$ is embedded in the incoming null hypersurface $\uC_{\us}$ of the double null foliation, then $\dd{\Xlt}^i \ddpartial_i \big( \ddcircDelta \ufl{2,t} \big)$ vanishes, thus the regularity of  $\dd{\uflt}$ can be improved. In fact, this improvement is a simple corollary of proposition \ref{pro 3.3}, since if $\ufl{2,t} \equiv \os_2$ being constant, then
\begin{align*}
\ddslashd \dd{\uflt}= - \ddslashd \ufl{1,t},
\end{align*}
thus proposition \ref{pro 3.3} already gives the estimate of $\ddslashd{\dd{\uflt}}$. We state the following improvement for $\dd{\uflt}$ under the additional assumption.
\begin{proposition}\label{pro 6.2}
Under the same setting of proposition \ref{pro 6.1}, if we assume additionally that $\ufl{2,s=0}$ or $\ufl{1,s=0}$ is constant and
\begin{align*}
\leVe \dslashd \dd{\ufl{s=0}} \riVe^{n,p} \leq \ufrakd_o r_0,
\end{align*}
then there exist a small positive constant $\delta$ and constants $c_o,c_m$ all depending on $n,p$, such that if $\epsilon, \udelta_o, \udelta_m, \delta_o$ are suitably bounded that $\epsilon, \udelta_o, \epsilon \udelta_m, \delta_o \leq \delta$, then
\begin{align*}
&
\leVe \ddslashd \dd{\uf} \riVe^{n,p}
\leq
c_o \ufrakd_o r_0,
\\
&
\leve \overline{\dd{\uf}}^{\circg} - \overline{\dd{\ufl{s=0}}}^{\circg} \rive
\leq
c_m (|\os|/r_0 + \delta_o)( \udelta_o + \epsilon \udelta_m) \ufrakd_{o} r_0.
\end{align*}
\end{proposition}
\begin{proof}
Without loss of generality, assume that $\ufl{1,s=0}$ is a constant. Choose $\delta$ in proposition \ref{pro 3.3}. Since $\dslashd \dd{\ufl{s=0}}= \dslashd \ufl{2,s=0}$, we have
\begin{align*}
\leVe \dslashd \ufl{2,s=0} \riVe^{n,p} \leq \min\{ \udelta_o, \ufrakd_o \} r_0.
\end{align*}
Substituting the above bound of $\dslashd \ufl{2,s=0}$ to estimates \eqref{eqn 3.6} in proposition \ref{pro 3.3}, the proposition is proved.
\end{proof}

\begin{remark}
Propositions \ref{pro 6.1} and \ref{pro 6.2} both require the $L^{\infty}$ bounds of the metric components up to $(n+1)$-th order derivatives, the same as proposition \ref{pro 3.3}.
\end{remark}

\section{Perturbation of the outgoing null expansion}\label{sec 7}
In this section, we consider the perturbation of the outgoing null expansion between two spacelike surfaces. 

Adopt the setting of two spacelike surfaces $\Sigmal{a}, a=1,2$ as in section \ref{sec 6}. Let $\uls{a}{( \ddtr \ddchi' )}$ be the outgoing null expansion of $\Sigmal{a}$. The difference of two outgoing null expansions is the corresponding perturbation, which is denoted by $\dd{\ddtr \ddchi'}$
\begin{align*}
\dd{\ddtr \ddchi'} = \uls{a=2}{( \ddtr \ddchi' )} - \uls{a=1}{( \ddtr \ddchi' )}.
\end{align*}
We shall estimate $\dd{\ddtr \ddchi'}$ in terms of the bounds of $\dd{\ufl{s=0}}$, $\dd{f}$. The procedure to obtain the estimate is parallel to the one in section \ref{sec 5}. We first estimate the perturbations of the differential $\uls{a}{(\dslashd \uh )}|_{\Sigmal{a}}$ and the Hessian $\uls{a}{( \circnabla^2 \uh )}|_{\Sigmal{a}}$ in order to estimate $\dd{\ddtr \ddchi'}$.

\subsection{Perturbations of differential and Hessian of parametrisation function $\uhl{a}$}\label{subsec 7.1}
Let $\uls{a}{\ucalH}$ be the incoming null hypersurfaces where $\Sigmal{a}$ is embedded, and $\uhl{a}$ be the parameterisation function of $\uls{a}{\ucalH}$. Since the differential and the Hessian of $\uhl{a}$ are on different surfaces, it is necessary to clarify the precise meaning of their perturbations. Here we use the coordinate $\vartheta$ in the double null coordinate system to match points on $\Sigmal{1}$ and $\Sigmal{2}$: matching two points with the same $\vartheta$ coordinate gives a diffeomorphism between $\Sigmal{a}$. By this diffeomorphism, we can compare the tensors on two surfaces $\Sigmal{a}$.

The perturbations of $\uls{a}{(\dslashd \uh )}|_{\Sigmal{a}}$ and $\uls{a}{( \circnabla^2 \uh )}|_{\Sigmal{a}}$ are their corresponding differences
\begin{align*}
&
\dd{ (\dslashd \uh )|_{\Sigma} } 
=
\uls{a=2}{ (\dslashd \uh )}|_{\Sigmal{2}}
-
\uls{a=1}{ (\dslashd \uh )}|_{\Sigmal{1}},
\quad
\dd{ (\circnabla^2 \uh )|_{\Sigma} } 
=
\uls{a=2}{ (\circnabla^2 \uh )}|_{\Sigmal{2}}
-
\uls{a=1}{ (\circnabla^2 \uh )}|_{\Sigmal{1}},
\end{align*}
We use the rotational vector field derivatives as tools, same as in section \ref{sec 5}. Denote the rotational vector field derivatives of $\uhl{a}$ by
\begin{align*}
\uhl{a}_{\dR,k}= \dR_k \uhl{a},
\quad
\uhl{a}_{\dR,lk} = \dR_l \dR_k \uhl{a}.
\end{align*}
and their perturbations on $\Sigmal{a}$ by
\begin{align*}
\dd{\uh_{\dR,k}|_{\Sigma}} = \uhl{a=2}_{\dR,k}|_{\Sigmal{2}} - \uhl{a=1}_{\dR,k}|_{\Sigmal{1}},
\quad
\dd{\uh_{\dR,lk}|_{\Sigma}} = \uhl{a=2}_{\dR,lk}|_{\Sigmal{2}} - \uhl{a=1}_{\dR,lk}|_{\Sigmal{1}}.
\end{align*}
The rotational vector field components of the perturbations $ \dd{ \big(\dslashd \uh\big)|_{\Sigma}}$ and $ \dd{ \big(\circnabla^2 \uh\big)|_{\Sigma}}$ can be expressed in terms of $\dd{\uh_{\dR,k}|_{\Sigma}}$ and $\dd{ \uh_{\dR,lk}|_{\Sigma}}$ by
\begin{align*}
&
\big[ \dd{ \big(\dslashd \uh\big)|_{\Sigma}} \big]_{\dR,k}
=
\dd{\uh_{\dR,k}|_{\Sigma}},
\\
&
\big[ \dd{ \big(\circnabla^2 \uh\big)|_{\Sigma}} \big]_{\dR,lk}
=
\dd{\uh_{\dR,lk}|_{\Sigma}}  + \epsilon_{lij} x_k x_i \dd{\uh_{\dR,j}|_{\Sigma}}.
\end{align*}

Introduce the family of surfaces $\{\uls{a}{S}_t\}$ on each $\uls{a}{\ucalH}$. Denote the restrictions of $\uhl{a}_{\dR,k}$, $\uhl{a}_{\dR,lk}$ on $\uls{a}{S}_t$ by $\uhl{a,t}_{\dR,i}$, $\uhl{a,t}_{\dR,ij}$. Their differences are denoted by
\begin{align*}
\dd{\uhlt_{\dR,k}} = \uhl{a=2,t}_{\dR,k} - \uhl{a=1,t}_{\dR,k},
\quad
\dd{\uhlt_{\dR,lk}} = \uhl{a=2,t}_{\dR,lk} - \uhl{a=1,t}_{\dR,lk}.
\end{align*}
Since $\uhl{a,t}_{\dR,k}$, $\uhl{a,t}_{\dR,lk}$ satisfy equations \eqref{eqn 5.3}, \eqref{eqn 5.7}, we can derive the equations satisfied by $\dd{\uhlt_{\dR,k}}$, $\dd{\uhlt_{\dR,lk}}$.
\begin{align}
&
\partial_t \dd{\uhlt_{\dR,k}}
=
\Xl{1,t}_{\uh}^i \ddpartial_i \dd{\uhlt_{\dR,k}}
+
\dd{\Xlt_{\uh}^i} \ddpartial_i \big( \uhl{2,t}_{\dR,k} \big)
+
\dd{\relt_{\uh,k}},
\label{eqn 7.1}
\\
&
\partial_t \dd{\uhlt_{\dR,lk}}
=
\Xl{1,t}_{\uh}^i \ddpartial_i \dd{\uhlt_{\dR,lk}}
+
\dd{\Xlt_{\uh}^i} \ddpartial_i \big( \uhl{2,t}_{\dR,lk} \big)
+
\dd{\relt_{\uh,lk}},
\label{eqn 7.2}
\end{align}
where $\dd{\Xlt_{\uh}}$, $\dd{\relt_{\uh,k}}$, $\dd{\relt_{\uh,lk}}$ are the perturbations of the corresponding quantities
\begin{align*}
&
\dd{\Xlt_{\uh}^i} = \Xl{a=2,t}_{\uh}^i - \Xl{a=1,t}_{\uh}^i,
\\
&
\dd{\relt_{\uh,k}} = \rel{a=2,t}_{\uh,k} - \rel{a=1,t}_{\uh,k},
\quad
\dd{\relt_{\uh,lk}} = \rel{a=2,t}_{\uh,lk} - \rel{a=1,t}_{\uh,lk}.
\end{align*}
We shall integrate equations \eqref{eqn 7.1} \eqref{eqn 7.2} to obtain the estimates of $\dd{\uhlt_{\dR,k}}$ and $\dd{\uhlt_{\dR,lk}}$. Before stating the results, it is worth to point out the roles of terms $\dd{\Xlt_{\uh}^i} \ddpartial_i \big( \uhl{2,t}_{\dR,k} \big)$ and $\dd{\Xlt_{\uh}^i} \ddpartial_i \big( \uhl{2,t}_{\dR,lk} \big)$: due to these two terms, the regularities of $\dd{\uhlt_{\dR,k}}$ and $\dd{\uhlt_{\dR,lk}}$ are one order less than $\uhl{a,t}_{\dR,k}$ and $\uhl{a,t}_{\dR,lk}$ respectively. This effect is similar to the one caused by the term $\dd{\Xlt}^i \ddpartial_i \big( \ddcircDelta \ufl{2,t} \big)$ in equation \eqref{eqn 6.2}.

\begin{proposition}\label{pro 7.1}
Let $\Sigmal{a}, a=1,2$ be two spacelike surfaces in $(M,g)$. Suppose that $\Sigmal{a}$ has the second parameterisation $(\ufl{a,s=0}, \fl{a})$. Assume that the parameterisation functions satisfy the following estimates
\begin{align*}
&
\leVe \dslashd \ufl{a,s=0} \riVe^{n+1,p} 
\leq
\udelta_o r_0,
\quad
\leve \overline{\ufl{a,s=0}}^{\circg} \rive
\leq 
\udelta_m r_0,
\\
&
\Vert \ddslashd \fl{a} \Vert^{n+1,p} \leq  \delta_o r_0,
\quad
\overline{\fl{a}}^{\circg} = \os_a,
\quad
\os = \max \{ |\os_1|, |\os_2| \},
\end{align*}
and the perturbation functions satisfy
\begin{align*}
&
\leVe \dslashd \dd{\ufl{s=0}} \riVe^{n,p}  \leq  \ufrakd_o r_0,
\quad
\leve \overline{\dd{\ufl{s=0}}}^{\circg} \rive \leq \ufrakd_m r_0,
\\
&
\Vert \ddslashd \dd{f} \Vert^{n,p} \leq \frakd_o r_0,
\quad
\leve \overline{\dd{f}}^{\circg} \rive \leq \frakd_m r_0,
\end{align*}
where $n\geq 2, p>2$ or $n\geq 3, p>1$.

There exist a small positive constant $\delta$ and constants $c'_{\uh}$ both depending on $n,p$, such that if $\epsilon, \udelta_o, \udelta_m, \delta_o, \ufrakd_o, \ufrakd_m, \frakd_o, \frakd_m$ are suitably bounded that
\begin{align*}
\epsilon,\ \udelta_o,\ \epsilon \udelta_m,\ \delta_o,\ \ufrakd_o,\ \epsilon \ufrakd_m,\ \frakd_o,\ \frakd_m \leq \delta,
\end{align*}
then the perturbations $\dd{\big( \dslashd \uh \big)|_{\Sigma}}$ and $\dd{\big( \circnabla^2 \uh \big)|_{\Sigma}}$ satisfy
\begin{align}
\begin{aligned}
&
\leVe \dd{( \dslashd \uh )|_{\Sigma}} \riVe^{n,p},
\leVe \dd{( \circnabla^2 \uh )|_{\Sigma}} \riVe^{n-1,p}
\\
&
\leq
c'_{\uh} \ufrakd_o r_0
+
c'_{\uh} (|\os|/r_0+\delta_o) (\epsilon \udelta_o + \udelta_o^2) \ufrakd_m r_0
+
c'_{\uh} ( \epsilon \udelta_m \udelta_o + \udelta_o^2) (\frakd_m + \frakd_o) r_0.
\end{aligned}
\label{eqn 7.3}
\end{align}
\end{proposition}
The proof is rather technical and involved, thus we present it in appendix \ref{appen pro 7.1}. Here we just remark that proposition \ref{pro 7.1} requires the $L^{\infty}$ bounds of the metric components up to $(n+2)$-th order derivatives.

Similar to proposition \ref{pro 6.2}, if one of $\ufl{a,s=0}$ is constant, then we can improve the regularities of $\dd{( \dslashd \uh )|_{\Sigma}}$, $\dd{( \circnabla^2 \uh )|_{\Sigma}}$.
\begin{proposition}\label{pro 7.2}
Under the same setting of proposition \ref{pro 7.1}, if we assume additionally that $\ufl{2,s=0}$ or $\ufl{1,s=0}$ is constant and 
\begin{align*}
\leVe \dslashd \dd{\ufl{s=0}} \riVe^{n+1,p} \leq \ufrakd_o r_0,
\end{align*}
then there exist a small positive constant $\delta$ and constant $c'_{\uh}$ both depending on $n,p$, such that
\begin{align*}
\leVe \dd{( \dslashd \uh )|_{\Sigma}} \riVe^{n+1,p},
\leVe \dd{( \circnabla^2 \uh )|_{\Sigma}} \riVe^{n,p}
\leq c'_{\uh} \ufrakd_o r_0.
\end{align*}
\end{proposition}
\begin{proof}
Without loss of generality, we can assume that $\ufl{1,s=0}$ is constant, then $\uhl{1}$ is also constant, thus
\begin{align*}
\dd{( \dslashd \uh )|_{\Sigma}} = \uls{a=2}{\big( \dslashd \uh \big)}|_{\Sigmal{2}},
\quad
\dd{( \circnabla^2 \uh )|_{\Sigma}} = \uls{a=2}{\big( \circnabla^2 \uh \big)}|_{\Sigmal{2}},
\end{align*}
Then the proposition follows from propositions \ref{pro 5.1}, \ref{pro 5.3} with the condition
\begin{align*}
\leVe \dslashd \ufl{2,s=0} \riVe^{n+1,p} \leq \min\{ \udelta_o,\ufrakd_o \} r_0.
\end{align*}
The proposition also requires the $L^{\infty}$ bounds of the metric components up to the $(n+2)$-th order derivatives.
\end{proof}

\subsection{Estimate of perturbation of the outgoing null expansion}
We already obtain the estimates for $\dd{\uf}$ in proposition \ref{pro 6.1} and for $\dd{( \dslashd \uh )|_{\Sigma}}$, $\dd{( \circnabla^2 \uh )|_{\Sigma}}$ in proposition \ref{pro 7.1}. In this subsection, we use these estimates to estimate the perturbation of the outgoing null expansion.

We adopt the setting of proposition \ref{pro 7.1}, and the notations $\ud_o, \ud_m$ in formulae \eqref{eqn 3.7}, $\uslashd_m$ in \eqref{eqn 5.8}, $\ubfd_o, \uslashbfd_m, \ubfd_m$ in \eqref{eqn 6.4}, $\ubfd_{\uh}$ in \eqref{appen pro 7.1 eqn 1}. Introduce another notation $\ud_{\uh}$ that
\begin{align}
\ud_{\uh}= \max\{ c_{\uh}, c_{\uh,2} \} \udelta_o.
\label{eqn 7.4}
\end{align}
Then we assume that the following estimates hold:
\begin{align}
\begin{aligned}
&
\big\Vert \ddslashd \ufl{a} \big\Vert^{n,p} \leq \ud_o r_0,
\quad
\big| \overline{\ufl{a}}^{\circg} - \overline{\ufl{a,s=0}}^{\circg} \big| \leq \uslashd_m r_0,
\\
&
\big\Vert \ddslashd \dd{\uf} \big\Vert^{n-1,p} \leq \ubfd_o r_0,
\quad
\big| \overline{\dd{\uf}}^{\circg} - \overline{\dd{\ufl{s=0}}}^{\circg} \big| \leq \uslashbfd_m r_0,
\\
&
\big\Vert \uls{a}{( \dslashd \uh )|_{\Sigmal{a}}} \big\Vert^{n+1,p},
\big\Vert \uls{a}{( \circnabla^2 \uh )|_{\Sigmal{a}}} \big\Vert^{n,p}
\leq
\ud_{\uh} r_0,
\\
&
\big\Vert \dd{( \dslashd \uh )|_{\Sigma}} \big\Vert^{n,p},
\big\Vert \dd{( \circnabla^2 \uh)|_{\Sigma}} \big\Vert^{n-1,p}
\leq
\ubfd_{\uh} r_0.
\end{aligned}
\label{eqn 7.5}
\end{align}
We estimate the perturbation of the outgoing null expansion with the above estimates.

\begin{proposition}\label{pro 7.3}
Under the setting of proposition \ref{pro 7.1} and assuming estimates \eqref{eqn 7.5},
\begin{enumerate}[label=\alph*.]
\item \label{pro 7.3.a}
the perturbation of the first order main part $\dd{\lo{\ddtr \ddchi'}}$ satisfies the estimate
\begin{align*}
\big\Vert \dd{\lo{\ddtr \ddchi'}} \big\Vert^{n-1,p}
\leq
\frac{c(n,p)}{r_0} ( \frakd_m+ \frakd_o )
+
\frac{c(n,p)(|\os|/r_0+\delta_o)}{r_0} (\ubfd_m + \ubfd_o ),
\end{align*}
\item \label{pro 7.3.b}
the perturbations of the high order remainder parts $\dd{\hir{1}{\ddtr \ddchi'}}$, $\dd{\hir{2}{\ddtr \ddchi'}}$ satisfy the estimates
\begin{align*}
\big\Vert \dd{\hir{1}{\ddtr \ddchi'}} \big\Vert^{n-1,p}
\leq&
\frac{c(n,p) \epsilon}{r_0} (\frakd_m+ \frakd_o + \ubfd_m + \ubfd_o +\ubfd_{\uh}),
\\
\big\Vert \dd{\hir{2}{\ddtr \ddchi'}} \big\Vert^{n-1,p}
\leq&
\frac{c(n,p) (\epsilon + \delta_o + \ud_{\uh}) \delta_o}{r_0} \frakd_m
+
\frac{c(n,p) (\epsilon + \delta_o + \ud_{\uh} |\os|/r_0)}{r_0} \frakd_o
\\
&
+
\frac{c(n,p) (\epsilon + \delta_o + \ud_{\uh} |\os|/r_0) \delta_o}{r_0} (\ubfd_m+\ubfd_o)
\\
&
+
\frac{c(n,p) (\epsilon + \delta_o + |\os|/r_0) \delta_o}{r_0} \ubfd_{\uh}.
\end{align*}
\end{enumerate}
\end{proposition}
\begin{proof}[Proof of proposition \ref{pro 7.3}.\ref{pro 7.3.a}]
Recall that
\begin{align*}
\lo{\ddtr \ddchi'}
=
\tr \chi'_S|_{\Sigma} - 2 ( r_S|_{\Sigma} )^{-2} \circDelta f,
\end{align*}
therefore
\begin{align*}
\dd{\lo{\ddtr \ddchi'}} = \dd{\tr \chi'_{S}|_{\Sigma}} - 2 ( r_S|_{\Sigmal{2}} )^{-2} \circDelta \dd{f} - 2 \dd{ ( r_S|_{\Sigma} )^{-2} } \circDelta ( \fl{a=1} ).
\end{align*}
By the estimates of $\dd{\uf}$ and $\dd{f}$, we have
\begin{align*}
&
\big\Vert \circDelta \dd{f}  \big\Vert^{n-1,p} 
\leq
\frakd_o r_0,
\\
&
\left.
\begin{aligned}
\big\Vert \dd{\tr \chi'_{S}|_{\Sigma}} \big\Vert^{n,p}
\\
\big\Vert \dd{ \big( r_S|_{\Sigma}\big)^{-2} } \big\Vert^{n,p}
\end{aligned}
\right\}
\leq
\frac{c(n,p)}{r_0} (\frakd_m + \frakd_o) + \frac{c(n,p)(|\os|/r_0 + \delta_o)}{r_0} ( \ubfd_m + \ubfd_o ),
\end{align*}
Then the estimate of $\dd{\lo{\ddtr \ddchi'}}$ follows.
\end{proof}

\begin{proof}[Proof of proposition \ref{pro 7.3}.\ref{pro 7.3.b}]
Recall the estimates of $\hir{1}{\ddtr \ddchi'}$, $\hir{2}{\ddtr \ddchi'}$ in proposition \ref{pro 5.4}
\begin{align*}
&
\leVe \hir{1}{\ddtr \ddchi'} \riVe^{n,p}
\leq
c(n,p) \epsilon r_0^{-1},
\\
&
\leVe \hir{2}{\ddtr \ddchi'} \riVe^{n,p}
\leq
c(n,p) ( \epsilon + \delta_o +  \ud_{\uh} |\os|/r_0) \delta_o r_0^{-1}.
\end{align*}
Then applying rules \eqref{appen pro 7.1 eqn 2} \eqref{appen pro 7.1 eqn 3} in appendix \ref{appen pro 7.1} to the right hand side of the above estimates, we obtain the right hand side of the estimates of $\dd{\hir{1}{\ddtr \ddchi'}}$, $\dd{\hir{2}{\ddtr \ddchi'}}$.

We shall clarify the regularities of $\dd{\hir{1}{\ddtr \ddchi'}}$, $\dd{\hir{2}{\ddtr \ddchi'}}$ in their estimates. Both of them involve $\dd{\ddslashg}$ and $\dd{\ddchi'}$.

In $\dd{\ddslashg}$, the terms with the worst regularity are 
\begin{enumerate}[label=\emph{\roman*}.]
\item 
$\dd{( \dslashd \uh )|_{\Sigma}}$ in $\dd{\db}$, which is in the Sobolev space $\mathrm{W}^{n,p}$,
\item 
$\dd{\slashg}$ and $\dd{b}$ which are both in the Sobolev space $\mathrm{W}^{n,p}$, since $\dd{\uf}$ is in $\mathrm{W}^{n,p}$.
\end{enumerate}

In $\dd{\ddchi'}$, the terms with the worst regularity are
\begin{enumerate}[label=\emph{\roman*}.]
\item
$\dd{\circnabla^2 \uh}$ in $\dd{\duchi}$, $\dd{\dslashnabla \db}$, $\dd{\dpartial_{\us} \db}$, which is in the Sobolev space $\mathrm{W}^{n-1,p}$,
\item
$\dd{\circnabla^2 f}$ which is in $\mathrm{W}^{n-1,p}$, because $\dd{f}$ is in $\mathrm{W}^{n+1,p}$ in the assumption of the proposition (same as in proposition \ref{pro 7.1}).
\end{enumerate}
Therefore $\dd{\hir{1}{\ddtr \ddchi'}}$, $\dd{\hir{2}{\ddtr \ddchi'}}$ are in the Sobolev space $\mathrm{W}^{n-1,p}$.

We determine up to which order derivatives, the $L^{\infty}$ bounds of the metric components and the structure coefficients are required.

For the metric components:
\begin{enumerate}[label=\emph{\roman*}.]
\item the $L^{\infty}$ bounds of the metric components up to $(n+2)$-th order derivatives are sufficient for estimates \eqref{eqn 7.5}, by proposition \ref{pro 3.3}, \ref{pro 6.1} and \ref{pro 7.1},
\item The Sobolev norms $\Vert \cdot \Vert^{n-1,p}$ of $\dd{\slashg}$, $\dd{b}$ and 
\begin{align*}
\dd{\circnabla \slashg} \text{ in } \dd{\circtriangle}, 
\quad
\dd{\big(\circnabla, \partial_{\us}, \partial_s \big) \big(b, \Omega^2 \slashg^{-1} \big)} \text{ in } \dd{\dslashnabla \db} \text{ and }\dd{\dpartial_{\us} \db}
\end{align*}
are required. Thus the $L^{\infty}$ bounds of the metric components up to $(n+1)$-th order derivatives are required.
\end{enumerate}
Therefore, the $L^{\infty}$ bounds of the metric components up to $(n+2)$-th order derivatives are required.

For the structure coefficients: the Sobolev norms $\Vert \cdot \Vert^{n-1,p}$ of $\dd{(\chi', \uchi, \eta, \uomega)}$ are required, therefore the $L^{\infty}$ bounds of the structure coefficients up to $n$th order derivatives are required.
\end{proof}

\begin{remark}
Note that the requirements on the order of derivatives of the metric components and structure coefficients for proposition \ref{pro 7.3} are the same as the requirements for the estimates of $\hir{1}{\ddtr \ddchi'}, \hir{2}{\ddtr \ddchi'}$ in proposition \ref{pro 5.4}.
\end{remark}

\subsection{Improved estimate of perturbation of the outgoing null expansion}
The regularities in the estimates of $\dd{\uf}$, $\dd{( \dslashd \uh )|_{\Sigma}}$, $\dd{( \circnabla^2 \uh )|_{\Sigma}}$ are improved in propositions \ref{pro 6.2}, \ref{pro 7.2}. Similar improvement is available for the perturbation of the outgoing null expansion. We present it in this subsection.

In the following, we adopt the setting of proposition \ref{pro 7.2}. Introduce the notations
\begin{align}
\ubfd'_o = c_o \ufrakd_o,
\quad
\uslashbfd'_m = c_m (|\os|/r_0 + \delta_o)( \udelta_o + \epsilon \udelta_m ) \ufrakd_o,
\quad
\ubfd'_m = \ufrakd_m + \uslashbfd'_m
\quad
\ubfd'_{\uh} = c'_{\uh} \ufrakd_o.
\label{eqn 7.6}
\end{align}
We have the following improved estimates on $\dd{\uf}$, $\dd{( \dslashd \uh )|_{\Sigma}}$, $\dd{( \circnabla^2 \uh )|_{\Sigma}}$ comparing to \eqref{eqn 7.5} in the previous subsection 
\begin{align}
\begin{aligned}
&
\big\Vert \ddslashd \dd{\uf} \big\Vert^{n,p} \leq \ubfd'_o r_0,
\quad
\big| \overline{\dd{\uf}}^{\circg} - \overline{\dd{\ufl{s=0}}}^{\circg} \big| \leq \uslashbfd'_m r_0,
\quad
\text{from proposition \ref{pro 6.2},}
\\
&
\big\Vert \dd{\big( \dslashd \uh \big)|_{\Sigma}} \big\Vert^{n+1,p},
\big\Vert \dd{\big( \circnabla^2 \uh \big)|_{\Sigma}} \big\Vert^{n,p}
\leq
\ubfd'_{\uh} r_0,
\quad
\text{from proposition \ref{pro 7.2}.}
\end{aligned}
\label{eqn 7.7}
\end{align}
We state the improved estimate of the perturbation of the outgoing null expansion.
\begin{proposition}\label{pro 7.5}
Given the setting of proposition \ref{pro 7.2}, assume estimates \eqref{eqn 7.7} and
\begin{align*}
\big\Vert \ddslashd \dd{f} \big\Vert^{n+1,p} \leq \frakd_o r_0,
\end{align*}
 then
\begin{enumerate}[label=\alph*.]
\item \label{pro 7.5.a}
the perturbation of the first order main part $\dd{\lo{\ddtr \ddchi'}}$ satisfies the estimate
\begin{align*}
\big\Vert \dd{\lo{\ddtr \ddchi'}} \big\Vert^{n,p}
\leq
\frac{c(n,p)}{r_0} ( \frakd_m+ \frakd_o )
+
\frac{c(n,p)(|\os|/r_0+\delta_o)}{r_0} (\ubfd'_m + \ubfd'_o );
\end{align*}
\item \label{pro 7.5.b}
the perturbations of the high order remainder parts $\dd{\hir{1}{\ddtr \ddchi'}}$, $\dd{\hir{2}{\ddtr \ddchi'}}$ satisfy the estimates
\begin{align*}
\big\Vert \dd{\hir{1}{\ddtr \ddchi'}} \big\Vert^{n,p}
\leq&
\frac{c(n,p) \epsilon}{r_0} (\frakd_m+ \frakd_o + \ubfd'_m + \ubfd'_o +\ubfd'_{\uh}),
\\
\big\Vert \dd{\hir{2}{\ddtr \ddchi'}} \big\Vert^{n,p}
\leq&
\frac{c(n,p) (\epsilon + \delta_o + \ud_{\uh}) \delta_o}{r_0} \frakd_m
+
\frac{c(n,p) (\epsilon + \delta_o + \ud_{\uh} |\os|/r_0)}{r_0} \frakd_o
\\
&
+
\frac{c(n,p) (\epsilon + \delta_o + \ud_{\uh} |\os|/r_0) \delta_o}{r_0} (\ubfd'_m+\ubfd'_o)
\\
&
+
\frac{c(n,p) (\epsilon + \delta_o + |\os|/r_0) \delta_o}{r_0} \ubfd'_{\uh}.
\end{align*}
\end{enumerate}
\end{proposition}
Note that the form of the improved estimates is almost identical to the one of estimates in proposition \ref{pro 7.3}, while the regularities are all improved by one order. The proof follows the same route as the proof of proposition \ref{pro 7.3}, thus we omit it here.
\begin{remark}
We determine up to which order derivatives, the $L^{\infty}$ bounds of the metric components and the structure coefficients are required for the improved estimates.

For the metric components, the $L^{\infty}$ bounds up to $(n+2)$-th order derivatives are required, the same as proposition \ref{pro 7.3}.

For the structure coefficients, the $L^{\infty}$ bounds up to $(n+1)$-th order derivatives are required, one order higher than proposition \ref{pro 7.3}, since the Sobolev norms $\Vert \cdot \Vert^{n,p}$ of $\dd{(\chi', \uchi, \eta, \uomega)}$ are required for the improved estimates.
\end{remark}

\section{Linearised perturbation of parameterisation of spacelike surface}\label{sec 8}
In this section, we shall construct a linearised perturbation of the parameterisation of spacelike surfaces. 

We adopt the notations in section \ref{sec 6}. Let $\Sigmal{a}, a=1,2$ be two spacelike surfaces in $(M,g)$. In section \ref{sec 6}, we obtained the perturbation $\dd{\uf}=\ufl{2}-\ufl{1}$ from the perturbations $\dd{\ufl{s=0}}, \dd{f}$. In this section, we will construct an appropriate linearisation $\bdd{\uf}$ for the perturbation $\dd{\uf}$, and estimate the corresponding error $\er{\uf} = \dd{\uf} - \bdd{\uf}$.

\subsection{Equation of linearised perturbation of the first parameterisation function}
Recall equations \eqref{eqn 6.1} \eqref{eqn 6.2} of the perturbation functions $\dd{\uflt}$
\begin{align}
&
\partial_t \dd{\uflt} 
= 
\dd{\Flt},
\tag{\ref{eqn 6.1}}
\\
&
\partial_t \big( \ddcircDelta \dd{\uflt} \big)
=
\Xl{1,t}^i \ddpartial_i \big( \ddcircDelta \dd{\uflt} \big)
+\dd{\Xlt}^i \ddpartial_i \big( \ddcircDelta \ufl{2,t} \big) 
+ \dd{\relt}.
\tag{\ref{eqn 6.2}}
\end{align}
We shall construct a system of linear equations from equations \eqref{eqn 6.1} \eqref{eqn 6.2}, then define the linearised perturbation function $\bdd{\uf}$ as the solution of the linear system. Introduce a family of functions $\bdd{\uflt}$ where $t\in[0,1]$, and the following linear system of $\bdd{\uflt}$
\begin{align}
&
\overline{\bdd{\uflt}}^{\circg} = \overline{\bdd{\ufl{s=0}}}^{\circg},
\label{eqn 8.1}
\\
&
\partial_t \big( \ddcircDelta \bdd{\uflt} \big) 
= 
\Xl{1,t}^i \ddpartial_i \big( \ddcircDelta \bdd{\uflt} \big) 
- \overline{\Xl{1,t}^i \ddpartial_i \big( \ddcircDelta \bdd{\uflt} \big)}^{\circg}.
\label{eqn 8.2}
\end{align}
We set the initial data of the above system to be
\begin{align*}
\bdd{\ufl{t=0}}= \dd{\ufl{s=0}} = \ufl{2,s=0} - \ufl{1,s=0}.
\end{align*}
Then define the linearised perturbation $\bdd{\uf}$ as the solution $\bdd{\ufl{t=1}}$. We summarise the above construction in the following definition.
\begin{definition}\label{def 8.1}
Let $\Sigmal{a},a=1,2$ be two spacelike surfaces in $(M,g)$. Suppose $\Sigmal{a}$ has the second parameterisation $( \ufl{a,s=0}, \fl{a})$ and the first parameterisation $(\ufl{a}, \fl{a})$. We define the linearised perturbation of the first parameterisation, $(\bdd{\uf}, \bdd{f})$, from $\Sigmal{1}$ to $\Sigmal{2}$ as follows:
\begin{align*}
\bdd{f} = \dd{f} = \fl{2}-\fl{1},
\quad
\bdd{\uf} = \bdd{\ufl{t=1}},
\end{align*}
where $\bdd{\uflt}$ solves equations \eqref{eqn 8.1}, \eqref{eqn 8.2} with the initial data $\bdd{\ufl{t=0}} = \dd{\ufl{s=0}}$.
\end{definition}

We explain briefly why the above linearised perturbation $\bdd{\uf}$ is appropriate. If we assume formally that the size of $\dd{\uflt}$, $\dd{f}$ is $\mathbf{d}$ and the size of $\epsilon \ufl{a,t}$, $\ddslashd \ufl{a,t}$, $\ddslashd\fl{a}$ is $\delta$, then the following terms in equation \eqref{eqn 6.2}
\begin{align*}
\dd{\Xlt}^i \ddpartial_i \big( \ddcircDelta \ufl{2,t} \big),
\quad
\dd{\relt},
\end{align*}
will be of size $\delta \cdot \mathbf{d}$, which is one magnitude smaller than the size of the perturbation. Thus we omit these terms when constructing the linear system for $\bdd{\uflt}$ by allowing an error of size $\delta \mathbf{d}$. The above explanation will be made rigorous when we estimate the error $\er{\uf}=\dd{\uf}-\bdd{\uf}$ in the next subsection.

We estimate the solution $\bdd{\uflt}$ of equations \eqref{eqn 8.1} \eqref{eqn 8.2}.
\begin{lemma}\label{lem 8.2}
Assume that the vector field $\Xl{1,t}$ satisfies the estimate
\begin{align*}
\leVe \Xl{1,t} \riVe^{n,p} \leq k,
\end{align*}
where $n\geq 2,p>2$ or $n\geq 3,p>1$. There exist constants $c(p)$, $c(n,p)$ such that the solution $\bdd{\uflt}$ of equation \eqref{eqn 8.2} satisfies the estimate
\begin{align*}
\leVe \ddcircDelta \bdd{\uflt} \riVe^{m,p} \leq c(p) \exp(c(n,p)k) \leVe \ddcircDelta \bdd{\ufl{t=0}} \riVe^{m,p}.
\end{align*}
for all $0\leq m \leq n$.
\end{lemma}
\begin{proof}

Consider a slightly different equation
\begin{align*}
\partial_t G_t = \Xl{1,t}^i \ddpartial_i G_t.
\end{align*}
Then by Gronwall's inequality, we have $\leVe G_t \riVe^{m,p} \leq  \exp(c(n,p)k) \leVe G_{t=0} \riVe^{m,p}$. Let $G_t$ be the solution with the initial condition $G_{t=0} = \ddcircDelta \bdd{\ufl{t=0}}$. Then we have $\ddcircDelta \bdd{\uflt} = G_t - \overline{G_t}^{\circg}$. Therefore
\begin{align*}
\leVe \ddcircDelta \bdd{\uflt} \riVe^{m,p} \leq  \leVe G_t \riVe^{m,p} + \leVe \overline{G_t}^{\circg} \riVe^{m,p} \leq c(p) \leVe G_t \riVe^{m,p} \leq c(p) \exp(c(n,p) k) \leVe G_{t=0} \riVe^{m,p}. 
\end{align*} 
The lemma is proved.
\end{proof}

\subsection{Estimate of error of the linearised parametrisation perturbation}
The error of the linearised perturbation $\bdd{\uflt}$ is the difference between $\dd{\uflt}$ and $\bdd{\uflt}$,
\begin{align*}
\er{\uflt} = \dd{\uflt} - \bdd{\uflt},
\quad
\er{\uf} = \dd{\uf} - \bdd{\uf}.
\end{align*}
Taking the difference of equations \eqref{eqn 6.2} and \eqref{eqn 8.2}, we derive the equation for $\er{\uflt}$,
\begin{align}
\partial_t \big( \ddcircDelta \er{\uflt} \big)
=
\Xl{1,t}^i \ddpartial_i \big( \ddcircDelta \er{\uflt} \big)
- \overline{\Xl{1,t}^i \ddpartial_i \big( \ddcircDelta \bdd{\uflt} \big)}^{\circg}
+\dd{\Xlt}^i \ddpartial_i \big( \ddcircDelta \ufl{2,t} \big) 
+ \dd{\relt}.
\label{eqn 8.3}
\end{align}
Adopting the notations $\ud_o, \uslashd_m, \ud_m, \ubfd_o, \uslashbfd_m, \ubfd_m$ in formulae \eqref{eqn 3.7} \eqref{eqn 5.8} \eqref{eqn 6.4}, we assume the following estimates in \eqref{eqn 7.5},
\begin{align}
\begin{aligned}
&
\big\Vert \ddslashd \ufl{a} \big\Vert^{n,p} \leq \ud_o r_0,
\quad
\big| \overline{\ufl{a}}^{\circg} - \overline{\ufl{a,s=0}}^{\circg} \big| \leq \uslashd_m r_0,
\\
&
\big\Vert \ddslashd \dd{\uf} \big\Vert^{n-1,p} \leq \ubfd_o r_0,
\quad
\big| \overline{\dd{\uf}}^{\circg} - \overline{\dd{\ufl{s=0}}}^{\circg} \big| \leq \uslashbfd_m r_0.
\end{aligned}
\tag{\ref{eqn 7.5}}
\end{align}
We estimate the error $\er{\uflt}$ by integrating equation \eqref{eqn 8.3}.

\begin{proposition}\label{pro 8.3}
Let $\Sigmal{a}, a=1,2$ be two spacelike surfaces in $(M,g)$. Assume that $\Sigmal{a}$ has the second parameterisation $(\ufl{a,s=0}, \fl{a})$. Suppose that the parameterisation functions satisfy the estimates
\begin{align*}
&
\leVe \dslashd \ufl{a,s=0} \riVe^{n,p} 
\leq
\udelta_o r_0,
\quad
\leve \overline{\ufl{a,s=0}}^{\circg} \rive
\leq 
\udelta_m r_0,
\\
&
\leVe \ddslashd \fl{a} \riVe^{n+1,p} \leq  \delta_o r_0,
\quad
\overline{\fl{a}}^{\circg} = \os_a,
\quad
\os = \max \{ |\os_1|, |\os_2| \},
\end{align*}
and the perturbation functions satisfy
\begin{align*}
&
\leVe \dslashd \dd{\ufl{s=0}} \riVe^{n-1,p}  \leq  \ufrakd_o r_0,
\quad
\leve \overline{\dd{\ufl{s=0}}}^{\circg} \rive \leq \ufrakd_m r_0,
\\
&
\leVe \ddslashd \dd{f} \riVe^{n,p} \leq \frakd_o r_0,
\quad
\leve \overline{\dd{f}}^{\circg} \rive \leq \frakd_m r_0,
\end{align*}
where $n\geq 2, p>2$ or $n\geq 3, p>1$.

There exists a small positive constant $\delta$ depending on $n,p$, such that if $\epsilon$, $\udelta_o$, $\udelta_m$, $\delta_o$, $\ufrakd_o$, $\ufrakd_m$, $\frakd_o$, $\frakd_m$ are suitably bounded that $\epsilon,\ \udelta_o,\ \epsilon \udelta_m,\ \delta_o,\ \ufrakd_o,\ \epsilon \ufrakd_m,\ \frakd_o,\ \frakd_m \leq \delta$, then the error $\er{\uf}$ satisfies the following estimates
\begin{align}
\begin{aligned}
\leve \overline{\er{\uf}}^{\circg} \rive,
\leVe \ddslashd \er{\uf} \riVe^{n-1,p}
\leq
&
c(n,p) \big(|\os|/r_0+\delta_o\big) 
\left[ (\udelta_o^2 + \epsilon \udelta_o) \ufrakd_m
+
(\udelta_o + \epsilon \udelta_m) \ufrakd_o \right] r_0
\\
&
+
c(n,p) (\udelta_o^2 + \epsilon \udelta_o \udelta_m) (\frakd_m + \frakd_o) r_0.
\end{aligned}
\label{eqn 8.4}
\end{align}
\end{proposition}
\begin{proof}
We choose $\delta$ sufficiently small such that propositions \ref{pro 3.3} and \ref{pro 6.1} hold. The estimate of $\overline{\er{\uf}}^{\circg}$ follows directly from estimates \eqref{eqn 7.5} and $\overline{\er{\uf}}^{\circg} = \overline{\dd{\uf}}^{\circg} - \overline{\dd{\ufl{s=0}}}^{\circg}$.

By estimate \eqref{eqn 3.8} and lemma \ref{lem 8.2}, we have
\begin{align*}
&
\Vert \Xl{1,t} \Vert^{n,p} \leq c(n,p) (|\os|/r_0 + \delta_o)(\ud_o + \epsilon \ud_m) \leq c(n,p),
\\
&
\begin{aligned}
\leve \overline{\Xl{1,t}^i \ddpartial_i \big( \ddcircDelta \bdd{\uflt} \big)}^{\circg} \rive
=&
\leve \overline{ \ddcircdiv \Xl{1,t} \cdot  \ddcircDelta \bdd{\uflt} }^{\circg} \rive
\leq
\leVe \ddcircdiv \Xl{1,t} \riVe_{L^{\infty}} \cdot \leVe \ddcircDelta \bdd{\uflt}  \riVe_{L^1}
\\
\leq&
c(n,p) (|\os|/r_0 + \delta_o)(\ud_o + \epsilon \ud_m) \ufrakd_o r_0.
\end{aligned}
\end{align*} 
Together with estimate \eqref{eqn 6.6} of the integral
$\int_0^t \leVe \dd{\Xlt}^i \ddpartial_i \big( \ddcircDelta \ufl{2,t} \big) + \dd{\relt} \riVe^{n-2,p} \d t$,
we obtain that
\begin{align*}
\leVe \ddslashd \er{\uf} \riVe^{n-1,p}
\leq
&
c(n,p) (|\os|/r_0 + \delta_o)(\ud_o + \epsilon \ud_m) \ufrakd_o r_0
\\
&+
c(n,p) \big(|\os|/r_0+\delta_o\big) 
\left[ (\ud_o^2 + \epsilon \ud_o) \ubfd_m
+
(\ud_o + \epsilon \ud_m) \ubfd_o \right] r_0
\\
&
+
c(n,p) (\ud_o^2 + \epsilon \ud_o \ud_m) (\frakd_m + \frakd_o) r_0.
\end{align*}
Substituting $\ud_o, \ud_m, \ubfd_o, \uslashbfd_m, \ubfd_m$, we prove the proposition.
\end{proof}
\begin{remark}
The assumption of proposition \ref{pro 8.3} is the same as proposition \ref{pro 6.1}. Proposition \ref{pro 8.3} requires the $L^{\infty}$ bounds of the metric components up to $(n+1)$-th order derivatives, since propositions \ref{pro 3.3} and \ref{pro 6.1} require so.
\end{remark}

\subsection{Improved estimate of error of the linearised parametrisation perturbation}
We show that the estimate of $\er{\uf}$ can be improved given the additional condition that $\ufl{2,s=0}$ is constant. It is similar to the improvement for the estimate of $\dd{\uf}$ in proposition \ref{pro 6.2}.
\begin{proposition}\label{pro 8.5}
Under the setting of proposition \ref{pro 8.3}, we assume additionally that $\ufl{2,s=0}$ is constant and
\begin{align*}
\leVe \dslashd \dd{\ufl{s=0}} \riVe^{n,p} \leq \ufrakd_o r_0.
\end{align*}

There exists a small positive constant $\delta$ depending on $n,p$, that if $\epsilon$, $\udelta_o$, $\udelta_m$, $\delta_o$, $\ufrakd_o$, $\ufrakd_m$, $\frakd_o$, $\frakd_m$ are suitably bounded such that
$\epsilon,\ \udelta_o,\ \epsilon \udelta_m,\ \delta_o,\ \ufrakd_o,\ \epsilon \ufrakd_m,\ \frakd_o,\ \frakd_m \leq \delta$,
then the error $\er{\uf}$ satisfies the following improved estimates comparing with proposition \ref{pro 8.3}
\begin{align}
\begin{aligned}
\leve \overline{\er{\uf}}^{\circg} \rive,
\leVe \ddslashd \er{\uf} \riVe^{n,p}
\leq
&
c(n,p) \big(|\os|/r_0+\delta_o\big) (\udelta_o + \epsilon \udelta_m) \ufrakd_o r_0.
\end{aligned}
\label{eqn 8.5}
\end{align}
\end{proposition}
\begin{proof}
We choose $\delta$ sufficiently small such that proposition \ref{pro 3.3} holds. Then proposition \ref{pro 6.2} also holds as stated in its proof. The estimate of $\er{\uf}$ follows directly from proposition \ref{pro 6.2} and $\overline{\er{\uf}}^{\circg} = \overline{\dd{\uf}}^{\circg} - \overline{\dd{\ufl{s=0}}}^{\circg}$. 

We integrate equation \eqref{eqn 8.3} to estimate $\ddslashd \er{\uf}$. Since $\ufl{2,s=0}$ is constant, $\dd{\Xlt}^i \ddpartial_i \big( \ddcircDelta \ufl{2,t} \big)$ vanishes and $\dd{\relt} = -\rel{1,t}$. Note $\dslashd \ufl{1,s=0} = - \dslashd \dd{\ufl{s=0}}$, then
\begin{align*}
\leVe \dslashd \ufl{1,s=0}  \riVe^{n,p} \leq \min\{\udelta_o, \ufrakd_o\} r_0.
\end{align*} 
Applying estimate \eqref{eqn 3.8} to $-\rel{1,t}$, we obtain that
\begin{align*}
\leVe \dd{\relt} \riVe^{n-1,p} = \leVe \rel{1,t} \riVe^{n-1,p} \leq c(n,p)(|\os|/r_0 + \delta_o) (\udelta_o + \epsilon \udelta_m ) \ufrakd_o r_0.
\end{align*}
$\leve \overline{\Xl{1,t}^i \ddpartial_i \Big( \ddcircDelta \bdd{\uflt} \Big)}^{\circg} \rive$ satisfies the same estimate as in proposition \ref{pro 8.3}.
Then the estimate of $\ddslashd \er{\uf}$ follows from equation \eqref{eqn 8.3} and the estimates of $\dd{\relt}$ and $\leve \overline{\Xl{1,t}^i \ddpartial_i \big( \ddcircDelta \bdd{\uflt} \big)}^{\circg} \rive$.
\end{proof}
\begin{remark}
Proposition \ref{pro 8.5} requires the $L^{\infty}$ bounds of the metric components up to $(n+1)$-th order derivatives, the same as proposition \ref{pro 3.3}.
\end{remark}

\section{Linearised perturbation of the outgoing null expansion}\label{sec 9}
In section \ref{sec 7}, we study the perturbation of  the outgoing null expansion. In this section, we shall construct a linearised perturbation of the outgoing null expansion.

Recall the following decomposition of the outgoing null expansion in section \ref{sec 4},
\begin{align*}
&
\lo{\ddtr \ddchi'}
=
\tr \chi'_S|_{\Sigma} - 2 \left( r_S|_{\Sigma}\right)^{-2} \circDelta f,
\\
&
\hi{\ddtr \ddchi'} = \ddtr \ddchi' - \lo{\ddtr \ddchi'}.
\end{align*}
When constructing the linearised perturbation of $\ddtr \ddchi'$, we will neglect the high order remainder part $\hi{\ddtr \ddchi'}$, as its perturbation is at least one magnitude smaller than the perturbation of the surface, as shown in proposition \ref{pro 7.5}.

The idea behind the linearised perturbation of $\ddtr \ddchi'$ is that the error should be one magnitude smaller than the perturbation of the surface. We shall verify this by estimating the error.

\subsection{Construction of linearised perturbation of $\ddtr \ddchi'$}
We use $\bdd{\ddtr \ddchi'}$ to denote the linearised perturbation of $\ddtr \ddchi'$. The first order main part of $\ddtr \ddchi'$ is 
\begin{align*}
&
\lo{\ddtr \ddchi'}
=
\tr \chi'_S|_{\Sigma} - 2 ( r_S|_{\Sigma} )^{-2} \circDelta f.
\end{align*}
We use the linearised perturbation $(\bdd{\uf}, \bdd{f})$ in definition \ref{def 8.1} to construct $\bdd{\ddtr \ddchi'}$.

\begin{definition}\label{def 9.1}
Let $\Sigmal{a},a=1,2$ be two spacelike surfaces in $(M,g)$. Suppose $\Sigmal{a}$ has the second parameterisation $( \ufl{a,s=0}, \fl{a})$. We define the linearised perturbation $\bdd{\ddtr \ddchi'}$ of the outgoing null expansion from $\Sigmal{1}$ to $\Sigmal{2}$ as follows:
\begin{align}
\bdd{\ddtr \ddchi'} 
= 
( \partial_{\us} \tr \chi'_S )|_{\Sigmal{1}} \cdot  \bdd{\uf} 
+ 
(\partial_s \tr \chi'_S )|_{\Sigmal{1}} \cdot \bdd{f} 
- 
2 ( r_S|_{\Sigmal{1}})^{-2} \circDelta \big( \bdd{f} \big),
\label{eqn 9.1}
\end{align}
where $\bdd{\uf}, \bdd{f}$ are the linearised perturbations in definition \ref{def 8.1}, and $( \partial_{\us} \tr \chi'_S )|_{\Sigmal{1}}$, $(\partial_s \tr \chi'_S )|_{\Sigmal{1}}$, $r_S|_{\Sigmal{1}}$ are the geometric quantities in the Schwarzschild spacetime on $\Sigmal{1}$.
\end{definition}

In the above definition, we only take the first order main part $\lo{\ddtr \ddchi'}$ into account. We also neglect the perturbation of $r_S|_{\Sigma}$ in the term $( r_S|_{\Sigma} )^{-2} \circDelta f$. The reason to omit the high order remainder part $\hi{\ddtr \ddchi'}$ is briefly explained in the beginning of the section. The same reason applies to the negligibility of $\dd{\left( r_S|_{\Sigma}\right)^{-2}} \circDelta f$.

\subsection{Estimate of error of the linearised perturbation $\bdd{\ddtr \ddchi'}$}
Let $\Sigmal{a}$ be two spacelike surfaces in $(M,g)$. The perturbation of the outgoing null expansions from $\Sigmal{1}$ to $\Sigmal{2}$ 
is
\begin{align*}
\dd{\ddtr \ddchi'} = \uls{2}{\big( \ddtr \ddchi' \big)} - \uls{1}{\big( \ddtr \ddchi' \big)}.
\end{align*}
Denote the error of the linearised perturbation $\bdd{\ddtr \ddchi'}$ from $\Sigmal{1}$ to $\Sigmal{2}$ by $\er{\ddtr \ddchi'}$,
\begin{align*}
\er{\ddtr \ddchi'} = \dd{\ddtr \ddchi'} - \bdd{\ddtr \ddchi'}.
\end{align*}
We have the following decomposition of $\er{\ddtr \ddchi'}$ from the decomposition of $\ddtr \ddchi'$
\begin{align}
\er{\ddtr \ddchi'} = 
\Big( \dd{\lo{\ddtr \ddchi'}} - \bdd{\ddtr \ddchi'} \Big)
+ 
\dd{\hi{\ddtr \ddchi'}}.
\label{eqn 9.2}
\end{align}
Summarise the estimates of the parameterisation functions as follows: with the notations of $\ud_o, \ud_m,  \uslashd_m,\ubfd_o, \uslashbfd_m, \ubfd_m, \ubfd_{\uh}, \ud_{\uh}$ in formulae \eqref{eqn 3.7} \eqref{eqn 5.8} \eqref{eqn 6.4} \eqref{appen pro 7.1 eqn 1} \eqref{eqn 7.4},
\begin{align}
\begin{aligned}
&
\big\Vert \ddslashd \ufl{a} \big\Vert^{n,p} \leq \ud_o r_0,
\quad
\big| \overline{\ufl{a}}^{\circg} - \overline{\ufl{a,s=0}}^{\circg} \big| \leq \uslashd_m r_0,
\\
&
\big\Vert \ddslashd \dd{\uf} \big\Vert^{n-1,p} \leq \ubfd_o r_0,
\quad
\big| \overline{\dd{\uf}}^{\circg} - \overline{\dd{\ufl{s=0}}}^{\circg} \big| \leq \uslashbfd_m r_0,
\\
&
\big\Vert \uls{a}{( \dslashd \uh )|_{\Sigmal{a}}} \big\Vert^{n+1,p},
\big\Vert \uls{a}{( \circnabla^2 \uh )|_{\Sigmal{a}}} \big\Vert^{n,p}
\leq
\ud_{\uh} r_0,
\\
&
\big\Vert \dd{( \dslashd \uh )|_{\Sigma}} \big\Vert^{n,p},
\big\Vert \dd{( \circnabla^2 \uh )|_{\Sigma}} \big\Vert^{n-1,p}
\leq
\ubfd_{\uh} r_0.
\end{aligned}
\tag{\ref{eqn 7.5}}
\end{align}
Furthermore introduce the notation $\ubfe$ by
\begin{align}
\begin{aligned}
\ubfe 
=&
c(n,p) \big(|\os|/r_0+\delta_o\big) 
\left[ (\udelta_o^2 + \epsilon \udelta_o) \ufrakd_m
+
(\udelta_o + \epsilon \udelta_m) \ufrakd_o \right]
\\
&
+
c(n,p) (\udelta_o^2 + \epsilon \udelta_o \udelta_m) (\frakd_m + \frakd_o),
\end{aligned}
\label{eqn 9.3}
\end{align}
then estimates \eqref{eqn 8.4} for $\er{\uf}$ can be rewritten as
\begin{align}
\leve \overline{\er{\uf}}^{\circg} \rive,
\leVe \ddslashd \er{\uf} \riVe^{n-1,p},
\leVe \er{\uf} \riVe^{n,p}
\leq
\ubfe r_0.
\label{eqn 9.4}
\end{align}
We shall estimate $\er{\ddtr \ddchi'}$ by above estimates \eqref{eqn 7.5} \eqref{eqn 9.4}.
\begin{proposition}\label{pro 9.2}
Let $\Sigmal{a}, a=1,2$ be two spacelike surfaces in $(M,g)$. Assume that $\Sigmal{a}$ has the second parameterisation $(\ufl{a,s=0}, \fl{a})$. Suppose that the parameterisation functions satisfy the following estimates
\begin{align*}
&
\leVe \dslashd \ufl{a,s=0} \riVe^{n+1,p} 
\leq
\udelta_o r_0,
\quad
\leve \overline{\ufl{a,s=0}}^{\circg} \rive
\leq 
\udelta_m r_0,
\\
&
\leVe \ddslashd \fl{a} \riVe^{n+1,p} \leq  \delta_o r_0,
\quad
\overline{\fl{a}}^{\circg} = \os_a,
\quad
\os = \max \{ |\os_1|, |\os_2| \},
\end{align*}
and the perturbation functions satisfy
\begin{align*}
&
\leVe \dslashd \dd{\ufl{s=0}} \riVe^{n,p}  \leq  \ufrakd_o r_0,
\quad
\leve \overline{\dd{\ufl{s=0}}}^{\circg} \rive \leq \ufrakd_m r_0,
\\
&
\leVe \ddslashd \dd{f} \riVe^{n,p} \leq \frakd_o r_0,
\quad
\leve \overline{\dd{f}}^{\circg} \rive \leq \frakd_m r_0,
\end{align*}
where $n\geq 2, p>2$ or $n\geq 3, p>1$.

There exist a small positive constant $\delta$ depending on $n,p$ and constants $c(n,p)$, such that if $\epsilon, \udelta_o, \udelta_m, \delta_o, \ufrakd_o, \ufrakd_m, \frakd_o, \frakd_m$ are suitably bounded that
\begin{align*}
\epsilon,\ \udelta_o,\ \epsilon \udelta_m,\ \delta_o,\ \ufrakd_o,\ \epsilon \ufrakd_m,\ \frakd_o,\ \frakd_m \leq \delta,
\end{align*}
then the error $\er{\ddtr \ddchi'}$ satisfies the following estimate,
\begin{align}
\begin{aligned}
\leVe \er{\ddtr \ddchi'} \riVe^{n-1,p}
\leq
&
\frac{c(n,p)(\epsilon +\delta_o + \udelta_o^2 |\os|/r_0) }{r_0} \frakd_m
+
\frac{c(n,p)(\epsilon +\delta_o + \udelta_o |\os|/r_0 )}{r_0} \frakd_o
\\
&
+
\frac{c(n,p)\big(\epsilon + \delta_o^2 + \delta_o |\os|/r_0 + \udelta_o^2 (|\os|/r_0)^2 \big)}{r_0} \ufrakd_m
\\
&
+
\frac{c(n,p)\big(\epsilon +\delta_o^2 + \delta_o |\os|/r_0 + \udelta_o (|\os|/r_0)^2 \big)}{r_0} \ufrakd_o
\\
&
+
\frac{c(n,p)(|\os|/r_0 + \delta_o)}{r_0} (\ufrakd_m + \ufrakd_o)^2
\\
&
+
\frac{c(n,p)}{r_0} (\ufrakd_m + \ufrakd_o + \frakd_m + \frakd_o)(\frakd_m + \frakd_o).
\end{aligned}
\label{eqn 9.5}
\end{align}
\end{proposition}
\begin{proof}
It is sufficient to estimate $\dd{\lo{\ddtr \ddchi'}} - \bdd{\ddtr \ddchi'}$ and $\dd{\hi{\ddtr \ddchi'}}$ seperately by formula \eqref{eqn 9.2}. Assume that $\delta$ is sufficiently small such that propositions \ref{pro 3.3}, \ref{pro 5.1}, \ref{pro 5.3}, \ref{pro 7.1}, \ref{pro 8.3} are true, i.e. estimates \eqref{eqn 7.5} and \eqref{eqn 9.4} hold.

The estimate of $\dd{\hi{\ddtr \ddchi'}}$ is already obtain by proposition \ref{pro 7.3}. We estimate $\dd{\lo{\ddtr \ddchi'}} - \bdd{\ddtr \ddchi'}$ in the following.
\begin{align*}
\dd{\lo{\ddtr \ddchi'}} - \bdd{\ddtr \ddchi'}
=
&
\underbrace{\dd{\tr \chi'_S} 
- \big( \partial_{\us} \tr \chi'_S \big)|_{\Sigmal{1}} \cdot  \bdd{\uf}
- \big(\partial_s \tr \chi'_S\big)|_{\Sigmal{1}} \cdot \bdd{f}}_{\mathbf{I}}
\\
&
\underbrace{-2 \dd{(r_S|_{\Sigma})^{-2}}  \circDelta \big( \fl{2} \big)}_{\mathbf{II}}.
\end{align*}
For $\mathbf{I}$, we have
\begin{align*}
\dd{\tr \chi'_S} 
=
\int_0^1 \partial_{\us} \tr \chi'_S \big( \ufl{1} + t \dd{\uf}, \fl{1} \big) \d t  \cdot \dd{\uf}
+
\int_0^1 \partial_s \tr \chi'_S \big( \ufl{2}, \fl{1} + t \dd{f} \big) \d t \cdot \dd{f}.
\end{align*}
Note that $\dd{f} = \bdd{f}$, then
\begin{align*}
\mathbf{I}
=&
\underbrace{\int_0^1 
\left[ \partial_{\us} \tr \chi'_S \big( \ufl{1} + t \dd{\uf}, \fl{1} \big) - \partial_{\us} \tr \chi'_S \big( \ufl{1}, \fl{1} \big) \right]
\d t 
\cdot \dd{\uf}}_{\mathbf{I}_a}
+
\underbrace{\big( \partial_{\us} \tr \chi'_S \big)|_{\Sigmal{1}} \cdot \er{\uf}}_{\mathbf{I}_b}
\\
&
+
\underbrace{\int_0^1 
\left[ \partial_s \tr \chi'_S \big( \ufl{2}, \fl{1} + t \dd{f} \big)  - \partial_s \tr \chi'_S \big( \ufl{1}, \fl{1} \big) \right] 
\d t \cdot \dd{f}}_{\mathbf{I}_c}.
\end{align*}
We estimate $\mathbf{I}_a$, $\mathbf{I}_b$, $\mathbf{I}_c$ in the following.
\begin{align*}
&\begin{aligned}
\Vert \mathbf{I}_a \Vert^{n-1,p}
\leq&
\sup_{t\in [0,1]} \leVe \partial_{\us}^2 \tr \chi'_S \big( \ufl{1} + t \dd{\uf}, \fl{1} \big)  \riVe^{n-1,p} \big( \leVe \dd{\uf} \riVe^{n-1,p} \big)^2
\\
\leq&
\frac{c(n,p) (|\os|/r_0 + \delta_o)}{r_0} (\ubfd_m + \ubfd_o)^2,
\end{aligned}
\\
&
\Vert \mathbf{I}_b \Vert^{n-1,p} 
\leq
\leVe \partial_{\us} \tr \chi'_S \riVe^{n-1,p} \leVe \er{\uf} \riVe^{n-1,p}
\leq
\frac{c(n,p)(|\os|/r_0 + \delta_o)}{r_0} \ubfe,
\\
&\begin{aligned}
\Vert \mathbf{I}_c \Vert^{n-1,p}
\leq&
\sup_{t\in [0,1]} \leVe \partial_{\us} \partial_s \tr \chi'_S \big( \ufl{1} + t \dd{\uf}, \fl{1} \big)  \riVe^{n-1,p} \cdot \leVe \dd{\uf} \riVe^{n-1,p} \cdot \leVe \dd{f} \riVe^{n-1,p}
\\
&
+
\sup_{t\in [0,1]} \leVe \partial_s^2 \tr \chi'_S \big( \ufl{2}, \fl{1} + t \dd{f} \big)  \riVe^{n-1,p}  
\cdot \big( \leVe \dd{f} \riVe^{n-1,p} \big)^2
\\
\leq
&
\frac{c(n,p)}{r_0} (\ubfd_m + \ubfd_o + \frakd_m + \frakd_o)(\frakd_m + \frakd_o).
\end{aligned}
\end{align*}
Therefore
\begin{align*}
\Vert \mathbf{I} \Vert^{n-1,p}
\leq
&
\frac{c(n,p)(|\os|/r_0 + \delta_o)}{r_0} 
\left[ \ubfe + (\ubfd_m + \ubfd_o)^2 \right]
+
\frac{c(n,p)}{r_0} (\ubfd_m + \ubfd_o + \frakd_m + \frakd_o)(\frakd_m + \frakd_o).
\end{align*}
For $\mathbf{II}$, we have
\begin{align*}
&
\begin{aligned}
\dd{(r_S|_{\Sigma} )^{-2}}
= 
\int_0^1 
-2 \left[
\frac{\partial_{\us} r_S}{r_S^3} 
\big( \ufl{1} + t \dd{\uf}, \fl{1} \big) \cdot \dd{\uf} 
+
\frac{\partial_s r_S}{r_S^3} \big( \ufl{2}, \fl{1} + t \dd{f} \big) \cdot \dd{f} \right] \d t,
\end{aligned}
\\
&
\leVe \dd{(r_S|_{\Sigma} )^{-2}} \riVe^{n-1,p}
\leq
\frac{c(n,p) (|\os|/r_0 + \delta_o) }{r_0^2} (\ubfd_m + \ubfd_o) + \frac{c(n,p)}{r_0^2} ( \frakd_m + \frakd_o),
\end{align*}
therefore
\begin{align*}
\Vert \mathbf{II} \Vert^{n-1,p} 
\leq
\frac{c(n,p) (|\os|/r_0 + \delta_o) \delta_o }{r_0} (\ubfd_m + \ubfd_o) + \frac{c(n,p) \delta_o}{r_0} ( \frakd_m + \frakd_o).
\end{align*}
Assembling the above estimates, we obtain that
\begin{align*}
\leVe \dd{\ddtr \ddchi'} \riVe^{n-1,p}
\leq
\Vert \mathbf{I} \Vert^{n-1,p} + \Vert \mathbf{II} \Vert^{n-1,p} + \leVe \dd{\hi{\ddtr \ddchi'}} \riVe^{n-1,p}.
\end{align*}
Substituting $\ud_o$, $\ud_m$, $\uslashd_m$, $\ubfd_o$, $\uslashbfd_m$, $\ubfd_m$, $\ubfd_{\uh}$, $\ud_{\uh}$, $\ubfe$ in formulae \eqref{eqn 3.7} \eqref{eqn 5.8} \eqref{eqn 6.4} \eqref{appen pro 7.1 eqn 1} \eqref{eqn 7.4} \eqref{eqn 9.3}, we prove the estimate of $\dd{\ddtr \ddchi'}$ in the proposition.
\end{proof}

\begin{remark}
Proposition \ref{pro 9.2} requires the $L^{\infty}$ bounds of the metric components up to $(n+2)$-th order derivatives, same as proposition \ref{pro 7.3} for the estimate of $\dd{\hi{\ddtr \ddchi'}}$.

It also requires the $L^{\infty}$ bounds of the structure coefficients up to $n$th order derivative, again same as proposition \ref{pro 7.3}.

Note that we donot need the estimate of $\er{\uf}$ with the highest regularity obtained in proposition \ref{pro 8.3} in the proof. The estimate of $\Vert \er{\uf} \Vert^{n-1,p}$ is sufficient for the above proposition, while we obtain the better estimate of $\Vert \er{\uf} \Vert^{n,p}$ in proposition \ref{pro 8.3}.
\end{remark}

\subsection{Improved estimate of error of the linearised perturbation $\bdd{\ddtr \ddchi'}$}
We can improve the estimate of $\er{\ddtr \ddchi'}$ if assuming additional conditions as in proposition \ref{pro 7.5} on the improved estimate of $\dd{\hi{\ddtr \ddchi'}}$.

Summarise the improved estimates of the parameterisation functions as follows: with the notations $\ubfd'_o$, $\uslashbfd'_m$, $\ubfd'_m$, $\ubfd'_h$ in equations \eqref{eqn 7.6},
\begin{align*}
\begin{aligned}
&
\big\Vert \ddslashd \dd{\uf} \big\Vert^{n,p} \leq \ubfd'_o r_0,
\quad
\big| \overline{\dd{\uf}}^{\circg} - \overline{\dd{\ufl{s=0}}}^{\circg} \big| \leq \uslashbfd'_m r_0,
\\
&
\big\Vert \dd{\big( \dslashd \uh \big)|_{\Sigma}} \big\Vert^{n+1,p},
\big\Vert \dd{\big( \circnabla^2 \uh \big)|_{\Sigma}} \big\Vert^{n,p}
\leq
\ubfd'_{\uh} r_0.
\end{aligned}
\tag{\ref{eqn 7.7}}
\end{align*}
Introduce the notation $\ubfe'$ by
\begin{align*}
\ubfe' = c(n,p) (|\os|/r_0 + \delta_o)(\udelta_o + \epsilon \udelta_m) \ufrakd_o r_0,
\end{align*}
then the improved estimate \eqref{eqn 8.5} for $\er{\uf}$ can be rewritten as
\begin{align}
\leve \overline{\er{\uf}}^{\circg} \rive,
\leVe \ddslashd \er{\uf} \riVe^{n,p},
\leVe \er{\uf} \riVe^{n+1,p}
\leq
&
\ubfe' r_0.
\label{eqn 9.6}
\end{align}
The improved estimate of $\er{\ddtr \ddchi'}$ follows from the above improved estimates \eqref{eqn 7.7} \eqref{eqn 9.6}.
\begin{proposition}\label{pro 9.4}
Under the setting of proposition \ref{pro 9.2}, we assume additionally that $\ufl{2,s=0}$ or $\ufl{1,s=0}$ is constant and 
\begin{align*}
\leVe \dslashd \dd{\ufl{s=0}} \riVe^{n+1,p} \leq \ufrakd_o r_0,
\quad
\leVe \ddslashd \dd{f} \riVe^{n+1,p} \leq \frakd_o r_0.
\end{align*}

There exists a small positive constant $\delta$ depending on $n,p$ and constants $c(n,p)$, such that if $\epsilon, \udelta_o, \udelta_m, \delta_o, \ufrakd_o, \ufrakd_m, \frakd_o, \frakd_m$ are suitably bounded that
\begin{align*}
\epsilon,\ \udelta_o,\ \epsilon \udelta_m,\ \delta_o,\ \ufrakd_o,\ \epsilon \ufrakd_m,\ \frakd_o,\ \frakd_m \leq \delta,
\end{align*}
\begin{enumerate}[label=\alph*.]
\item \label{pro 9.4.a}
then $\er{\ddtr \ddchi'}$ satisfies the estimate of the same form as \eqref{eqn 9.5} with the Sobolev norm improved to $\leVe \er{\ddtr \ddchi'} \riVe^{n,p}$;

\item \label{pro 9.4.b}
and if we assume that it is $\ufl{2,s=0}$ being constant, then $\er{\ddtr \ddchi'}$ satisfies the following improved estimate
\begin{align*}
\leVe \er{\ddtr \ddchi'} \riVe^{n,p}
\leq
&
\frac{c(n,p) ( \epsilon + \delta_o^2 + \delta_o \udelta_o)}{r_0} \frakd_m
+
\frac{c(n,p) ( \epsilon + \delta_o + \udelta_o |\os|/r_0)}{r_0} \frakd_o
\\
&
+
\frac{c(n,p) (\epsilon + \delta_o^2 + \delta_o \udelta_o |\os|/r_0)}{r_0} \ufrakd_m
\\
&
+
\frac{c(n,p) \big(\epsilon + \delta_o^2 + \delta_o |\os|/r_0 + \udelta_o (|\os|/r_0)^2 \big)}{r_0} \ufrakd_o
\\
&
+
\frac{c(n,p)(|\os|/r_0 + \delta_o)}{r_0} (\ufrakd_m + \ufrakd_o)^2
\\
&
+
\frac{c(n,p)}{r_0} (\ufrakd_m + \ufrakd_o + \frakd_m + \frakd_o)(\frakd_m + \frakd_o).
\end{align*}
\end{enumerate}
\end{proposition}
\begin{proof}[Proof of proposition \ref{pro 9.4}.\ref{pro 9.4.a}]
The proof follows exactly the proof of proposition \ref{pro 9.2} by improving the Sobolev norm from $\Vert \cdot \Vert^{n-1,p}$ to $\Vert \cdot \Vert^{n,p}$ and replacing $\ubfd_m, \ubfd_o, \ubfd_{\uh}$ by $\ubfd'_m, \ubfd'_o, \ubfd'_{\uh}$. Note that estimate \eqref{eqn 9.4} of $\er{\uf}$ is sufficient for the improved Sobolev norm. The choice of $\delta$ requires that propositions \ref{pro 3.3}, \ref{pro 5.1}, \ref{pro 5.3}, \ref{pro 7.2}, \ref{pro 8.3} are true.
\end{proof}

\begin{proof}[Proof of proposition \ref{pro 9.4}.\ref{pro 9.4.b}]
The proof follows the same scheme as the proof of proposition \ref{pro 9.2}. In additional to the above proof of \ref{pro 9.4}.\ref{pro 9.4.a}, we replace $\ubfe$ in the proof of proposition \ref{pro 9.2} by $\ubfe'$ and note that the term $\mathbf{II}$ vanishes, since $\circDelta\ufl{2}=0$. The choice of $\delta$ requires that the improved estimate \eqref{eqn 9.6} of $\er{\uf}$ holds, thus it additionally requires that the proposition \ref{pro 8.5} is true.
\end{proof}

\begin{remark}
Proposition \ref{pro 9.4} requires the $L^{\infty}$ bounds of the metric components up to $(n+2)$-th order derivatives, the same as proposition \ref{pro 9.2}. 

While it requires the $L^{\infty}$ bounds of the structure coefficients up to $(n+1)$-th order derivatives, one order higher than the requirement of proposition \ref{pro 9.2}. It is the requirement of proposition \ref{pro 7.3} on the improved estimate of $\dd{\hi{\ddtr \ddchi'}}$.
\end{remark}

\subsection{Invertibility of the linearised perturbation $\bdd{\ddtr \ddchi'}$}
We view $\bdd{\ddtr \ddchi'}$ in definition \ref{def 9.1} as a linear operator defined by equation \eqref{eqn 9.1}
\begin{align*}
\begin{aligned}
\bdd{\ddtr \ddchi'}|_{\Sigmal{1}} ( \bdd{\uf}, \bdd{f} )
\footnotemark
= &
( \partial_{\us} \tr \chi'_S )|_{\Sigmal{1}} \cdot  \bdd{\uf} 
\\
&
+ 
(\partial_s \tr \chi'_S )|_{\Sigmal{1}} \cdot \bdd{f} 
- 
2 ( r_S|_{\Sigmal{1}} )^{-2} \circDelta ( \bdd{f} ).
\end{aligned}
\tag{\ref{eqn 9.1}}
\end{align*}
\footnotetext{We emphasis where the linearised perturbation $\bdd{\ddtr \ddchi'}$ is defined by the subscript.}
$\tr \chi'_S$ is calculated from the values of $\Omega_S^2$ in formulae \eqref{eqn 2.1} and $\tr \chi_S$ in \eqref{eqn 2.2}
\begin{align*}
\tr \chi'_S 
= \Omega_S^{-2} \tr \chi_S 
=  \frac{2s}{r(r_0+s)},
\end{align*}
and we calculate $\partial_{\us} \tr \chi'_S$, $\partial_s \tr \chi'_S$ directly
\begin{align*}
&
\partial_{\us} \tr \chi'_S 
= 
-\frac{2s}{r^2(r_0+s)} \partial_{\us} r
=
-\frac{2s(r-r_0)}{r^3(r_0+s)},
\\
&
\partial_s \tr \chi'_S 
= 
-\frac{2s}{r^2(r_0+s)} \partial_{s} r 
+ \frac{2r_0}{r(r_0+s)^2}
=
-\frac{2(r-r_0)}{r^3}
+ \frac{2r_0}{r(r_0+s)^2}.
\end{align*}
From the above calculations, the linearised perturbation $\bdd{\ddtr \ddchi'}$ at the surface $\Sigma_{\us,s=0}$ is
\begin{align*}
\begin{aligned}
\bdd{\ddtr \ddchi'}|_{\Sigma_{\us,s=0}} ( \bdd{\uf}=0, \bdd{f} )
= &
2 r_0^{-2} \big[ \bdd{f} - \circDelta ( \bdd{f} ) \big].
\end{aligned}
\end{align*}
On a general surface $\Sigma_{\us,s}$ in $M$, we have
\begin{align*}
\bdd{\ddtr \ddchi'}|_{\Sigma_{\us,s}} ( \bdd{\uf}=0, \bdd{f} )
= &
2 r^{-2} \Big[ \frac{r_0}{r} + \frac{r_0 r}{(r_0+s)^2} + (-1 - \circDelta) \Big] ( \bdd{f} ).
\end{align*}
By this explicit formula, we obtain that the linear operator $\bdd{\ddtr \ddchi'}|_{\Sigma_{\us,s}} ( \bdd{\uf}=0, \cdot )$ is an isomorphism from $\mathrm{W}^{n+2,p}$ to $\mathrm{W}^{n,p}$.\footnote{One can derive that $\frac{r_0 r}{(r_0+s)^2} >\frac{2}{3}, \frac{r_0}{r}-1< \frac{1}{3}$ in $M$ by elementary estimates.}

We can estimate the difference between $\bdd{\ddtr \ddchi'}|_{\Sigmal{1}}$ and $\bdd{\ddtr \ddchi'}|_{\Sigma_{\us,s}}$.
\begin{proposition}\label{pro 9.6}
Let $\Sigma$ be a spacelike surface in $(M,g)$. Assume that $\Sigma$ has the second parameterisation $(\ufl{s=0}, f)$. Suppose that the parameterisation functions satisfy the following estimates
\begin{align*}
&
\leVe \dslashd \ufl{s=0} \riVe^{n,p} 
\leq
\udelta_o r_0,
\quad
\leve \overline{\ufl{s=0}}^{\circg} \rive
\leq 
\udelta_m r_0,
\\
&
\leVe \ddslashd f \riVe^{n+1,p} \leq  \delta_o r_0,
\quad
\overline{f}^{\circg} = \os,
\end{align*}
where $n\geq 2, p$ or $n\geq 3, p>1$.

There exist a small positive constant $\delta$ depending on $n,p$ and constants $c(n,p)$, such that if $\epsilon, \udelta_o, \udelta_m, \delta_o$ are suitably bounded that
$\epsilon, \udelta_o, \epsilon \udelta_m, \delta_o \leq \delta$,
then the difference between $\bdd{\ddtr \ddchi'}|_{\Sigma}$ and $\bdd{\ddtr \ddchi'}|{\raisebox{0.9ex}{$_{\Sigma(\us=\overline{\uf}^{\circg},s=\os)}$}}$ satisfies the following estimate: for $0\leq m \leq n$,
\begin{align}
\begin{aligned}
&
\leVe 
\big[ 
\bdd{\ddtr \ddchi'}|_{\Sigma} 
- 
\bdd{\ddtr \ddchi'}|{\raisebox{0.9ex}{$_{\Sigma(\us=\overline{\uf}^{\circg},s=\os)}$}}
\big]
\big( \bdd{\uf}, \bdd{f} \big)
\riVe^{m,p}
\\
&
\leq
\frac{c(n,p) [( |\os|/r_0 + \delta_o )\udelta_o + \delta_o]}{r_0^2} 
\big[ ( |\os|/r_0 + \delta_o ) \leVe \bdd{\uf} \riVe^{m,p}
+
\leVe \bdd{f} \riVe^{m+2,p} \big].
\end{aligned}
\label{eqn 9.7}
\end{align}
\end{proposition}
\begin{proof}
Choose $\delta$ sufficiently small such that proposition \ref{pro 3.3} on the estimate of $\uf$ hold. Then estimate \eqref{eqn 9.7} follows from
\begin{align*}
&
\partial_s \partial_{\us} \tr \chi'_S
=
-\frac{2(r-r_0)}{r^5(r_0+s)^2}
\left[ r_0 r^2 + (3r_0-2r)(r_0+s)^3 \right],
\\
&
\partial_s^2 \tr \chi'_S
=
-\frac{4r_0}{r(r_0+s)^3}
- \frac{2s}{(r-r_0)} \cdot \frac{1}{r^5(r_0+s)} \left[ r_0 r^2 + (3r_0-2r)(r_0+s)^3 \right].
\end{align*}
and the differential of $r$ in formulae \eqref{eqn 2.2}.
\end{proof}
The above proposition implies that $\bdd{\ddtr \ddchi'}|_{\Sigma}$ is a small perturbation of $\bdd{\ddtr \ddchi'}|_{\Sigma_{\us,s=0}}$. As a corollary, we show that the linear operator $\bdd{\ddtr \ddchi'}|_{\Sigma} \big( \bdd{\uf}=0, \cdot \big)$ is an isomorphism from $\mathrm{W}^{n+2,p}$ to $\mathrm{W}^{n,p}$.
\begin{corollary}\label{coro 9.7}
Under the setting of proposition \ref{pro 9.6}, there exist a small positive constant $\delta$ depending on $n,p$ and a constant $c(n,p)$, such that if $\epsilon, \udelta_o, \udelta_m, \delta_o$ are suitably bounded that
$\epsilon,\, \udelta_o,\, \epsilon \udelta_m,\, \delta_o \leq \delta$,
then $\bdd{\ddtr \ddchi'}|_{\Sigma} \big( \bdd{\uf}=0, \cdot \big)$ is an isomorphism from $W^{m+2,p}$ to $W^{m,p}$ for all $0\leq m\leq n$, with the following estimate
\begin{align*}
\frac{c(n,p)^{-1}}{r_0^2} \leVe \bdd{f} \riVe^{m+2,p}
\leq
\leVe \bdd{\ddtr \ddchi'}|_{\Sigma} \big( 0, \bdd{f} \big) \riVe^{m,p}
\leq
\frac{c(n,p)}{r_0^2} \leVe \bdd{f} \riVe^{m+2,p}.
\end{align*}
\end{corollary}
The proof is straightforward, simply by applying estimate \eqref{eqn 9.7}.

\section{Construction of marginally trapped surfaces}\label{sec 10}
In this section, we shall construct marginally trapped surfaces in the perturbed Schwarzschild spacetime in definition \ref{def 2.3}.

\subsection{Translate the construction to a problem of analysis}
We introduced methods to parametrise spacelike surfaces in section \ref{sec 3}, and obtained the formula of the outgoing null expansion in section \ref{sec 4}. Then we define a map, denoted by $\bft$, from the parameterisation to the outgoing null expansion of the corresponding spacelike surface
\begin{align*}
\bft:
\quad
(\ufl{s=0}, f)
\quad
\overset{\text{a.}}{\longrightarrow}
\quad
\Sigma
\quad
\overset{\text{b.}}{\longrightarrow}
\quad
\ddtr \ddchi',
\end{align*}
\begin{enumerate}[label=\alph*.]
\item $\Sigma$ is the surface with the second parameterisation $(\ufl{s=0}, f)$,
\item $\ddtr \ddchi'$ is the outgoing null expansion of $\Sigma$.
\end{enumerate}
Therefore, in order to construct a marginally trapped surface, it is sufficient to solve the equation
\begin{align}
\ddtr \ddchi'(\Sigma) = \bft (\ufl{s=0}, f)= 0,
\label{eqn 10.1}
\end{align}
and show that the incoming null expansion $\ddtr \dduchi$ of $\Sigma$ is non-positive. Sometimes we also use $\ddtr \ddchi'(\ufl{s=0}, f)$ to denote $\bft(\ufl{s=0},f)$ to emphasis the concrete geometric meaning of $\bft$.

We proved that the linearised perturbation $\bdd{\ddtr \ddchi'}$ is invertible for the second slot $\bdd{f}$ in  corollary \ref{coro 9.7}. Thus intuitively by the implicit function theorem, the equation $\bft (\ufl{s=0}, f)= 0$ should have a unique solution $f$ for given $\ufl{s=0}$.

We shall use the linearised perturbation $\bdd{\ddtr \ddchi'}$ to construct approximating solutions to solve equation \eqref{eqn 10.1}. For preparations, we introduce some notations: suppose that $\Sigma$ has the second parameterisation $(\ufl{s=0}, f)$,
\begin{enumerate}[label=\alph*.]
\item $\delta \bft [\ufl{s=0}, f] $: we use this notation to denote $\bdd{\ddtr \ddchi'}|_{\Sigma}$, 
\begin{align*}
\delta \bft [\ufl{s=0}, f] ( \bdd{\ufl{s=0}}, \bdd{f} )
= 
\bdd{\ddtr \ddchi'}|_{\Sigma}  ( \bdd{\ufl{s=0}}, \bdd{f} ).
\end{align*}

\item $\upartial \bft [\ufl{s=0}, f]$: we use this notation to denote $\bdd{\ddtr \ddchi'}|_{\Sigma}$ in the first slot $\bdd{\ufl{s=0}}$,
\begin{align*}
\upartial \bft [\ufl{s=0}, f] ( \bdd{\ufl{s=0}} ) 
=
\delta \bft [\ufl{s=0}, f] ( \bdd{\ufl{s=0}}, 0 )
= 
\bdd{\ddtr \ddchi'}|_{\Sigma}  ( \bdd{\ufl{s=0}}, 0 ).
\end{align*}

\item $\partial \bft [\ufl{s=0}, f]$: we use this notation to denote $\bdd{\ddtr \ddchi'}|_{\Sigma}$ in the second slot $\bdd{f}$,
\begin{align*}
\partial \bft [\ufl{s=0}, f] ( \bdd{f} )
=
\delta \bft [\ufl{s=0}, f] ( 0 , \bdd{f} )
=
\bdd{\ddtr \ddchi'}|_{\Sigma} ( 0, \bdd{f} ).
\end{align*}
\end{enumerate}

The sequence of approximating solutions $\{\fl{k}\}_{k\in \mathbb{N}}$ of equation \eqref{eqn 10.1} is constructed as follows: given the parameterisation function $\ufl{s=0}$
\begin{enumerate}[label=\alph*.]
\item
set $\fl{k=0} =0$;
\item
solve the linear equation successively:
\begin{align}
\partial \bft [ \ufl{s=0}, \fl{k}] \big( \fl{k+1} - \fl{k} \big)
= - \bft\big( \ufl{s=0}, \fl{k} \big). 
\label{eqn 10.2}
\end{align}
\end{enumerate}
The goal is to prove that $\{\fl{k}\}_{k \in \mathbb{N}}$ exists and converges, and the limit solves $\bft(\ufl{s=0}, f)=0$.

The convergence of $\{\fl{k}\}_{k\in \mathbb{N}}$ shall be proved in two steps: suppose that $\ufl{s=0}$ belongs to the Sobolev space $\mathrm{W}^{n+2,p}$, then
\begin{enumerate}[label=\roman*.]
\item we first prove that $\{ \fl{k} \}_{k\in \mathbb{N}}$ exists and is bounded in $\mathrm{W}^{n+2,p}$,
\item then prove that the sequence converges in the weaker Sobolev space $\mathrm{W}^{n+1,p}$.
\end{enumerate}
The latter will imply that the limit solves the equation $\bft(\ufl{s=0}, f)=0$.

\subsection{Existence and boundedness of the sequence of approximating solutions}
We need to solve the linear equation \eqref{eqn 10.2} in the construction of $\{\fl{k}\}_{k\in \mathbb{N}}$. It is solvable if the linear operator $\partial \bft[ \ufl{s=0}, \fl{k}]$ is invertible. A sufficient condition is given in corollary \ref{coro 9.7}, which requires suitable bounds on $\fl{k}$ in $\mathrm{W}^{n+2,p}$. Thus we shall prove that the existence of $\{\fl{k}\}_{k\in \mathbb{N}}$, and also show its boundedness at the same time.

Introduce the following notations of function spaces:
\begin{align}
&
\dot{\mathrm{W}}^{m,p}(\udelta_o r_0) = \big\{ \uf:\ \Vert \d \uf \Vert^{m-1,p} \leq \udelta_o r_0  \big\},
\quad
\mathrm{M}(\udelta_m r_0) = \big\{ \uf: \vert \overline{\uf}^{\circg} \vert \leq \udelta_m r_0  \big\}.
\label{eqn 10.3}
\end{align}

\begin{lemma}\label{lem 10.1}
Let $\ufl{s=0}$ be a function in $\dot{\mathrm{W}}^{m,p}(\udelta_o r_0) \cap \mathrm{M}(\udelta_m r_0)$, where $n \geq 2, p >2$ or $n\geq 3, p>1$. There exist a small positive constant $\delta$ and a constant $c$, both depending on $n,p$ and satisfying $c \delta <\kappa$, such that if $\epsilon, \udelta_o, \udelta_m$ are suitably bounded that
$\epsilon, \udelta_o, \epsilon \udelta_m \leq \delta$,
then the sequence $\{\fl{k}\}_{k \in \mathbb{N}}$ exists and satisfies the following estimate
\begin{align}
\Vert \fl{k} \Vert^{n+2,p}
\leq
c \epsilon r_0
<
\kappa r_0.
\label{eqn 10.4}
\end{align}
\end{lemma}
\begin{proof}
We shall prove the lemma by induction. Choose $\delta$ sufficiently small such that proposition \ref{pro 9.4}, proposition \ref{pro 9.6}, corollary \ref{coro 9.7} all apply.

Since $\fl{k=0}=0$, it clearly satisfies the estimate in the lemma. Now assume that $\fl{k}$ exists and satisfies estimate \eqref{eqn 10.4}. We prove that $\fl{k+1}$ exists. It is sufficient to show that $\partial\bft [\ufl{s=0}, \fl{k}]$ is invertible. This follows directly from corollary \ref{coro 9.7}.

We show that $\fl{k+1}$ satisfies estimate \eqref{eqn 10.4}. Transform equation \eqref{eqn 10.2} as follows:
\begin{align*}
\partial\bft [\ufl{s=0}, \fl{k}] \big(\fl{k+1}\big) 
=
&
\underbrace{\Big[\bft \big(\overline{\ufl{s=0}}^{\circg}, 0 \big) - \bft \big(\ufl{s=0}, \fl{k} \big) - \delta\bft [\ufl{s=0}, \fl{k}] \big(\overline{\ufl{s=0}}^{\circg} - \ufl{s=0}, - \fl{k} \big)\Big]}_{i}
\\
&
-
\underbrace{\bft \big( \overline{\ufl{s=0}}^{\circg}, 0 \big)}_{ii}.
\end{align*}
By corollary \ref{coro 9.7}, $\fl{k+1}$ is bounded by
\begin{align*}
\Vert \fl{k+1} \Vert^{n+2,p}
\leq
c(n,p) r_0^2 \cdot \Vert  i \Vert^{n,p} + c(n,p) r_0^2 \cdot \Vert ii \Vert^{n,p}.
\end{align*}
Note that the term $i$ is actually the error of the linearised perturbation $\bdd{\ddtr \ddchi'}$, thus we apply proposition \ref{pro 9.4}.\ref{pro 9.4.b} and substitute the following 
\begin{align*}
\delta_o + |\os|/r_0 \leq c\epsilon,
\quad
\frakd_m + \frakd_o \leq c\epsilon, 
\quad
\ufrakd_m=0,
\quad
\ufrakd_o = \udelta_o,
\end{align*}
therefore
\begin{align*}
\Vert i \Vert^{n,p} 
\leq
&
\frac{c(n,p) (\epsilon + c \epsilon )}{r_0} c\epsilon 
+
\frac{c(n,p) ( \epsilon + c^2 \epsilon^2)}{r_0} \udelta_o
+
\frac{c(n,p) c\epsilon}{r_0} \udelta_o^2
+
\frac{c(n,p)}{r_0} (\udelta_o + c \epsilon) c \epsilon
\\
\leq
&
\frac{c(n,p)(1+ c + c^2) \delta}{r_0} \cdot \epsilon.
\end{align*}
The term $ii$ satisfies the following estimate by definition \ref{def 2.3}
\begin{align*}
\Vert ii \Vert^{n,p} \leq \frac{c(n,p) \epsilon}{r_0}.
\end{align*}
Therefore we obtain the estimate of $\fl{k+1}$
\begin{align*}
\Vert \fl{k+1} \Vert^{n+2,p}
\leq
c(n,p) [ 1 + (1+ c + c^2) \delta ] \epsilon.
\end{align*}
If we choose $c=2c(n,p)$ and $\delta$ suitably small that $ (1+ c + c^2) \delta< 1$, then $\fl{k+1}$ satisfies estimate \eqref{eqn 10.4}. Hence the induction argument proves the lemma.
\end{proof}
\begin{remark}
Lemma \ref{lem 10.1} requires the $L^{\infty}$ bounds of the metric components up to $(n+2)$-th order derivatives, and the structure coefficients up to $(n+1)$-th order derivatives. These requirements are the same as proposition \ref{pro 9.4}.
\end{remark}

\subsection{Convergence of the sequence of approximating solutions}
We shall prove that the sequence of approximating solutions $\{\fl{k}\}_{k\in\mathbb{N}}$ converges in $\mathrm{W}^{n+1,p}$. The idea is to show that it is a contractive sequence.
\begin{lemma}\label{lem 10.3}
Let $\ufl{s=0}$ be a function in $\dot{\mathrm{W}}^{n+2,p}(\udelta_o r_0) \cap \mathrm{M}(\udelta_m r_0)$, where $n \geq 2, p >2$ or $n\geq 3, p>1$. There exist a small positive constants $\delta$ depending on $n,p$, such that if $\epsilon, \udelta_o, \udelta_m$ are suitably bounded that
$\epsilon,\, \udelta_o,\, \epsilon \udelta_m \leq \delta$,
then the sequence of approximating solutions $\{\fl{k}\}_{k \in \mathbb{N}}$ is contractive.
\end{lemma}
\begin{proof}
Choose $\delta$ sufficiently small such that lemma \ref{lem 10.1} applies, then the sequence $\{\fl{k}\}_{k\in \mathbb{N}}$ exists and
$\Vert \fl{k} \Vert^{n+2,p} \leq c \epsilon r_0$.
Taking the difference of equation \eqref{eqn 10.2} in two successive steps, we obtain
\begin{align*}
\partial \bft[\ufl{s=0}, \fl{k+1}] \big( \fl{k+2} - \fl{k+1} \big)
=
\underbrace{- \Big[\bft\big(\ufl{s=0}, \fl{k+1}\big) - \bft\big(\ufl{s=0}, \fl{k}\big) - \partial \bft[\ufl{s=0}, \fl{k}] \big( \fl{k+1} - \fl{k} \big) \Big]}_{r.h.s.}.
\end{align*}
Note the right hand side is simply the error of the linearised perturbation $\bdd{\ddtr \ddchi'}$, thus we choose $\delta$ sufficiently small such that proposition \ref{pro 9.2} applies. Substituting
\begin{align*}
\delta_o + |\os|/r_0 \leq c \epsilon \leq c \delta,
\quad
\frakd_m + \frakd_o \leq \frac{\Vert \fl{k+1} - \fl{k} \Vert^{n+1,p}}{r_0} \leq 2 c \delta,
\quad
\udelta_o\leq \delta,
\quad
\ufrakd_m = \ufrakd_o =0,
\end{align*}
in estimate \eqref{eqn 9.5} of $\er{\ddtr \ddchi'}$ in proposition \ref{pro 9.2}, we obtain
\begin{align*}
\Vert r.h.s. \Vert^{n-1,p}
\leq
\frac{c(n,p)\delta}{r_0^2} \Vert \fl{k+1} - \fl{k} \Vert^{n+1,p}.
\end{align*}
Since $\partial\bft [\ufl{s=0}, \fl{k+1}]$ is invertible ensured in lemma \ref{lem 10.1}, we have
\begin{align*}
\frac{\Vert \fl{k+2} - \fl{k+1} \Vert^{n+1,p}}{r_0}
\leq
c(n,p)\delta \frac{ \Vert \fl{k+1} - \fl{k} \Vert^{n+1,p}}{r_0}
\end{align*}
Therefore we choose $\delta$ sufficiently small such that $c(n,p) \delta\leq \frac{1}{2}$, then the sequence $\{\fl{k}\}_{k\in \mathbb{N}}$ is contractive in $\mathrm{W}^{n+1,p}$.
\end{proof}
\begin{remark}
The requirements on the $L^{\infty}$ bounds of metric components and structure coefficients for above lemma are the same as for lemma \ref{lem 10.1}.
\end{remark}

\subsection{Limit of approximating solutions is a true solution}
By lemmas \ref{lem 10.1} and \ref{lem 10.3}, we obtain a sequence of approximating solutions $\{ \fl{k} \}$ of equation \eqref{eqn 10.1}, bounded in $\mathrm{W}^{n+2,p}$ and converging in $\mathrm{W}^{n+1,p}$. Denote the limit of $\{ \fl{k} \}$ by $\fl{\bft}$,
\begin{align*}
\fl{\bft} = \lim_{k\rightarrow +\infty} \fl{k}
\quad \text{ in } \mathrm{W}^{n+1,p}.
\end{align*}
Then $\fl{\bft}$ is also bounded in $\mathrm{W}^{n+2,p}$. We show that $\fl{\bft}$ solves equation \eqref{eqn 10.1}
\begin{align*}
\bft( \ufl{s=0}, \fl{\bft} ) = 0.
\tag{\ref{eqn 10.1}}
\end{align*}
\begin{lemma}\label{lem 10.5}
Under the assumptions of lemmas \ref{lem 10.1}, \ref{lem 10.3}, the limit $\fl{\bft}$ of the sequence $\{ \fl{k} \}_{k \in \mathbb{N}}$ solves equation \eqref{eqn 10.1}.
\end{lemma}
\begin{proof}
Since $\fl{k+1} - \fl{k} \rightarrow 0$ in $\mathrm{W}^{n+1,p}$, we have that by equation \eqref{eqn 10.2},
\begin{align*}
- \bft(\ufl{s=0}, \fl{k} ) = \partial \bft[\ufl{s=0}, \fl{k}] \big( \fl{k+1} - \fl{k} \big)
\quad
\rightarrow 0
\quad
\text{ in }
\mathrm{W}^{n-1,p}.
\end{align*}
By proposition \ref{pro 7.3}, we deduce that the map $\bft$ is continuous relative to the Sobolev norm $\Vert \cdot \Vert^{n+1,p}$. Therefore since $\fl{k} \rightarrow \fl{\bft}$ in $\mathrm{W}^{n+1,p}$ and $\bft(\ufl{s=0}, \fl{k}) \rightarrow 0$ in $\mathrm{W}^{n-1,p}$, we have
\begin{align*}
\bft( \ufl{s=0}, \fl{\bft} ) = \lim_{k\rightarrow +\infty} \bft( \ufl{s=0}, \fl{k} ) =0,
\end{align*}
i.e. $\fl{\bft}$ solves equation \eqref{eqn 10.1}. The lemma is proved.
\end{proof}

\subsection{Local uniqueness of the solution}
We show that the previous obtained solution $\fl{\bft}$ of equation \eqref{eqn 10.1} is unique in a small neighbourhood of $0$ in $\mathrm{W}^{n+2,p}$. We use the notation $\mathrm{W}^{n+2,p}(r)$ to denote the closed ball of radius $r$ at the origin $0$ in $\mathrm{W}^{n+2,p}$.
\begin{lemma}
Let $\ufl{s=0}$ be a function in $\dot{\mathrm{W}}^{n+2,p}(\udelta_o r_0) \cap \mathrm{M}(\udelta_m r_0)$, where $n \geq 2, p >2$ or $n\geq 3, p>1$. There exist a small positive constant $\delta$ and a constant $c$ depending on $n,p$, such that if $\epsilon, \udelta_o, \udelta_m$ are suitably bounded that
$\epsilon,\, \udelta_o,\, \epsilon \udelta_m \leq \delta$,
then the limit $\fl{\bft}$ of the sequence of approximating solutions $\{ \fl{k} \}_{k\in \mathbb{N}}$ is the unique solution of equation \eqref{eqn 10.1} $\bft(\ufl{s=0}, f)=0$ in $\mathrm{W}^{n+2,p}(c \delta r_0)$.
\end{lemma}
\begin{proof}
Let $\delta$ be as in lemmas \ref{lem 10.1}, \ref{lem 10.3}. The uniqueness follows from the invertibility of the linearised map $\partial \bft [\ufl{s=0}, f]$ in corollary \ref{coro 9.7} and error estimate in proposition \ref{pro 9.2}: let $f$, $f'$ be both solutions of equation \eqref{eqn 10.1} in $\mathrm{W}^{n+2,p}$, then
\begin{align*}
\Vert f' - f \Vert^{n+1,p} 
\leq &
c(n,p) r_0^2\, \leVe \partial \bft [\ufl{s=0}, f] ( f'- f ) \riVe^{n-1,p}
\\
= &
c(n,p) r_0^2\, \leVe \bft\big(\ufl{s=0}, f'\big) - \bft\big(\ufl{s=0}, f'\big) - \partial \bft [\ufl{s=0}, f] ( f'- f ) \riVe^{n-1,p}
\\
\leq &
c(n,p) \delta \Vert f' - f \Vert^{n+1,p}.
\end{align*}
The derivation of the above is identical to the derivation in the proof of lemma \ref{lem 10.3} with the same constants $c(n,p)$. Since $c(n,p) \delta \leq \frac{1}{2}$ as in lemma \ref{lem 10.3}, we conclude that $\Vert f' - f \Vert^{n+1,p}=0$. This implies the uniqueness of the solution. 
\end{proof}

\subsection{Parametrisation map of marginally trapped surfaces}
We already prove the existence and local uniqueness of the solution of equation \ref{eqn 10.1} $\bft\big(\ufl{s=0}, f \big)=0$ when $\ufl{s=0}$ is given. It naturally gives rise to a map from $\ufl{s=0}$ to the corresponding solution $f$.
\begin{definition}
Choose the positive constant $\delta$ as in lemmas \ref{lem 10.1}, \ref{lem 10.3}. For constants $\epsilon, \udelta_o, \udelta_m$ suitably bounded by
$\epsilon,\, \udelta_o,\, \epsilon \udelta_m \leq \delta$,
define the solution map $\bfs$ of equation \eqref{eqn 10.1}
\begin{align*}
\bfs:
\quad
\dot{\mathrm{W}}^{n+2,p}(\udelta_o r_0) \cap \mathrm{M}(\udelta_m r_0) \
\rightarrow \
\mathrm{W}^{n+2,p}(c \delta r_0),
\quad
\bft \big( \ufl{s=0}, \bfs ( \ufl{s=0} ) \big)=0.
\end{align*}
\end{definition}

It is clear now that the surface with the second parameterisation $\big( \ufl{s=0}, \bfs(\ufl{s=0} ) \big)$ has the vanishing outgoing null expansion. We shall show that it is a marginally trapped surface. It is sufficient to show that the incoming null expansion is negative. This follows from formula \eqref{eqn 4.1}, and estimate of the parameterisation function $\uh$ in proposition \ref{pro 3.2}.
\begin{lemma}\label{lem 10.8}
There exists a positive constant $\delta$, such that the spacelike surface $\Sigma$ with the second parameterisation $\big( \ufl{s=0}, \bfs ( \ufl{s=0} ) \big)$ has negative future incoming null expansion.
\end{lemma}
\begin{proof}
Recall formula \eqref{eqn 4.1} of the incoming null expansion $\dtr \duchi$,
\begin{align*}
\dtr \duchi =& 
\tr \uchi 
- 2\Omega^2 \slashDelta \uh - \Omega^2 \vert\dslashd \uh \vert_{\slashg}^2 \tr \chi -4\Omega^2 \ueta \cdot \dslashd \uh 
- 4 \Omega^2 \omega \vert \dslashd \uh \vert^2_{\slashg} + 4 \Omega^2 \chi ( \slashnabla \uh, \slashnabla \uh),
\tag{\ref{eqn 4.1}}
\end{align*}
Therefore
\begin{align*}
- \dtr \duchi 
\leq
&
- \tr \uchi_S + \vert \tr \uchi - \tr \uchi_S \vert
+ 2\Omega^2 \vert \slashDelta \uh \vert 
+ \Omega^2 \vert \dslashd \uh \vert_{\slashg}^2 \cdot \vert \tr \chi \vert
+ 4\Omega^2 \vert \ueta \vert \cdot \vert \dslashd \uh \vert
\\
&
+ 4 \Omega^2 \vert \omega \vert \cdot \vert \dslashd \uh \vert^2_{\slashg} 
+ 4 \Omega^2 \vert \chi \vert \cdot \vert \dslashd \uh \vert^2,
\end{align*}
Then by estimate of $\uh$ in proposition \ref{pro 3.2}, we obtain that
\begin{align*}
- \dtr \duchi
\leq
- \tr \uchi_S
+
\frac{c \delta}{r_0}
<
-\frac{2 e^{-0.1}}{1.1} \cdot \frac{1}{r_0} + \frac{c \delta}{r_0}.
\end{align*}
The last inequality follows from
\begin{align*}
\min_{M_{\kappa, \underline{\kappa}}} \tr \uchi_S \geq \tr \uchi_S(s=0.1 r_0,\us = -0.1 r_0),
\quad
r(s=0.1 r_0,\us = -0.1 r_0) \leq 1.1 r_0,
\end{align*} 
and the formula \eqref{eqn 2.2} of $\tr \uchi_S$. We choose $\delta$ sufficiently small such that $c \delta \leq \frac{2 e^{-0.1}}{1.1}$, then $\dtr \duchi$ is negative. The lemma is proved.
\end{proof}
Because of lemma \ref{lem 10.8}, we shall also call $\bfs$ the parameterisation map of marginally trapped surfaces. Thus we can summarise the result in this section in the following theorem, which is the main theorem of this paper on marginally trapped surfaces in a perturbed Schwarzschild spacetime.
\begin{theorem}\label{thm 10.9}
There exist a positive constant $\delta$ and a constant $c$, both depending on $n,p$ and satisfying $c\delta < \kappa$, such that if $\epsilon, \udelta_o, \udelta_m$ are suitably bounded that
$\epsilon,\, \udelta_o,\, \epsilon \udelta_m \leq \delta$,
then there exists a unique map $\bfs$,
\begin{align*}
\bfs:
\quad
\dot{\mathrm{W}}^{n+2,p} ( \udelta_o r_0) \cap \mathrm{M}(\udelta_m r_0)
\ \rightarrow \ 
\mathrm{W}^{n+2,p} (c\epsilon r_0)
\subset
\mathrm{W}^{n+2,p} (c \delta r_0),
\quad
\ufl{s=0}
\mapsto
\bfs ( \ufl{s=0} ),
\end{align*}
such that the spacelike surface $\Sigma$ with the second parameterisation $\big( \ufl{s=0}, \bfs (\ufl{s=0}) \big)$ is a marginally trapped surface, where the future outgoing null expansion vanishes and the future incoming null expansion is negative. We call $\bfs$ the parameterisation map of marginally trapped surfaces.
\end{theorem}
\begin{remark}
The above theorem requires the same on the $L^{\infty}$ bounds of metric components and structure coefficients as in lemma \ref{lem 10.1}: the $L^{\infty}$ bounds of the metric components up to $(n+2)$-th order derivatives, and the structure coefficients up to $(n+1)$-th order derivatives.
\end{remark}

We briefly explain the geometric meaning of the above theorem to conclude this section: given an incoming null hypersurface $\ucalH$ near some $\uC_{\us}$, the parameterisation map $\bfs$ tells where a marginally trapped surface $\Sigma$ lies in $\ucalH$.
\begin{figure}[H]
\begin{center}
\begin{tikzpicture}
\draw[dashed] (-1,0)
to [out=70, in=180] (0,0.5)
to [out=0,in=110] (1,0);
\draw (1,0)
to [out=-70,in=0] (0,-0.8)
to [out=180,in=-110] (-1,0); 
\node[below] at (0,-0.7) {\tiny $\Sigma_{0,0}$}; 
\draw[dashed] (-1,0) to [out=70,in=-110] (-0.7,0.9);
\draw (-1,0) to [out=-110,in=70] (-1.85,-2.4);
\draw[dashed] (1,0) to [out=110,in=-70] (0.7,0.9);
\draw[->] (1,0) to [out=-70,in=110] (1.85,-2.4) node[right] {\small $s$}; 
\node[above right] at (1.5,-2) {\tiny $\uC_{\us=0}$};
\draw[dashed] (-1,0) to [out=-45,in=135] (-0.5,-0.5);
\draw (-1,0) to [out=135, in= -45] (-2,1);
\draw[dashed] (1,0) to [out=-135,in=45] (0.5,-0.5);
\draw[->] (1,0) to [out=45, in= -135] (1.9,0.9) node[right] {\tiny $C_{s=0}$}to [out=45,in=-135] (2.3,1.3) node[right] {\small $\us$}; 
\draw[dashed] (-1.5,0.5) to [out=135,in=45] (1.5,0.5);
\draw (1.5,0.5) to [out=-135,in=0] (0,0) to [out=180,in=-45] (-1.5,0.5); 
\node[below] at (0,0.1) {\scriptsize $\Sigma_0$}; 
\draw[dashed] (-1.4,1.3) to [out=-110,in=70] (-1.6,0.7);
\draw (-1.6,0.7) to [out=-110,in=70] (-2.75,-2.4);
\draw[dashed] (1.4,1.3) to [out=-70,in=110] (1.6,0.7);
\draw[->] (1.6,0.7) to [out=-70,in=110] (2.75,-2.4) node[right]{\small $s$}; 
\node[right] at (2.4,-1.5) {\scriptsize $\ucalH$}; 
\draw[dashed] (-2.3,-2.2+1) to [out=70,in=180] (-0.5,-2.5+1) node[below] {\scriptsize $\Sigma$} to [out=0,in=110] (2.3,-2.2+1);
\draw (2.3,-2.2+1) to [out=-70,in=0] (1,-2+0.9) to [out=180,in=-110] (-2.3,-2.2+1); 
\draw[->] (-1.73,0.73) to [out=-110,in=70] (-2.45,-1.2);
\node[left] at (-2.3,-0.8) {$\bfs(\ufl{s=0})$};
\draw[->] (-1.05,-0.1) to [out=135,in=-45] (-1.65,0.5);
\node[below] at (-1.5,0.2) {\tiny $\ufl{s=0}$};
\end{tikzpicture}
\end{center}
\caption{The parameterisation map $\bfs$.}
\end{figure}

\section{Some properties of the parameterisation map}
In this section, we study the continuity of the parameterisation map $\bfs$ of marginally trapped surfaces, and discuss the linearisation of $\bfs$.

\subsection{Continuity of the parameterisation map}
Let $\ufl{a,s=0}, a=1,2$ be two functions in $\dot{\mathrm{W}}^{n+2,p} (\udelta_o r_0) \cap \mathrm{M}(\udelta_m r_0)$, and $\fl{a} = \bfs (\ufl{a,s=0} )$. Denote
\begin{align*}
&
\dd{\ufl{s=0}} = \bdd{\ufl{s=0}} = \ufl{2,s=0} - \ufl{1,s=0}, 
\\
&
\dd{f} = \fl{2} - \fl{1} = \bfs ( \ufl{2,s=0} ) - \bfs ( \ufl{1,s=0} ).
\end{align*}
We obtain an estimate of $\dd{f}$.
\begin{proposition}\label{pro 11.1}
Let $\ufl{a,s=0}, a=1,2$ be two functions in $\dot{\mathrm{W}}^{n+2,p} (\udelta_o r_0) \cap \mathrm{M}(\udelta_m r_0)$. Suppose their difference satisfies the following estimates
\begin{align*}
\leVe \dslashd \bdd{\ufl{s=0}} \riVe^{n,p}  \leq  \ufrakd_o r_0,
\quad
\leve \overline{\bdd{\ufl{s=0}}}^{\circg} \rive \leq \ufrakd_m r_0.
\end{align*}
There exist a positive constant $\delta$, and a constant $c(n,p)$ both depending on $n,p$, such that if $\epsilon, \udelta_o, \delta_m, \ufrakd_o, \ufrakd_m$ are suitably bounded that
$\epsilon,\, \udelta_o,\ \epsilon\udelta_m,\ \ufrakd_o,\ \epsilon \ufrakd_m \leq \delta$,
then $\dd{f} =\bfs \big( \ufl{2,s=0} \big) - \bfs \big( \ufl{1,s=0} \big)$ satisfies the estimate
\begin{align*}
\Vert \dd{f} \Vert^{n+1,p} \leq c(n,p) \epsilon \Vert \bdd{\ufl{s=0}} \Vert^{n+1,p}.
\end{align*}
Therefore the parametrisation map $\bfs$ is continuous from $\dot{\mathrm{W}}^{n+2,p} (\udelta_o r_0) \cap \mathrm{M}(\udelta_m r_0)$ to $\mathrm{W}^{n+2,p} ( c \epsilon r_0)$, but in the weaker Sobolev norm $\Vert \cdot \Vert^{n+1,p}$.
\end{proposition}
\begin{proof}
Let $\Sigmal{a}$ be the marginally trapped surface with the second parameterisation $\big( \ufl{a,s=0}, \fl{a} \big)$. Consider the perturbation of the outgoing null expansion from $\Sigmal{1}$ to $\Sigmal{2}$,
\begin{align*}
\er{\ddtr \ddchi'}
=&
\bft \big(  \ufl{2,s=0}, \fl{1} \big) - \bft \big(  \ufl{2,s=0}, \fl{1} \big) 
-
\partial \bft [\ufl{1,s=0}, \fl{1} ] \big( \bdd{\ufl{s=0}} \big)
-
\upartial \bft [\ufl{1,s=0}, \fl{1} ] \big( \dd{f} \big)
\\
=&
-
\partial \bft [\ufl{1,s=0}, \fl{1} ] \big( \bdd{\ufl{s=0}} \big)
-
\upartial \bft [\ufl{1,s=0}, \fl{1} ] \big( \dd{f} \big).
\end{align*}
Therefore we obtain that
\begin{align*}
\upartial \bft[\ufl{1,s=0}, \fl{1} ] ( \dd{f} )
=
- \er{\ddtr \ddchi'} 
- \partial \bft [\ufl{1,s=0}, \fl{1} ] \big( \bdd{\ufl{s=0}} \big),
\end{align*}
thus by corollary \ref{coro 9.7},
\begin{align*}
\Vert \dd{f} \Vert^{m+2,p}
\leq
c(n,p) r_0^2 \Big[ \Vert \er{\ddtr \ddchi'} \Vert^{m,p} + \leVe \partial \bft [\ufl{1,s=0}, \fl{1} ] \big( \bdd{\ufl{s=0}} \big) \riVe^{m,p} \Big].
\end{align*}
We apply proposition \ref{pro 9.2} to estimate $\er{\ddtr \ddchi'}$. Note that $\dd{f}$ satisfies a rough estimate from theorem \ref{thm 10.9},
\begin{align*}
\Vert \dd{f} \Vert^{n+2,p} \leq 2 \epsilon r_0.
\end{align*}
Therefore in estimate \eqref{eqn 9.5} of $\er{\ddtr \ddchi'}$, we can set that
\begin{align*}
\frakd_m + \frakd_o=2c \epsilon,
\quad
\delta_o + |\os|/r_0 \leq c \epsilon \leq c \delta' \leq \delta,
\quad
\ufrakd_m \leq 2 \underline{\kappa} \leq 0.2,\footnotemark
\quad
\ufrakd_o \leq 2 \delta,
\end{align*}
\footnotetext{Recall the notion of $\underline{\kappa}$ in definition \ref{def 2.1}.}
then
\begin{align*}
\Vert \er{\ddtr \ddchi'} \Vert^{n-1,p}
\leq
\frac{c(n,p)\delta}{r_0^2} \Vert \dd{f} \Vert^{n+1,p}
+
\frac{c(n,p) \epsilon }{r_0} (\ufrakd_m + \ufrakd_o ),
\end{align*}
For the term $\partial \bft[\ufl{1,s=0}, \fl{1}] \big( \bdd{\ufl{s=0}} \big)$, 
\begin{align*}
\leVe \big( \partial_{\us} \tr \chi'_S \big)|_{\Sigmal{1}} \riVe^{n,p}
\leq
\frac{c(n,p) \epsilon}{r_0^2},
\end{align*}
thus
\begin{align*}
\leVe \partial \bft[\ufl{1,s=0}, \fl{1}] \big( \bdd{\ufl{s=0}} \big) \riVe^{m,p}
=
\leVe \big( \partial_{\us} \tr \chi'_S \big)|_{\Sigmal{1}} \cdot \bdd{\ufl{s=0}} \riVe^{m,p}
\leq
\frac{c(n,p) \epsilon}{r_0^2} \Vert \bdd{\ufl{s=0}} \Vert^{m,p}.
\end{align*}
Therefore we obtain that
\begin{align*}
\Vert \dd{f} \Vert^{n+1,p}
\leq
c(n,p) \delta \Vert \dd{f} \Vert^{n+1,p} + c(n,p) \epsilon \Vert \bdd{\ufl{s=0}} \Vert^{n+1,p}.
\end{align*}
Choose $\delta$ sufficiently small that $c(n,p) \delta \leq \frac{1}{2}$, then we obtain that
\begin{align*}
\Vert \dd{f} \Vert^{n+1,p}
\leq
c(n,p) \epsilon \Vert \bdd{\ufl{s=0}} \Vert^{n+1,p},
\end{align*}
thus the proposition is proved.
\end{proof}
\begin{remark}
The above proposition shows the continuity of the parameterisation map $\bfs$ in the weaker Sobolev norm. The loss of regularity follows directly from proposition \ref{pro 9.2} on the estimate of $\er{\ddtr \ddchi'}$. In the special case that $\ufl{1,s=0}$ or $\ufl{2,s=0}$ is constant, we have improved estimate of $\er{\ddtr \ddchi'}$ in proposition \ref{pro 9.4}, thus we can show the continuity of $\bfs$ at such constant functions without loss of regularity.
\end{remark}
\begin{proposition}\label{pro 11.3}
Under the same setting of proposition \ref{lem 10.3}, we assume further that $\ufl{1,s=0}$ or $\ufl{2,s=0}$ is constant, and
\begin{align*}
\Vert \dslashd \bdd{\ufl{s=0}} \Vert^{n+1,p} \leq \ufrakd_o r_0.
\end{align*}
Then we can improve the conclusion of proposition \ref{lem 10.3} to that $\dd{f}$ satisfies the estimate
\begin{align*}
\Vert \dd{f} \Vert^{n+2,p} \leq c(n,p) \epsilon \Vert \bdd{\ufl{s=0}} \Vert^{n+2,p}.
\end{align*}
Therefore the parameterisation map $\bfs$ is continuous at $\ufl{s=0} \equiv \mathrm{const.}$ in the Sobolev norm $\Vert \cdot \Vert^{n+2,p}$.
\end{proposition}
\begin{proof}
The proof follows the same path as in the proof of proposition \ref{pro 11.1}, simply replacing the estimate of $\er{\ddtr \ddchi'}$ using proposition \ref{pro 9.4}.
\end{proof}

\subsection{A linearisation of the parameterisation map}
As the parameterisation map $\bfs$ is the solution map of equation \eqref{eqn 10.1} $\bft\big( \ufl{s=0}, \bfs\big( \ufl{s=0} \big) \big)=0$, it is natural to use the linearisation of $\bft$ to construct a linearisation $\delta \bfs$ of $\bfs$: formally
\begin{align*}
0 = \delta \bft = \upartial \bft \big( \bdd{\ufl{s=0}} \big) + \partial \bft \circ  \delta \bfs \big( \bdd{\ufl{s=0}} \big)
\
\Rightarrow
\
\delta \bfs \big( \bdd{\ufl{s=0}} \big) 
=
- ( \partial \bft )^{-1} \circ \upartial \bft \big( \bdd{\ufl{s=0}} \big),
\end{align*}
thus we use the above formula as a linearisation of the parameterisation map $\bfs$.

\begin{definition}\label{def 11.4}
For the parameterisation map $\bfs$ of marginally trapped surfaces in theorem \ref{thm 10.9}, define the following linear map as the linearisation of $\bfs$ at $\ufl{s=0}$, denoted as $\delta \bfs [\ufl{s=0}]$: let $f=\bfs \big( \ufl{s=0} \big)$, then
\begin{align*}
\delta \bfs [\ufl{s=0}] \big( \bdd{\ufl{s=0}} \big) 
=
- ( \partial \bft [\ufl{s=0}, f])^{-1} \circ \upartial \bft [\ufl{s=0}, f] \big( \bdd{\ufl{s=0}} \big).
\end{align*}
\end{definition}
Note by corollary \ref{coro 9.7}, the above inverse operator $( \partial \bft [\ufl{s=0}, f])^{-1}$ is well-defined. We estimate the operator norm of the linearisation map $\delta \bfs$.
\begin{proposition}\label{proposition 11.5}
Let $\ufl{s=0}$ be a function in $\dot{\mathrm{W}}^{n+2,p} (\udelta_o r_0) \cap \mathrm{M}(\udelta_m r_0)$, and $f=\bfs\big(\ufl{s=0} \big)$ that $\Vert f \Vert^{n+2,p} \leq c\epsilon r_0$. There exists a constant $c(n,p)$ that for $0\leq m \leq n$,
\begin{align*}
\leVe \delta \bfs [\ufl{s=0}] \big( \bdd{\ufl{s=0}} \big) \riVe^{m+2,p}
\leq
c(n,p) \epsilon \Vert \bdd{\ufl{s=0}} \Vert^{m+2,p}.
\end{align*}
\end{proposition}
\begin{proof}
Let $\Sigma$ be the marginally trapped surface with the second parameterisation $\big( \ufl{s=0}, f \big)$. We have that
\begin{align*}
\leVe \big( \partial_{\us} \tr \chi'_S \big)|_{\Sigma} \riVe^{n,p}
\leq
\frac{c(n,p) \epsilon}{r_0^2},
\end{align*}
thus
\begin{align*}
\leVe \upartial \bft [ \ufl{s=0}, f ] \big( \bdd{\ufl{s=0}} \big) \riVe^{m,p}
\leq
\frac{c(n,p) \epsilon}{r_0^2} \Vert \bdd{\ufl{s=0}} \Vert^{m,p},
\end{align*}
therefore by corollary \ref{coro 9.7},
\begin{align*}
\leVe \delta \bfs [\ufl{s=0}] \big( \bdd{\ufl{s=0}} \big) \riVe^{m+2,p}
\leq
c(n,p) r_0^2 \leVe \upartial \bft [ \ufl{s=0}, f ] \big( \bdd{\ufl{s=0}} \big) \riVe^{m,p}
\leq
c(n,p) \epsilon \Vert \bdd{\ufl{s=0}} \Vert^{m,p}.
\end{align*}
\end{proof}

\subsection{Error of the linearisation of the parameterisation map}
In this section, we study the error of the linearisation $\bfs$ constructed in above subsection. Let $\ufl{a,s=0}, a=1,2$ be two functions in $\dot{\mathrm{W}}^{n+2,p} (\udelta_o r_0) \cap \mathrm{M}(\udelta_m r_0)$, and $\fl{a} = \bfs(\ufl{a,s=0} \big)$. We denote $\delta \bfs [ \ufl{1,s=0} ] \big( \bdd{\ufl{s=0}} \big)$ by $\bdd{f}$, thus the error of the linearisation map $\delta \bfs$, denoted by $\er{\bfs} [ \ufl{1,s=0} ] \big( \bdd{\ufl{s=0}} \big)$, is the difference of $\dd{f}$ with $\bdd{f}$
\begin{align*}
\er{\bfs} [ \ufl{1,s=0} ] \big( \bdd{\ufl{s=0}} \big)
=&
\dd{f} - \bdd{f}
\\
=&
\bfs \big( \ufl{2,s=0} \big) - \bfs \big( \ufl{1,s=0} \big)
-
\delta \bfs [ \ufl{1,s=0} ] \big( \bdd{\ufl{s=0}} \big).
\end{align*}

We can estimate the error of the linearisation map $\delta \bfs$. For the sake of brevity, use $\bfe$ to denote $\er{\bfs} [\ufl{1,s=0} ] \big( \bdd{\ufl{s=0}} \big)$. Note that
\begin{align*}
\er{\ddtr \ddchi'}
=
-
\partial \bft [\ufl{1,s=0}, \fl{1} ] \big( \bdd{\ufl{s=0}} \big)
-
\upartial \bft [\ufl{1,s=0}, \fl{1} ] \big( \dd{f} \big),
\end{align*}
and
\begin{align*}
0 = 
\partial \bft [\ufl{1,s=0}, \fl{1} ] \big( \bdd{\ufl{s=0}} \big)
+
\upartial \bft [\ufl{1,s=0}, \fl{1} ] \big( \bdd{f} \big),
\end{align*}
thus we obtain that
\begin{align*}
\upartial \bft[\ufl{1,s=0}, \fl{1} ] \big( \bfe \big) = - \er{\ddtr \ddchi'}.
\end{align*}
The above equation gives the following estimate of the error $\bfe$
\begin{align*}
\Vert \bfe \Vert^{n+1,p} \leq c(n,p) \epsilon \Vert \bdd{\ufl{s=0}} \Vert^{n+1,p},
\end{align*}
which is at the same order as $\dd{f}$ and $\bdd{f}$. Thus from the point of view of analysis, the linearisation map $\delta \bfs$ is not a good linearisation of $\bfs$, although by more careful analyses, we can show that the constant $c(n,p)$ above could be finer than the constant in the estimate of $\dd{f}$. The reason is that when we construct the linearisation of $\ddtr \ddchi'$ in section \ref{sec 9}, we allow the error of such sizes: $\epsilon \cdot (\ufrakd_m + \ufrakd_o)$. Such kind of error comes from $\dd{\hir{1}{\ddtr \ddchi'}}$ by replacing the geometric quantities in $(M,g)$ by the corresponding quantities in the Schwarzschild spacetime. 

The simplified construction of the linearisation of $\dd{\ddtr \ddchi'}$ in section \ref{sec 9} is sufficient for the construction of marginally trapped surfaces, while it is too rough to give an appropriate linearisation of the parameterisation map $\bfs$. To construct an appropriate linearisation of $\bfs$, we shall use the exact geometric quantities in the perturbed Schwarzschild spacetime $(M,g)$ rather than their approximations in the Schwarzschild spacetime. In this way, we will improve the estimate of the error $\dd{\ddtr \ddchi'}$ by improving the terms $\epsilon \ufrakd_m, \epsilon \ufrakd_o$ in estimate \eqref{eqn 9.5}. However this goes beyond the scope of this paper, thus we stop here.

\section*{Acknowledgements}
\noindent
This paper generalises the result in the author's thesis \cite{L2} on marginally trapped surfaces in a perturbed Schwarzschild black hole. The author is grateful to Demetrios Christodoulou for his constant encouragement and generous guidance. The author also thanks Alessandro Carlotto and Lydia Bieri for many helps on the refinements of the manuscript.

\appendix
\section{Derivation of equation \eqref{eqn 3.3}}\label{appen eqn 3.3}
We adopt the notations in section \ref{sec 3}. From formulae for the normal null vectors and tangent vectors of $S_t$ in \cite{L4}, we have
\begin{align*} 
&
\uL_{S_t}= \uL + \uvarepsilon L + \uvarepsilon^i \partial_i = \partial_s + \uvarepsilon \partial_{\us} + \left( \uvarepsilon^i + \uvarepsilon b^i \right) \partial_i,
\\
&
\partial_{i,S_t} = \partial_i + \flt_i \partial_s + \uflt_i \partial_{\us} = \left( \delta_i^j - \flt_i b^j \right) \partial_j + \flt_i \uL + \uflt_i L,
\end{align*}
where $\uflt_i, \flt_i$ are the partial derivatives of $\uflt, \flt$ and $\uvarepsilon$ is given by the following formulae
\begin{align*}
&
\underline{\varepsilon}^k = \underline{e}^k + \underline{\varepsilon} e^k,
\quad
\underline{\varepsilon} = \frac{ -|\underline{e}|^2}{(2\Omega^2 + e\cdot \underline{e}) + \sqrt{(2\Omega^2 + e\cdot \underline{e})^2 -|e|^2 |\underline{e}|^2}},
\\
&
|e|^2 = \slashg_{ij}e^ie^j,
\quad
|\underline{e}|^2 = \slashg_{ij} \underline{e}^i \underline{e}^j, 
\quad
e\cdot \underline{e} =\slashg_{ij} e^i \underline{e}^j,
\\
&
e^k =-2\Omega^2 \cdot \flt_i \left(B^{-1}\right)_j^i \left(\slashg^{-1}\right)^{jk}, 
\quad
\underline{e}^k = -2\Omega^2 \cdot \uflt_i \left(B^{-1}\right)_j^i \left( \slashg^{-1} \right)^{jk},
\\
&
B_i^j= \delta_i^j - \flt_i b^j.
\end{align*}

On the other hand, since the family of surfaces $\{ S_t \}$ has the first parameterisation $(\uflt, \flt)$ for $S_t$, then the deformation vector field $V_t$ for the family is
\begin{align*}
V_t = \partial_t \uflt \cdot \partial_{\us} + \partial_t \flt \cdot \partial_s.
\end{align*}
Since $\{ S_t \}$ is embedded in the incoming null hypersurface $\ucalH$, then $V_t$ is tangential to $\ucalH$, thus $V_t$ is a linear combination of $\uL_{S_t}, \partial_{i, S_t}$. Assume
\begin{align*}
V_t = c \uL_{S_t} + c^i \partial_{i,S_t},
\end{align*}
then we obtain the following system of linear equations for $c, c^i$
\begin{align*}
\left\{
\begin{aligned}
&
\partial_t \uflt = c \uvarepsilon + c^i \cdot \uflt_i,
\\
&
\partial_t \flt = c + c^i \cdot \flt_i,
\\
&
0 = c \left( \uvarepsilon^i + \uvarepsilon b^i \right) + c^i.
\end{aligned}
\right.
\end{align*}
The above system is overdetermined. We solve $c,c^i$ from the last two equations and substitute to the first equation, then we obtain an equation for $\uflt$ which is the compatible condition for this system having a solution. We obtain that
\begin{align*}
\left\{
\begin{aligned}
&
\partial_t \flt = \left[ 1- \left( \uvarepsilon^i + \uvarepsilon b^i \right) \flt_i \right] c,
\\
&
\partial_t \uflt = \left[ \uvarepsilon - \left( \uvarepsilon^i + \uvarepsilon b^i \right) \uflt_i \right] c.
\end{aligned}
\right.
\end{align*}
Therefore, substituting $\flt= t f$, we derive equation \eqref{eqn 3.3}
\begin{align*}
\partial_t \uflt = f \cdot \left[ 1- \left( b^i + \uvarepsilon^i \right) t f_i \right]^{-1} \cdot \left[ \uvarepsilon - \left( b^i + \uvarepsilon^i \right) \partial_i \uflt \right].
\tag{\ref{eqn 3.3}}
\end{align*}

\section{Proof of proposition \ref{pro 3.3}}\label{appen pro 3.3}
We fill the details to complete the proof of proposition \ref{pro 3.3}. It is sufficient to prove the estimates \eqref{eqn 3.8} of $F, \Xlt, \relt$ to verify the proof sketch.
\begin{proof}
Following the proof sketch, we assume that
$\delta \leq \frac{1}{2},
(c_o + c_{m,m} + c_{m,o} ) \delta \leq 1$,
and $\uflt$ satisfies estimates \eqref{eqn 3.6} for $t \in [0,t_a]$. 

In order to estimate $F, \Xlt, \relt$, we estimate $ f$, $t b^i f_i$, $t e^i f_i$, $t \ue^i f_i$, $\uvarepsilon$, $b^i\, \uflt_i$, $e^i\, \uflt_i$, $\ue^i\, \uflt_i$ on $S_t$ first. Here we need the following estimate of $\ddcircnabla a$ where $a$ is a background quantity like $\vec{b}, \Omega, \slashg$ and their derivatives with respect to $\circnabla, \partial_s, \partial_{\us}$,
\begin{align*}
&
\ddcircnabla_i a = \circnabla_i a + \uflt_i \cdot \partial_{\us} a + t f_i \cdot \partial_s a,
\\
&
\Vert \ddcircnabla a \Vert^{w,p} 
\leq
\Vert \circnabla_i a \Vert^{w,p}_{S_t}
+
c(n,p) \Vert \ddslashd \uflt \Vert^{n,p} \cdot \Vert \partial_{\us} a \Vert^{w,p}_{S_t}
+ 
c(n,p) t \Vert \ddslashd f \Vert^{n,p} \cdot \Vert \partial_s a \Vert^{w,p}_{S_t}
\\
&\hspace{40pt}
\leq
\Vert \circnabla_i a \Vert^{w,p}_{S_t}
+
c(n,p) \ud_o r_0 \cdot \Vert \partial_{\us} a \Vert^{w,p}_{S_t}
+ 
c(n,p) t \delta_o r_0 \cdot \Vert \partial_s a \Vert^{w,p}_{S_t},
\end{align*}
where $w \leq n$. Therefore by induction arguments, we can show that for $w \leq n$
\begin{align*}
& 
\Vert \circnabla^k \partial_{s}^l \vec{b} \Vert^{w+1,p}_{S_t} 
\leq 
\frac{c(n,p,k,l)\epsilon }{r_0} \left( \ud_o +  \ud_m \right),
\\
& 
\Vert \circnabla^k \partial_s^l \partial_{\us}^m \vec{b} \Vert^{w+1,p}_{S_t} \leq \frac{c(n,p,k,l,m)\epsilon}{r_0},
\end{align*}
and
\begin{align*}
\Vert \circnabla^k \partial_s^l \partial_{\us}^m \Omega \Vert^{w+1,p}_{S_t} \leq c(n,p,k,l,m),
\quad
\Vert \circnabla^k \partial_s^l \partial_{\us}^m \slashg \Vert^{w+1,p}_{S_t} \leq c(n,p,k,l,m) r_0^2.
\end{align*}
Note that when estimate $\Vert \circnabla^k \partial_s^l \partial_{\us}^m a \Vert^{w+1,p}_{S_t}$, we use the $L^{\infty}$ bounds of $a$ and its derivatives with respect to $\circnabla, \partial_s, \partial_{\us}$ up to $(k+l+m+w+1)$-th order.

Now we can list the estimates of $ f$, $t b^i f_i$, $t e^i f_i$, $t \ue^i f_i$, $\uvarepsilon$, $b^i\, \uflt_i$, $e^i\, \uflt_i$, $\ue^i\, \uflt_i$ in $F$
\begin{align*}
&
\Vert f \Vert^{n+2,p} \leq ((4\pi)^{1/p} |\os|/r_0 + \delta_o) r_0, 
\\
&
\Vert t b^i f_i \Vert^{n+1,p}
\leq
c(n,p) \epsilon \left( \ud_o + \ud_m \right) \delta_o  t,
\\
&
\Vert t e^i f_i \Vert^{n+1,p}
\leq
c(n,p) \delta_o^2 t^2,
\\
&
\Vert t \ue^i f_i \Vert^{n,p}, \Vert e^i\, \uflt_i \Vert^{n,p}
\leq
c(n,p) \delta_o \ud_o t,
\\
&
\Vert \uvarepsilon \Vert^{n,p}, \Vert \ue^i\, \uflt_i \Vert^{n,p}
\leq
c(n,p) \ud_o^2 
\\
&
\Vert b^i\, \uflt_i \Vert^{n,p}
\leq
c(n,p) \epsilon \left( \ud_o + \ud_m \right) \ud_o.
\end{align*}
Therefore we obtain the estimate \eqref{eqn 3.8} of $F$ in the proof sketch,
\begin{align*}
\vert F \vert, \Vert F \Vert^{n,p}
\leq
c(n,p) (|\os|/r_0+ \delta_o) \left( \ud_o + \epsilon \ud_m \right) \ud_o r_0.
\end{align*}
By a similar argument, we obtain the estimate \eqref{eqn 3.8} of $\Xlt$ in the proof sketch
\begin{align*}
\Vert \Xlt \Vert^{n,p}
\leq
c(n,p) (|\os|/r_0+ \delta_o) \left( \ud_o + \epsilon \ud_m \right).
\end{align*}
Since $\Xlt$ only involves the first order derivatives of $\uflt$ and $f$, we can estimate $\Xlt$ up to $n$th derivatives in the above.

The estimate for $\relt$ is a bit more involved. We observe that the top order derivatives in $\relt$ are 2nd order derivatives of $\uflt$ and 3rd order derivatives of $f$. We have
\begin{align*}
\big\Vert \ddcircnabla^2 \uflt \big\Vert^{n-1,p} \leq c_o \udelta_o r_o,
\quad
\big\Vert \ddcircnabla^3 f \big\Vert^{n-1,p} \leq \delta_o r_0,
\end{align*}
then applying the following estimate on the product to the terms in $\relt$,
\begin{align*}
\Vert h_1 h_2 \Vert^{n-1,p} \leq c(n,p) \Vert h_1 \Vert^{n-1,p} \Vert h_2 \Vert^{n-1,p},
\quad
n-1\geq 1, p>2 \text{ or } n-1\geq 2, p>1,
\end{align*}
we obtain the estimate of $\Vert \relt \Vert^{n-1,p}$
\begin{align*}
\Vert \relt \Vert^{n-1,p} 
\leq 
c(n,p) (|\os|/r_0+ \delta_o) \left( \ud_o + \epsilon \ud_m \right) \ud_o r_0.
\end{align*}
Therefore we obtain the estimates \eqref{eqn 3.8} of $F, \Xlt, \relt$ in the proof sketch. The rest of the proof proceeds as in the proof sketch.

To conclude the proof, we determine up to which order derivatives the $L^{\infty}$ bounds of the metric components $\vec{b}, \Omega, \slashg$ are needed. In the estimates, we need the following norms of the metric components
\begin{align*}
&
\Vert \vec{b} \Vert_{S_t}^{n,p}, \quad \Vert \Omega \Vert_{S_t}^{n,p}, \quad \Vert \slashg \Vert_{S_t}^{n,p},
\\
&
\Vert \mathrm{D} a  \riVe_{S_t}^{n-1,p},
\
\Vert \mathrm{D}_1 \mathrm{D}_2 a  \riVe_{S_t}^{n-1,p},
\ \mathrm{D}, \mathrm{D_1}, \mathrm{D_2} = \circnabla, \partial_s, \partial_{\us}, 
\ a= \vec{b}, \Omega, \slashg,
\end{align*}
therefore the $L^{\infty}$ bounds of $\vec{b}, \Omega, \slashg$ up to $(n+1)$-th order derivatives are required in the proof.
\end{proof}

\section{Basics of rotational vector field derivatives}\label{appen R}
We construct the rotational vector fields on $(M,g_S)$ first. Consider an isometric embedding of $(\Sigma_{0,0}, \circg=r_0^{-2} \slashg_{S})$ to the 3-dimensional Euclidean space centring at the origin. Denote the rotational vector fields on $\Sigma_{0,0}$ by $R_i$ where
\begin{align*}
R_i = \sum_{j,k}\epsilon_{ijk} x_j \partial_k.
\end{align*}
$\epsilon_{ijk}$ is the permutation symbol of $\{1,2,3\}$.
Then extend $R_i$ to $M$ via the diffeomorphisms generated by $\partial_s$ and $\partial_{\us}$, i.e.
\begin{align*}
\lie_{\partial_s} R_i = \lie_{\partial_{\us}} R_i =0.
\end{align*}
Also extend functions $x_1, x_2, x_3$ to $M$ by requiring $\partial_s x_i = \partial_{\us} x_i =0$. We can also extend $R_i$ to any surface $\Sigma$ with first parameterisation $(\uf,f)$ since the parameterisation map $\vartheta \mapsto (\us,s,\vartheta) = ( \uf(\vartheta), f(\vartheta), \vartheta)$ pushes forward $R_i$ on $\Sigma_{0,0}$ to $\Sigma$.

We introduce the rotational vector components of a tensor: for a vector $\vec{v}$ and a one-form $\omega$, its $R_i$ component is\footnote{Use $R$ in the superscript or subscript to indicate it is the rotational vector component, instead of the coordinate component.}
\begin{align*}
v^{R,i} = \circg( \vec{v} , R_i ),
\quad
\omega_{R,i} = \omega(R_i),
\end{align*}
and for a general tensor $T$, its component can be defined inductively using the above rules, or explicitly
\begin{align*}
T_{R,i_1\cdots i_k}^{R,j_1,\cdots j_l} = \circg(T(R_{i_1}, \cdots, R_{i_k}), R_{j_1} \otimes \cdots \otimes R_{j_l}).
\end{align*}
Moreover, we introduce the mixed components of a tensor, i.e. partial coordinate components and partial rotational vector field components: for a tensor $T$, an example of mixed components is
\begin{align*}
T_{R,i_1,\cdots, \overline{i_s}, \cdots, i_k}^{R,j_1,\cdots, \overline{j_r}, \cdots, j_l}
=
\circg(T(R_{i_1}, \cdots, \partial_{i_s}, \cdots R_{i_k}), R_{j_1} \otimes \cdots \otimes \partial_{j_s} \otimes \cdots \otimes R_{j_l}),
\end{align*}
where we use overlined indices to denote the coordinate component indices.

We list the following properties on the rotational vector fields.
\begin{enumerate}
\item Contraction and inner product via rotational vector components:
\begin{align*}
&
T_i^i = T_{R,i}^{R,i},
\quad
\omega_i v^i = \omega_{R,i} \cdot v^{R,i},
\\
&
\circg( \vec{v}, \vec{w} ) 
= \circg_{R,ij} v^{R,i} w^{R,j}
= \sum_{i=1,2,3} v^{R,i} w^{R,i},
\\
&
\circg^{-1}(\omega, \sigma) 
= \big( \circg^{-1} \big)^{R,ij} \omega_{R,i} \sigma_{R,j}
= \sum_{i=1,2,3} \omega_{R,i} \sigma_{R,i},
\end{align*}

\item Let $\partial_{x_i}^{\parallel}$ be the orthogonal projection of $\partial_{x_i}$ to the tangent plane of $\mathbb{S}^2$.
\begin{align*}
\sum_{i=1,2,3} x_i R_i = 0,
\quad
R_i x_j = - \sum_{k=1,2,3} \epsilon_{ijk} x_k,
\quad
\partial_{x_i}^{\parallel} = - \sum_{j,k=1,2,3} \epsilon_{ijk}x_j R_k.
\end{align*}

\item Lie brackets and covariant derivatives of rotational vector fields:
\begin{align*}
[ R_i, R_j ] = - \sum_{k=1,2,3} \epsilon_{ijk} R_k,
\quad
\circnabla_{R_i} R_j = x_j \partial_{x_i}^{\parallel} = - \sum_{j,k=1,2,3} \epsilon_{imn} x_j x_m R_n.
\end{align*}

\item Lie derivatives with respect to $R_i$: for a vector field $\vec{v}$,
\begin{align*}
\left( \lie_{R_i} v \right)^{R,j}= [R_i ,v]^{R,j} 
= R_i ( v^{R,j} ) + \sum_{k=1,2,3} \epsilon_{ijk} v^{R,k};
\end{align*}
for a differential one form $\omega$,
\begin{align*}
\left( \lie_{R_i} \omega \right)_{R,j} 
= R_i (\omega_{R,j}) + \sum_{k=1,2,3} \epsilon_{ijk} \omega_{R,k};
\end{align*}
for a general tensor field $T$,
\begin{align*}
\left( \lie_{R_i} T \right)_{R,j_1 \cdots j_p}^{R, k_1 \cdots k_q}
=
R_i \left( T_{R, j_1 \cdots j_p}^{R, k_1 \cdots k_q}  \right)
+
\sum_{\substack{s=1,\cdots,p \\ t=1,2,3}} \epsilon_{i j_s t} T_{R,j_1 \cdots \underset{\hat{j_s}}{t} \cdots j_p}^{R, k_1 \cdots k_q}
+
\sum_{\substack{l=1,\cdots,q \\ r=1,2,3}} \epsilon_{i k_l r} T_{R,j_1 \cdots j_p}^{R, k_1 \cdots \overset{\hat{k_l}}{r} \cdots k_q}.
\end{align*}

\item Covariant derivatives with respect to $R_i$
\begin{align*}
&
( \circnabla_{R_i} v )^{R,j}
= R_i ( v^{R,j} ) + \sum_{m,n=1,2,3} \epsilon_{imn} x_j x_m v^{R,n};
\\
&
( \circnabla_{R_i} \omega )_{R,j}
= R_i ( \omega_{R,j} ) + \sum_{m,n=1,2,3} \epsilon_{imn} x_j x_m \omega_{R,n};
\end{align*}
and
\begin{align*}
\big( \circnabla_{R_i} T \big)_{R,j_1 \cdots j_p}^{R, k_1 \cdots k_q}
=&
R_i \left( T_{R, j_1 \cdots j_p}^{R, k_1 \cdots k_q}  \right)
+
\sum_{\substack{s=1,\cdots,p \\ m,t=1,2,3}} \epsilon_{i m t} x_{j_s} x_m T_{R,j_1 \cdots \underset{\hat{j_s}}{t} \cdots j_p}^{R, k_1 \cdots k_q}
\\
&
+
\sum_{\substack{l=1,\cdots,q \\ m, r=1,2,3}} \epsilon_{i m r} x_{k_l} x_m T_{R,j_1 \cdots j_p}^{R, k_1 \cdots \overset{\hat{k_l}}{r} \cdots k_q}.
\end{align*}

\item Hessian of a function in terms of $R_i$:
\begin{align*}
( \circnabla^2 f )_{R,ij}
=&
\circnabla_{R_i} \circnabla_{R_j} f
=
R_i R_j f + \sum_{m,n=1,2,3} \epsilon_{imn} x_j x_m R_n f.
\end{align*}

\item Sobolev norms: for a function $h$ on $\Sigma_{\us,s}$, define
\begin{align*}
\Vert h \Vert^{n,p}_{R} = \sum_{\substack{k=0,\cdots, n\\ i_1, \cdots ,i_k =1,2,3}}\left\{ \int_{\Sigma_{\us,s}} \leve R_{i_1} \cdots R_{i_k} h \rive^{p} \dvol_{\circg} \right\}^{1/p},
\quad
p\geq 1,
\end{align*}
then there exists a constant $c(n,p)$ such that
\begin{align*}
c(n,p)^{-1}\Vert h \Vert_R^{n,p} 
\leq 
\Vert h \Vert^{n,p} 
\leq
c(n,p) \Vert h \Vert_R^{n,p}.
\end{align*}
More generally, for a $(r,s)$-type tensor field $T$ on $\Sigma_{\us,s}$, define the norm
\begin{align*}
\Vert T \Vert^{n,p}_{R} 
= 
\sum_{\substack{l=0,\cdots, n\\ i_1, \cdots ,j_1 \cdots, k_1, \cdots=1,2,3}}\left\{ \int_{\Sigma_{\us,s}} \leve R_{i_1} \cdots R_{i_l} \left( T_{R, j_1 \cdots j_s}^{R, k_1 \cdots k_r}  \right)
\rive^{p} \dvol_{\circg} \right\}^{1/p},
\quad
p\geq 1,
\end{align*}
then there exists a constant $c(n,p,r,s)$ such that
\begin{align*}
c(n,p,r,s)^{-1}\Vert T \Vert_R^{n,p} 
\leq 
\Vert T \Vert^{n,p} 
\leq
c(n,p,r,s) \Vert T \Vert_R^{n,p}.
\end{align*}

\end{enumerate}

\section{Proof of proposition \ref{pro 6.1}}\label{appen pro 6.1}
It is sufficient to supply the proof of estimates \eqref{eqn 6.5} of $\dd{\Flt}$, $\dd{\Xlt}$, $\dd{\relt}$ in order to complete the proof of proposition \ref{pro 6.1}.
\begin{proof}
Estimates \eqref{eqn 6.5} can be obtained heuristically from estimates \eqref{eqn 3.8} of $\Flt$, $\Xlt$, $\relt$ in the proof of proposition \ref{pro 3.3} by the following procedure: introducing the rules
\begin{align*}
&
\begin{aligned}
&
|\os|/r_0
\rightarrow 
\frakd_m,
&&
\delta_o 
\rightarrow 
\frakd_o,
\\
&
\ud_m 
\rightarrow
\ubfd_m,
&&
\ud_o 
\rightarrow 
\ubfd_o,
\end{aligned}
\\
&
\epsilon^k r_0^l
\rightarrow 
\epsilon^k r_0^l \big( \ubfd_m + \ubfd_o + \frakd_m + \frakd_o \big),
\quad
k\in \mathbb{N}, l \in \mathbb{Z},
\end{align*}
and the following Leibniz rule.
\begin{center}
\begin{tikzcd}
\left[ \epsilon^k r_0^l \big( \ubfd_m + \ubfd_o + \frakd_m + \frakd_o \big) \right] \big( |\os|/r_0 \big)^i \delta_o^j \ud_m^r \ud_o^s
&
\epsilon^k r_0^l \left[ \big( |\os|/r_0 \big)^{i-1} \frakd_m \right] \delta_o^j \ud_m^r \ud_o^s
\\
\boxed{\epsilon^k r_0^l \big( |\os|/r_0 \big)^i \delta_o^j \ud_m^r \ud_o^s}
\arrow[u,"\epsilon^k r_0^l"]
\arrow[ur,"|\os|/r_0"]
\arrow[r,"\delta_o"]
\arrow[dr,"\ud_m"]
\arrow[d,"\ud_o"]
&
\epsilon^k r_0^l \big( |\os|/r_0 \big)^i \left[ \delta_o^{j-1} \frakd_o \right] \ud_m^r \ud_o^s
\\
\epsilon^k r_0^l \big( |\os|/r_0 \big)^i \delta_o^j \ud_m^r \left[ \ud_o^{s-1} \ubfd_o \right]
&
\epsilon^k r_0^l \big( |\os|/r_0 \big)^i \delta_o^j \left[ \ud_m^{r-1} \ubfd_m \right] \ud_o^s
\end{tikzcd}
\end{center}
Following the above pattern, estimates \eqref{eqn 6.5} can be deduced from estimates \eqref{eqn 3.8}. Check that the above pattern is valid for the perturbations of the background metric components $\vec{b}, \Omega, \slashg$ and their derivatives. The pattern is also valid for the parameterisation functions $f$, $\uflt$. Then the rigorous proof of this pattern for estimates \eqref{eqn 6.5} follows from the step-by-step arguments using the identity
\begin{align*}
\dd{h_1 \cdot h_2} = \dd{h_1} \cdot \hl{2}_2 + \hl{1}_1 \cdot \dd{h_2},
\end{align*}
and the following inequality on the Sobolev norm,
\begin{align}
\Vert h_1 \cdot h_2 \Vert^{m,p} \leq  c(n,p) \Vert h_1 \Vert^{m,p} \Vert h_2 \Vert^{n-1,p},
\quad
m\leq n-1,
\quad
\left\{
\begin{aligned}
& n-1 \geq 1, p>2,
\\
& \text{ or } 
\\
& n-1\geq 2, p>1.
\end{aligned}
\right.
\label{appen pro 6.1 eqn 1}
\end{align}
The regularities of $\dd{\Xlt}$ and $\dd{\relt}$ need to be treated carefully. For $\dd{\Xlt}$, the quantity with the worst regularity is $\ddslashd \dd{\uflt}$, whose regularity is in the Sobolev space $\mathrm{W}^{n-1,p}$. For $\dd{\relt}$, the quantity with the worst regularity is $\ddcircnabla^2 \dd{\uflt}$,  whose regularity is in $\mathrm{W}^{n-2,p}$. It is crucial that in $\dd{\relt}$, there exists no product of two terms with the worst regularity, thus the inequality \eqref{appen pro 6.1 eqn 1} can always be applied when estimating $\dd{\relt}$.

To conclude the proof, we determine up to which order derivatives, the $L^{\infty}$ bounds of metric components are required in the proof. In the estimate of $\Vert \dd{\Xlt} \Vert^{n-1,p}$, the following norms of metric components show up
\begin{align*}
\Vert \vec{b}, \Omega, \slashg \Vert^{n-1,p},
\quad
\Vert \dd{\vec{b}}, \dd{\Omega}, \dd{\slashg} \Vert^{n-1,p}.
\end{align*}
In the estimate of $\Vert \dd{\relt} \Vert^{n-2,p}$, the following norms of metric components show up
\begin{align*}
\Vert \circnabla^k \partial_s^l \partial_{\us}^m a \Vert^{n-2,p},
\quad
\Vert \dd{\circnabla^k \partial_s^l \partial_{\us}^m a} \Vert^{n-2,p},
\quad
a = \vec{b}, \Omega, \slashg,
\quad
k+l+m=0,1,2.
\end{align*}
Therefore the $L^{\infty}$ bounds of metric components up to $(n+1)$-th order derivatives are required, because of the need to estimate
$\Vert \dd{\circnabla^k \partial_s^l \partial_{\us}^m \big( \vec{b}, \Omega, \slashg \big)} \Vert^{n-2,p}, 
\
k+l+m=2$.
\end{proof}

\section{Proof of proposition \ref{pro 7.1}}\label{appen pro 7.1}
\begin{proof}
We first assume that $\delta$ is suitably small such that propositions \ref{pro 5.1}, \ref{pro 5.3} on $\uls{a}{ ( \dslashd \uh )}|_{\Sigmal{a}}$, $\uls{a}{ ( \circnabla^2 \uh )}|_{\Sigmal{a}}$, and proposition \ref{pro 6.1} on $\dd{\uf}$ holds. We further assume that $\delta \leq \frac{1}{2}$, $c'_{\uh} \delta \leq 1$ throughout the proof.

Adopt the notations $\ud_m, \ud_o$, $\ubfd_m, \ubfd_o$ in formulae \eqref{eqn 3.7}, \eqref{eqn 6.4}. Furthermore, introduce the notation $\ubfd_{\uh}$ to denote
\begin{align}
\ubfd_{\uh} 
= 
c'_{\uh} \ufrakd_o
+
c'_{\uh} (|\os|/r_0+\delta_o) (\epsilon \udelta_o + \udelta_o^2) \ufrakd_m
+
c'_{\uh} ( \epsilon \udelta_m \udelta_o + \udelta_o^2) (\frakd_m + \frakd_o).
\label{appen pro 7.1 eqn 1}
\end{align}

We apply the bootstrap argument to prove the proposition. Introduce the following bootstrap assumption.
\begin{assum}
Estimates \eqref{eqn 7.3} hold for $\dd{( \dslashd \uh )|_{S_t}}$, $\dd{( \circnabla^2 \uh )|_{S_t}}$ in the interval $[0,t_a]$.
\end{assum}

We can derive the following estimates of $\dd{\Xlt_{\uh}}$, $\dd{\relt_{\uh,k}}$, $\dd{\relt_{\uh,lk}}$, from the bootstrap assumption and the estimate of $\dd{\uf}$,
\begin{align*}
&
\begin{aligned}
\leVe \dd{\Xlt_{\uh}} \riVe^{n,p} 
\leq
&
c(n,p)(|\os|/r_0 + \delta_o) ( \epsilon + c_{\uh} \udelta_o) ( \ubfd_m + \ubfd_o)
+
c(n,p)(|\os|/r_0 + \delta_o) \ubfd_{\uh}
\\
&
+
c(n,p) ( \epsilon \ud_m + \epsilon \ud_o  + c_{\uh} \udelta_o) ( \frakd_m + \frakd_o ),
\end{aligned}
\\
&
\begin{aligned}
\leVe \dd{\relt_{\uh,k}} \riVe^{n,p}
\leq
&
c(n,p) (|\os|/r_0 + \delta_o)( \epsilon c_{\uh} \udelta_o + c_{\uh}^3 \udelta_o^3) (\ubfd_m + \ubfd_o) r_0
\\
&
+
c(n,p) (|\os|/r_0 + \delta_o) \big( \epsilon (\ud_m + \ud_o) + \epsilon c_{\uh} \udelta_o + c_{\uh}^2 \udelta_o^2 \big) \ubfd_{\uh} r_0
\\
&
+
c(n,p) \big( \epsilon ( \ud_m + \ud_o) c_{\uh} \udelta_o + \epsilon c_{\uh}^2 \udelta_o^2 + c_{\uh}^3 \udelta_o^3 \big) ( \frakd_m + \frakd_o) r_0,
\end{aligned}
\\
&
\begin{aligned}
\leVe \dd{\relt_{\uh,lk}} \riVe^{n-1,p}
\leq
&
c(n,p) (|\os|/r_0 + \delta_o)( \epsilon c_{\uh} \udelta_o + c_{\uh}^2 \udelta_o^2) (\ubfd_m + \ubfd_o) r_0
\\
&
+
c(n,p) (|\os|/r_0 + \delta_o) \big( \epsilon (\ud_m + \ud_o) + c_{\uh} \udelta_o \big) \ubfd_{\uh} r_0
\\
&
+
c(n,p) \big( \epsilon ( \ud_m + \ud_o) c_{\uh} \udelta_o + c_{\uh}^2 \udelta_o^2 \big) ( \frakd_m + \frakd_o) r_0.
\end{aligned}
\end{align*}
The proof of above estimates is at the end. Assuming them, we can integrate equations \eqref{eqn 7.1} \eqref{eqn 7.2},
\begin{align*}
\big\Vert \dd{\uhlt_{\dR,k}} \big\Vert^{n,p}
\leq
&
c(n,p) \left\{
\big\Vert \dslashd \dd{\ufl{s=0}} \big\Vert^{n,p}
\right.
\\
&
\phantom{c(n,p) \big\Vert \dslashd }
\left.
+
\int_0^t 
\left[ \big\Vert \dd{\Xlt_{\uh}} \big\Vert^{n,p} \big\Vert \ddslashd \uhl{2,t}_{\dR,k} \big\Vert^{n,p} + \big\Vert \dd{\relt_{\uh,k}} \big\Vert^{n,p} \right] \d t
\right\}
\\
\leq
&
c(n,p) r_0 
\left\{
\ufrakd_o
+
(|\os|/r_0 + \delta_o) (\epsilon\udelta_o + \udelta_o^2) ( \ubfd_m + \ubfd_o)
\right.
\\
&
\phantom{\ c(n,p) r_0  \{ \ufrakd_o }
+
(|\os|/r_0 + \delta_o) (\epsilon \ud_m + \epsilon \ud_o + \udelta_o) \ubfd_{\uh}
\\
&
\phantom{\ c(n,p) r_0  \{ \ufrakd_o }
\left.
+
(\epsilon \ud_m + \epsilon \ud_o + \udelta_o) \udelta_o (\frakd_m + \frakd_o )
\right\}
\\
\leq
&
c(n,p) \ufrakd_o r_0
+ 
c(n,p) ( |\os|/r_0 + \delta_o)(\epsilon \udelta_o + \udelta_o^2) \ufrakd_m r_0
\\
&
+
(\epsilon \udelta_m \udelta_o + \udelta_o^2) (\frakd_m + \frakd_o ) r_0,
\end{align*}
and similarly
\begin{align*}
\big\Vert \dd{\uhlt_{\dR,lk}} \big\Vert^{n-1,p}
\leq
&
c(n,p) \left\{
\big\Vert \dslashd \dd{\ufl{s=0}} \big\Vert^{n-1,p}
\right.
\\
&
\phantom{c(n,p) \big\Vert \dslashd }
\left.
+
\int_0^t 
\left[ \big\Vert \dd{\Xlt_{\uh}} \big\Vert^{n,p} \big\Vert \ddslashd \uhl{2,t}_{\dR,lk} \big\Vert^{n-1,p} + \big\Vert \dd{\relt_{\uh,lk}} \big\Vert^{n-1,p} \right] \d t
\right\}
\\
\leq
&
c(n,p) \ufrakd_o r_0
+ 
c(n,p) ( |\os|/r_0 + \delta_o)(\epsilon \udelta_o + \udelta_o^2) \ufrakd_m r_0
\\
&
+
(\epsilon \udelta_m \udelta_o + \udelta_o^2) (\frakd_m + \frakd_o ) r_0.
\end{align*}
Therefore we choose $c'_{\uh} > c(n,p)$, then estimates in the bootstrap assumption at $t=t_a$ can be strengthened to strict inequalities. Thus the bootstrap argument implies that the bootstrap assumption is valid for $t \in [0,1]$. Then the proposition is proved. 

We just need to verify estimates for $\dd{\Xlt_{\uh}}$, $\dd{\relt_{\uh,k}}$ and $\dd{\relt_{\uh,lk}}$. The following heuristic rules help us obtain these estimates:
\begin{align}
\begin{aligned}
&
\begin{aligned}
&
|\os|/r_0
\rightarrow
\frakd_m,
&&
\delta_o
\rightarrow
\frakd_o,
\\
&
\udelta_m
\rightarrow
\ud_m,
&&
\udelta_o
\rightarrow
\ud_o,
\end{aligned}
\\
&
\epsilon^k r_0^l
\rightarrow
\epsilon^k r_0^l \big( \ubfd_m + \ubfd_o + \frakd_m + \frakd_o + \ubfd_{\uh} \big),
\quad
k\in \mathbb{N}, l \in \mathbb{Z},
\\
&
c_{\uh} \udelta_o, c_{\uh,2} \udelta_o
\rightarrow
\ud_{\uh},
\end{aligned}
\label{appen pro 7.1 eqn 2}
\end{align}
where the first several ones are the same as in the proof of proposition \ref{pro 6.1}, and the Leibniz rule.
\begin{center}
\begin{equation}
\begin{tikzcd}
\epsilon^k r_0^l (|\os|/r_0)^i \delta_o^j \ud_m^r \ud_o^s \cdot (c_{\uh} \udelta_o)^q
\arrow[d,"c_{\uh} \udelta_o"]
\arrow[rrrrr,"\epsilon^k r_0^l (|\os|/r_0)^i \delta_o^j \ud_m^r \ud_o^s","\text{similar as in proof of proposition \ref{pro 6.1}}"{below}]
&&&&&
( \cdots ) \cdot (c_{\uh} \udelta_o)^q
\\
\epsilon^k r_0^l (|\os|/r_0)^i \delta_o^j \ud_m^r \ud_o^s \cdot (c_{\uh} \udelta_o)^{q-1} \ud_{\uh}
&&&&&
\end{tikzcd}
\label{appen pro 7.1 eqn 3}
\end{equation}
\end{center}
Similar Leibniz rule applies to terms $\epsilon^k r_0^l (|\os|/r_0)^i \delta_o^j \ud_m^r \ud_o^s \cdot (c_{\uh,2} \udelta_o)^q$.

The rigorous proof of estimates of $\dd{\Xlt_{\uh}}$, $\dd{\relt_{\uh,k}}$, $\dd{\relt_{\uh,lk}}$ is the same as estimates of $\dd{\Xlt}$, $\dd{\relt}$ in the proof of proposition \ref{pro 6.1}. We just point out some important aspects on the regularities of $\dd{\Xlt_{\uh}}$, $\dd{\relt_{\uh,k}}$, $\dd{\relt_{\uh,lk}}$ here.

First, $\dd{\uflt}$ is in the Sobolev space $\mathrm{W}^{n,p}$ by proposition \ref{pro 6.1}. Second, $\dd{ \uhlt_{\dR,k}}$ is also in the Sobolev space $\mathrm{W}^{n,p}$ by the bootstrap assumption. Thus in $\dd{\Xlt_{\uh}}$ and $\dd{\relt_{\uh,k}}$, the terms with the worst regularities are the perturbations of background metric quantities between $\uls{a}{S_t}$,
\begin{align*}
&
\dd{b},
\quad
\dd{[R_k,b]},
\quad
\dd{\partial_{\us} b},
\\
&
\dd{\Omega^2 \slashg^{-1}},
\quad
\dd{[R_k, \Omega^2 \slashg^{-1}]},
\quad
\dd{\partial_{\us} \big( \Omega^2 \slashg^{-1} \big)},
\end{align*}
and $\dd{\uhlt_{\dR,k}}$, which are all in the Sobolev space $\mathrm{W}^{n,p}$.

In $\dd{\relt_{\dR,lk}}$, the term with the worst regularity is $\dd{\uhlt_{\dR,lk}}$, which is in the Sobolev space $\mathrm{W}^{n-1,p}$ by the bootstrap assumption.

We determine up to which order derivatives, the $L^{\infty}$ bounds of the metric components are required in the proof. 
\begin{itemize}
\item
Firstly, the estimate of $\Vert \dd{\uflt} \Vert^{n,p}$ in proposition \ref{pro 6.1} requires the $L^{\infty}$ bounds of the metric components up to $(n+1)$-th order derivatives. 
\item
Secondly, the estimates of $\big\Vert \ddslashd \uhl{2,t}_{\dR,k} \big\Vert^{n,p}$ and $\big\Vert \ddslashd \uhl{2,t}_{\dR,lk} \big\Vert^{n-1,p}$ in propositions \ref{pro 5.1}, \ref{pro 5.3} require the $L^{\infty}$ bounds of the metric components up to $(n+2)$-th order derivatives. 
\item
Thirdly, estimates of $\big\Vert \dd{\Xlt_{\uh}} \big\Vert^{n,p}$, $\big\Vert \dd{\relt_{\uh,k}} \big\Vert^{n,p}$, $\big\Vert \dd{\relt_{\uh,lk}} \big\Vert^{n-1,p}$ requires the bounds for
\begin{align*}
\big\Vert \dd{\Xlt_{\uh}} \big\Vert^{n,p}:
&
\Vert \vec{b}, \Omega, \slashg \Vert^{n,p},
\Vert \dd{\vec{b}}, \dd{\Omega}, \dd{\slashg} \Vert^{n,p},
\\
\big\Vert \dd{\relt_{\uh,k}} \big\Vert^{n,p}:
&
\left\{
\begin{aligned}
&
\Vert \circnabla^k \partial_{\us}^m \big( \vec{b}, \Omega, \slashg \big) \Vert^{n,p},
\\
&
\Vert \dd{\circnabla^k \partial_{\us}^m \big( \vec{b}, \Omega, \slashg \big)} \Vert^{n,p},
\end{aligned}
\right.
\quad
k+m=0,1,
\\
\big\Vert \dd{\relt_{\uh,lk}} \big\Vert^{n-1,p}:
&
\left\{
\begin{aligned}
&
\Vert \circnabla^k \partial_{\us}^m \big( \vec{b}, \Omega, \slashg \big) \Vert^{n-1,p},
\\
&
\Vert \dd{\circnabla^k \partial_{\us}^m \big( \vec{b}, \Omega, \slashg \big)} \Vert^{n-1,p},
\end{aligned}
\right.
\quad
k+m=0,1,2.
\end{align*}
Therefore estimates of $\big\Vert \dd{\Xlt_{\uh}} \big\Vert^{n,p}$, $\big\Vert \dd{\relt_{\uh,k}} \big\Vert^{n,p}$, $\big\Vert \dd{\relt_{\uh,lk}} \big\Vert^{n-1,p}$ require the $L^{\infty}$ bounds of the metric components up to $(n+2)$-th order derivatives, because of the estimates of
\begin{align*}
&
\Vert \dd{\circnabla^k \partial_{\us}^m \big( \vec{b}, \Omega, \slashg \big)} \Vert^{n,p},
\quad
k+m=1,
\\
&
\Vert \dd{\circnabla^k \partial_{\us}^m \big( \vec{b}, \Omega, \slashg \big)} \Vert^{n-1,p},
\quad
k+m=2.
\end{align*}
\end{itemize}
Summarising the above, we conclude that the proof requires the $L^{\infty}$ bounds of the metric components up to $(n+2)$-th order derivatives.
\end{proof}

\end{document}